\documentclass[a4paper,12pt]{article}
\usepackage{epsfig} 
\usepackage{amsfonts}
\usepackage{amssymb}
\usepackage{amsmath}
\usepackage{latexsym}
\usepackage{ae}
\numberwithin{equation}{section}

\newcommand{\N}{{\mathbb N}}

\newcommand{\us}{{\, \underline{\sim} \, }}
\newcommand{\uo}{{\, \underline{(} \, }}
\newcommand{\uc}{{\, \underline{)} \, }}
\newcommand{\uk}{{\, \underline{,} \, }}
\newcommand{\ui}{{\, \underline{\to} \,}}
\newcommand{\un}{{\, \underline{\neg} \, }}
\newcommand{\ue}{{\, \underline{\leftrightarrow}\, }}
\newcommand{\ua}{{\, \underline{\&} \, }}
\newcommand{\uv}{{\, \underline{\vee}\, }}
\newcommand{\ual}{{\, \underline{\forall} \, }}
\newcommand{\uex}{{\, \underline{\exists} \, }}
\newcommand{\dokend}{\hfill \hbox{\vrule width 5pt height 5pt depth 0pt}}
\newcommand{\ml}{[M;{\cal L}]}
\newcommand{\mlp}{[M;{\cal L}']}
\newcommand{\mlz}{\Phi_{M;{\cal L}}^Z}

\begin{document}
\thispagestyle{empty}
\author{ {\normalsize Matthias Kunik
\footnote{matthias.kunik@mathematik.uni-magdeburg.de}}\\
\small Institute for Analysis and Numerics, Otto-von-Guericke University\\
\small PSF 4120 . D-39106 Magdeburg, Germany}
\title{Formal Mathematical Systems including a \\ 
Structural Induction Principle\\
\vspace*{2cm}
\small{A revised version of the\\ 
Preprint Nr. 31/2002\\
Fakult\"at f\"ur Mathematik\\
Otto-von-Guericke-Universit\"at Magdeburg}}
\date{\today}
\maketitle

languages and formal grammars, structural induction principle,\\
G\"odel's First and Second Incompleteness Theorem.\\


\thispagestyle{empty}
\newpage
\setcounter{page}{1}
\begin{abstract}
\noindent
We present a unified theory for formal mathematical systems
including recursive systems closely related to formal grammars, including 
the predi\-cate calculus as well as a formal induction principle.
We introduce recursive systems generating
the recursively enumerable relations between lists of terms,
the basic objects under consideration. A recursive system 
consists of axioms, which are special quantifier-free positive horn formu\-las, 
and of specific rules of inference.
Its extension to formal mathe\-matical systems leads to
a formal structural induction 
with respect to the axioms of the underlying recursive system.
This approach provides some new representation theorems without using artificial
and difficult interpretation techniques. 
Within this frame we will also derive versions of G\"odel's First and Second Incompleteness 
Theorems for a general class of axiomatized formal mathe\-matical systems.
\end{abstract}
\setcounter{section}{-1}
\section{Introduction}
In this work we have developed a natural general frame for the 
formal languages usually studied in theoretical computer science  
including the predicate calculus for 
completely formalized axiomatic theories.
We present elementary proof theory for formal mathe\-matical systems
which are extensions of recursive systems generating recursively enumerable relations
between lists of terms. The recursive systems are closely related 
to formal grammars, Post's production systems and rewriting systems, 
see for example the textbooks of Hopcroft \& Ullman \cite{HU} and 
Jantzen \cite{Jn} and Post's article \cite{Ps1}.
Some advantages of our approach are:
\begin{itemize}
\item The recursive systems can be studied by its own,
independent on ques\-tions concerning mathe\-matical logic.
\item The recursive systems are directely embedded into 
formal mathe\-matical systems, i.e. the strings of the languages 
usually generated by formal grammars or Post's production systems 
are the basic objects of the first order logic. Therefore
one is neither forced to use the encoding of these languages 
into a set of G\"odel numbers nor to use interpretations 
in other formalized theories like PA or ZFC for formal languages
dealing with strings in order to study an important part of meta\-mathe\-matics.
This approach leads to a class of axiomatized mathe\-matical systems
with straightforward proofs of G\"odel's First and Second Incompleteness Theorems.
\item The most common formal systems of mathe\-matical logic are covered
by this approach, since the theory is developed
for general restrictions of the arguments in the formulas. 
\item The formal mathe\-matical systems enable a 
\textit{formal induction principle} with respect to the axioms
of the underlying recursive systems, which generalizes
the usual induction principle for integer numbers.
\end{itemize}

In Section 1 we introduce the recursive systems
which are generalizations of the so-called elementary formal systems
studied in Smullyan \cite{Sm}. The recursive systems or elementary formal systems
may be regarded as variants of Post's production systems introduced in \cite{Ps1}, 
but they are better adapted for use in mathe\-matical logic and will enable us 
to generate in a simple way the recursively enumerable relations 
between lists of terms over a finite alphabet, using the R-axioms and the
R-rules of inference introduced in Section 1. The R-axioms of the recursive system
are special quantifier-free positive horn formulas, 
which play also an important role in logic programming. 
In addition, the recursive system contains R-axioms for the use of equations. 
The R-rules of inference provide the Modus Ponens Rule and a 
simple substitution mechanism in order to obtain conclusions from the R-axioms.
Resolution strategies in order to find formal proofs
for given formulas require an own study, for details see 
Lloyd \cite{Ll}. 
We present several examples and applications for recursive systems, 
ranging from the generation of natural numbers to the simulation of 
formal grammars important in computer science and linguistics. 

In Section 2 we construct a universal recursive system 
which simulates any other recursive system. 
Then we have not only recovered the methods which were already
developed by Church, Post and Smullyan in \cite{Chr}, \cite{Ps2}, \cite{Sm}, 
but will also use these results in Section 5
to obtain new representation theorems and straightforward proofs 
of G\"odel's First and Second Incompleteness Theo\-rems for a general class 
of axio\-matized mathe\-matical systems. We will also obtain a complexity result
for a special type of recursive systems and apply it to the universal
recursive system.

In Section 3 we embed a recursive system $S$ into a formal mathe\-matical
system $M$. This embedding is consistent in the sense that the R-axioms of $S$ will become 
special axioms of $M$ and that the R-rules of inference will be special
rules of inference in $M$. The advantage of this embedding is that we can develop
considerable portions of the theory of formal mathe\-matical systems directly
in the underlying recursive systems,
without using G\"odel numbering and arithmetization. Due to the structure of
this embedding we choose the rules of inference for $M$ as a variant 
of the classical Hilbert-style instead of Gentzen-style rules for sequences 
of formulas, see \cite{Gn1} and \cite{Gn2}. 
The formulas of the mathe\-matical system are written down in
Polish prefix notation, which simplifies the formal syntax. 

The formal structural induction in the mathe\-matical systems
is performed with respect to the axioms of the underlying recursive system. 
The formal induction principle for the natural numbers is a special case, 
but it is also possible to perform the structural induction 
for an arbitrary complicate constructive structure, for example 
the induction with respect to lists, terms, formulas, and so on.

We define the formal mathe\-matical systems
with restrictions in the argument lists in the formulas. The set of
restricted argument lists contains the variables and is closed
with respect to substitutions. With these definitions we have covered
the most common formal systems.

In Section 4 we obtain several results of elementary proof theory, for example
the Deduction Theorem, the generalization of new constants in formulas
and the formal proof by contradiction. Moreover, we prove the Z-Theorem as a 
general result for the manipulation of lists of terms in the formulas 
of a formal proof. As a by-product of the Z-Theorem we can characterize 
mathe\-matical systems with certain reduced structure, for example 
formal systems which describe only relations between variables and 
constant symbols rather then relations between lists of terms.

In Section 5 we give a simple proof for the consistency of special mathematical systems
which are built up from the axioms of their underlying recursive systems.
This result is not sufficient to prove the consistency of other
mathematical systems like PA, 
but we will state an interesting conjecture, 
namely Conjecture (5.4), and will prove
that it implies the consistency of PA and of some other mathematical systems.
Conjecture (5.4) states that, under certain restrictions 
of the argument lists, variable-free prime
formulas provable in a mathematical system whose basis-axioms
coincide with the basis R-axioms of the underlying recursive system
are already provable in this recursive system. 

We close in Section 6 with an outlook concerning a possible future work 
in logic. We hope that at some point in the future this theory may lead to a cooperation 
and new applications in (computer) linguistics.\\

\section{Definition of a recursive system S}

{\bf (1.1) The symbols}\\
Given are the following pairwise disjoint sets of symbols
\begin{enumerate}
\item[(a)] 
A \underline{finite} set $A_S$ of {\bf constant symbols} or {\bf operation symbols}, 
which may be empty.
\item[(b)]
A \underline{finite} set $P_S$ of {\bf predicate symbols}, 
which may be empty.
\item[(c)]
$X := [\,{\bf x_1}\,;\,{\bf x_2}\,;\,{\bf x_3}\,;\, ... \,]$, 
a denumerable, infinite alphabet of {\bf variable symbols}.
\begin{footnote}
{If we restrict our study to recursive systems, we may also replace $X$ by a finite set.}
\end{footnote}
\item[(d)]
$E_R := [\, \sim \, ; \, ( \, ;\, ) \, ; \, , \, ; \, \to \,]$, 
five symbols representing the equivalence (or equality), 
the brackets, the comma and the implication arrow.
\end{enumerate}
We may also assume that $A_S$ and $P_S$ are finite alphabets (then their symbols are arranged in a fixed order), respectively.

{\bf (1.2) (A$_S$-)lists and sublists}
\begin{enumerate}
\item[(a)]
$a \in A_S$ and $x \in X$ are lists. 
\item[(b)]
If $\lambda$ is a list and $f \in A_S$, then $f(\lambda)$ is a list.
$\lambda$ is a sublist of $f(\lambda)$.
\item[(c)]
If $\lambda$ and $\mu$ are lists, then also $\lambda \mu$.
$\lambda$ and $\mu$ are sublists of $\lambda \mu$.
\item[(d)]
Any list $\lambda$ is sublist of itself. If $\lambda$ is sublist of $\lambda'$
and if $\lambda'$ is sublist of $\lambda''$, then $\lambda$ is sublist of $\lambda''$.
\end{enumerate}

{\bf (1.3) Constants and operation terms} (with respect to A$_S$)
\begin{enumerate}
\item[(a)]
$a \in A_S$ is a constant.
\item[(b)]
If $\lambda$ is a list and $f \in A_S$, then 
$f(\lambda)$ is an operation term.
\end{enumerate}
Constants and operation terms will be called terms.\\

{\bf (1.4)  Elementary (A$_S$-)lists and (A$_S$-)terms} \\
Let $\lambda$ be a list and $t$ be a term.
If $\lambda$, $t$ are free from variables, then they are 
called elementary list and elementary term, respectively.

\begin{figure}[h] 
\unitlength1mm \hspace*{15.4mm} 
\epsfxsize=80,7mm \put(8,3){\epsffile{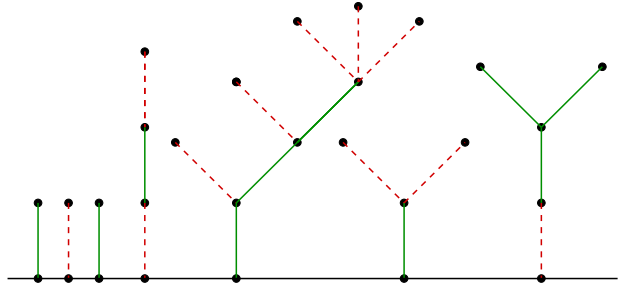}}
\caption{The elementary list $grgr(g(r))g(rg(rg(rrr)))g(rr)r(g(gg))$.}
\label{rotgruen}
\end{figure}

Figure \ref{rotgruen} illustrates the elementary list
composed on the elementary terms\\
$g$, $r$, $g$, $r(g(r))$, $g(rg(rg(rrr)))$, $g(rr)$ and $r(g(gg))$ for 
$A_S = \{\,g\,;\,r\,\}$. The solid lines are used for the symbol $g$ 
and the dashed ones for the symbol $r$.

\newpage

{\bf (1.5)  Prime R-formulas} (with respect to $A_S$ and $P_S$)
\begin{enumerate}
\item[(a)]
Let $\lambda$ and $\mu$ be lists. Then
$\sim \lambda,\mu $ is a prime R-formula, also called {\bf equation}.
$\lambda$ and $\mu$ are called argument lists of the equation.
\item[(b)]
For $p \in P_S$ and lists $\lambda_1$, $\lambda_2$, ... and so on
we define the prime R-formulas
\begin{center}
$p ~ ;\qquad p \, \lambda_1 ~;\qquad p \, \lambda_1,\lambda_2 ~;~...$\,.
\end{center}
$\lambda_1$, $\lambda_2$,... are called argument lists of these prime R-formulas.
\end{enumerate}

{\bf (1.6)  Elementary prime R-formulas} (with respect to $A_S$ and $P_S$)\\
are  prime R-formulas without variables.\\

{\bf (1.7)  R-formulas and R-subformulas} (with respect to $A_S$ and $P_S$)
\begin{enumerate}
\item[(a)]
Every prime R-formula is also an R-formula.
\item[(b)]
Let $F$ be a prime R-formula and $G$ be an R-formula. Then $\to F G$
is also an R-formula. $F$ and $G$ are R-subformulas of $\to F G$.
\item[(c)]
Every R-formula $F$ is R-subformula of itself.
If $F$ is R-subformula of $F'$ and if $F'$ is R-subformula of $F''$, 
then $F$ is R-subformula of $F''$.
\end{enumerate}
The last prime R-formula in an R-formula F is called the 
\underline{R-conclusion} of F, the other prime R-formulas in F are called the 
\underline{R-premises} of F.

{\bf (1.8)  Substitutions in R-formulas} (with respect to A$_S$ and $P_S$)\\
Let $F$ be an R-formula, $\lambda$ a list and $x \in X$. Then 
$F\,\frac{\lambda}{x}$ denotes the formula which results from $F$
by replacing everywhere in $F$ the variable $x$ by $\lambda$.
We may also write $\mbox{SbF}(F;\lambda;x)$ instead of $F\,\frac{\lambda}{x}$.
If $x \notin \mbox{var}(F)$, then $F\,\frac{\lambda}{x}=F$.\\

{\bf (1.9)  R-axioms of equality} (with respect to $A_S$ and $P_S$) 

Let $x,y \in X$ and $\lambda, \mu$ be any $A_S$-lists.
Then the following R-formulas are R-axioms of equality

\begin{tabular}{llll}
(a)  & $ \sim x,x ~ .$ &&\\
(b)  & $ \to ~\mbox{SbF}(\,\sim \lambda,\mu\,;\,x\,;\,y\,) ~            
         \to ~ \sim x,y ~
         \sim \lambda,\mu\,.$ &&\\
\end{tabular}

Let $p \in P_S$, $n \geq 1$ and $x_1,y_1,...,x_n,y_n \in X$.
Then the following R-formula is an R-axiom of equality

\begin{tabular}{llll}
(c)  & $\to ~ \sim x_1,y_1 ~...~ \to ~ \sim x_n,y_n ~ \to ~ 
p\,x_1,...,x_n~p\,y_1,...,y_n.$ &&\\
\end{tabular}\\

\underline{Remark:} Note that especially the R-axioms
$\to ~\sim x,x \to ~\sim x,y ~ \sim y,x $ and
$\to ~ \sim x,y ~ \to ~ \sim y,s ~ \sim x,s$ with $s \in X\setminus\{x\}$
result from (b).

{\bf (1.10)  A recursive system S} is given for fixed $X$ by 
$A_S$ and $P_S$ and by a finite list 
\[
B_S := [\,F_1\,;\,F_2\,; \, ... \, ;\,F_s\,]
\]
of R-formulas $F_1$, ..., $F_s$ with respect to $A_S$ and $P_S$, $s \geq 1$, 
which are called the {\bf basis R-axioms} of the recursive system S
which may be written as $S = [A_S; P_S; B_S]$.
We will in addition permit that $B_S$ may be empty.

The {\bf R-axioms} of the recursive system S are the R-axioms of equality 
and the basis R-axioms.

{\bf (1.11)  R-derivations, R-derivable R-formulas, rules of inference} 

An R-derivation in the recursive system S is a list $[F_1;...;F_l]$ of 
R-formulas $F_1,...,F_l$, including the empty list $[\,]$,
where the R-formulas $F_1,...,F_l$ are called the steps of the R-derivation, 
and is generated by the rules of inference
\begin{enumerate}
\item[(a)] \underline{Axiom Rule}:~ 
The empty list $[\,]$ is an R-derivation.
If $[\Lambda]$ is an R-derivation and $F$ an R-axiom, 
then $F$ is R-derivable and $[\Lambda;~F]$ is also an R-derivation.
\item[(b)] \underline{Modus Ponens Rule}:~
Let $[\Lambda]$ be an R-derivation, $F$, $G$ both R-formulas and 
$F$, $\to \,F \,G$ both steps of $[\Lambda]$.
Then $G$ is R-derivable and $[\Lambda;~G]$ is also an R-derivation.
\item[(c)] \underline{Substitution Rule}:~ 
Let $[\Lambda]$ be an R-derivation, F a step of $[\Lambda]$, $x$ a variable
and $\lambda$ a list.
Then $F\frac{\lambda}{x}$ is R-derivable and $[\Lambda;~F\frac{\lambda}{x}]$ 
is also an R-derivation.
\end{enumerate}
The set of all R-formulas, which are R-derivable from S, is denoted by $\Pi_R(S)$.
For $[\Lambda]=[\,]$ we put $[\Lambda;~F]=[F]$.\\

{\bf (1.12)  Recursively enumerable relations}

We fix a given recursive system $S = [A_S; P_S; B_S]$.
\begin{enumerate}
\item[(a)] Let $p \in P_S$ and $n \geq 0$ be an integer number.
With the given recursive system $S$ we define the $n-$ary recursively 
enumerable relation $R^{\,p,n}$ between \underline{elementary} $A_S-$lists 
$\lambda_1,\lambda_2,...,\lambda_n$ as follows:
\[
(\lambda_1,\lambda_2,...,\lambda_n) \in R^{\,p,n} ~ :\Leftrightarrow ~
p \, \lambda_1,\lambda_2,...,\lambda_n ~\mbox{is R-derivable in S.}
\]
The special case $\{\} \in R^{\,p,0}$ for $n=0$ means that
$p$ is R-derivable in S.
\item[(b)] The axioms of equality define an equivalence relation
$\equiv$ on the set of elementary $A_S-$lists $\lambda_1,\lambda_2$
as follows: $\lambda_1 \equiv \lambda_2$ if and only if
$\sim \lambda_1, \lambda_2$ is R-derivable in S. To the $n-$ary recursively 
enumerable relation $R^{\,p,n}$ between the elementary $A_S-$lists 
$\lambda_1,\lambda_2,...,\lambda_n$ in (a) there corresponds a relation
$R^{\,p,n}_{\,\sim}$ between the equivalence classes
$\langle \lambda_1 \rangle,
\langle \lambda_2 \rangle,...,\langle \lambda_n \rangle$ as follows:
\[
(\langle \lambda_1 \rangle,\langle \lambda_2 \rangle,...,
\langle \lambda_n \rangle) 
\in R^{\,p,n}_{\,\sim} ~
:\Leftrightarrow ~
(\lambda_1,\lambda_2,...,\lambda_n) \in R^{\,p,n}\,.
\]
The relation $R^{\,p,n}_{\,\sim}$ is also called recursively enumerable.
\end{enumerate}

\underline{Example 1:} ~ For given $A_S$, $P_S$, $X$ and
$\Box \in A_S$, $x \in X$ we consider a recursive system S 
which starts with the following two basis R-axioms:

\begin{tabular}{rcrc}
(1) & $\sim \Box\,x,x$ & \qquad \qquad \qquad \qquad 
(2) & $\sim x\,\Box,x \qquad \qquad ...$ \\
\end{tabular}\\
Here the symbol $\Box$ denotes the empty list in the formal system
and the two axioms above ensure that $\Box$ has no effect 
regarding the concatenation of lists. Therefore we can represent the empty
list in any recursive system.

\underline{Example 2:} ~ With 
$A_S := [\,a\,;\,b\,; \, f \,]$ and
$P_S := [\,W\,]$ 
we define a recursive system S by the following list 
of basis R-axioms, where $x, y \in X$ are distinct variables:

\begin{tabular}{llll}
(1)  & $W\,a$ &&\\ 
(2)  & $W\,b$ &&\\
(3)  & $\to ~ W\,x ~ \to ~ W\,y\,~ W\,xy$ &&\\
(4)  & $\sim\,f(a),a$ &&\\  
(5)  & $\sim\,f(b),b$ &&\\
(6)  & $\to ~ W\,x ~ \to ~ W\,y\,~ \sim\,f(xy),f(y)f(x)$\,. &&\\
\end{tabular}

The strings consisting on the symbols $a$ and $b$ are 
generated by the R-axioms (1)-(3). They are indicated by the
predicate symbol $W$, which is used only 1-ary here,
whereas $f$ denotes the operation which reverses the order
of such a string. For example, $\sim\,f(abaab),baaba$ is R-derivable, 
and equations like $\sim\,f(abaab),f(aab)ba$ and R-formulas
like \, $W\,f(aab)ba$ are also R-derivable.

But expressions like ~ $\to \,W x~\sim\,f(f(x)),x$ are clearly not R-derivable, 
whereas the latter R-formula will be provable in a mathematical system  
which contains $S$ as a recursive subsystem and which enables an 
induction principle with respect to the recursively enumerable relations 
represented in S. These mathematical systems will be defined in Section 3.

\underline{Example 3:} ~ With 
$A_S := [\,a\,;\,b\,]$ and $P_S := [\,W\,;\,C\,]$ 
we define a recursive system S by the following list 
of basis R-axioms, where $x, y, z \in X$ are distinct variables:

\begin{tabular}{llll}
(1)  & $W\,a$ &&\\ 
(2)  & $W\,b$ &&\\
(3)  & $\to ~ W\,x ~ \to ~ W\,y\,~ W\,xy$ &&\\
\end{tabular}\\
\begin{tabular}{llll}
(4)  & $\to ~ W\,x ~ \to ~ W\,y\,~ C\,x,y,xy$ &&\\  
(5)  & $\to ~ W\,x ~ \to ~ W\,y\,~ \to ~ W\,z\,~ C\,x,y,z,xyz$\,. &&\\  
\end{tabular}

The strings consisting on the symbols $a$ and $b$ are generated in (1)-(3) as before,
using the predicate symbol $W$, whereas in (4) and (5) we have used the
predicate symbol $C$ in order to represent the concatenation of two and three 
of these strings, respectively. 
This example demonstrates that it is possible to use the same predicate
symbol, here $C$, in order to represent different relations.

\underline{Example 4:} ~ With 
$A_S := [\,a\,]$ and
$P_S := [\,N\,;\,<\,;\,+\,;\,*\,]$ 
we define a recursive system S by the following list 
of basis R-axioms, where $x, y, z \in X$ are distinct variables:

\begin{tabular}{llll}
(1)  & $N\,a$&&\\
(2)  & $\to ~ N\,x \,~ N\,xa$&&\\
(3)  & $\to ~ N\,x ~ \to ~ N\,y\,~ <\,x,xy$&&\\
\end{tabular}
\begin{tabular}{llll}
(4)  & $\to ~ N\,x ~ \to ~ N\,y\,~ +\,x,y,xy$&&\\
(5)  & $\to ~ N\,y ~ *\,a,y,y$&&\\
(6)  & $\to ~ *\,x,y,z ~ *\,xa,y,zy$\,.&&\\
\end{tabular}

Here the positive integer numbers, indicated by the predicate symbol $N$,
are represented by $~a\,, ~aa\,, ~aaa\,, ...$ and so on. 
Let $\lambda$, $\mu$, $\nu$ be ($A_S$)-lists. Then $<\,\lambda,\mu$
is R-derivable if and only if $\lambda$ and $\mu$ represent positive integer numbers
and if the integer number represented by $\lambda$ is smaller then
the integer number represented by $\mu$.
Moreover, $+\,\lambda,\mu,\nu$ and $*\,\lambda,\mu,\nu$ are R-derivable
if and only if $\nu$ represents the positive integer number
which is the sum and the product of the two positive integer numbers represented
by $\lambda$ and $\mu$, respectively.\\

\underline{Example 5:} ~ 
With $A_S := [\, 0 \,;\, 1 \,;\, \Box \,;\, s \,;\, + \,;\, * ]$ and
$P_S := [\,N_0\,;\,NL_0^{\Box}\,]$ we define a recursive system S 
by the following list of basis R-axioms, where $x, y \in X$ 
are distinct variables:

\begin{tabular}{llll}
(1)   & $N_0\,0$ & & \\
(2)   & $\to ~ N_0\,x \,~ N_0\,s(x)$ & & \\
(3)   & $\sim 1,s(0)$ & & \\
\end{tabular}\\
\begin{tabular}{llll}
(4)   & $NL_0^{\Box}\,\Box$ & & \\
(5)   & $\to ~ N_0\,x \,~ NL_0^{\Box}\,x$ & & \\
(6)   & $\to ~ NL_0^{\Box}\,x ~ 
         \to ~ NL_0^{\Box}\,y\,~ NL_0^{\Box}\,x\,y$ & & \\
\end{tabular}\\
\begin{tabular}{llll}
(7)   & $\to ~ NL_0^{\Box}\,x ~ \sim x\,\Box,x$ & & \\
(8)   & $\to ~ NL_0^{\Box}\,x ~ \sim \Box\,x,x$ & & \\
(9)   & $\sim +(\Box),0$ & & \\
\end{tabular}\\
\begin{tabular}{llll}
(10)  & $\to ~ NL_0^{\Box}\,x ~ \sim +(0\,x),+(x)$ & & \\
(11)  & $\to ~ N_0\,x \,~ \to ~ NL_0^{\Box}\,y\, ~ 
             \sim +(s(x)\,y),s(+(x\,y))$ & & \\
\end{tabular}\\
\begin{tabular}{llll}
(12)  & $\sim *(\Box),1$ & & \\
(13)  & $\to ~ NL_0^{\Box}\,x \,~ \sim *(0\,x),0$ & & \\
(14)  & $\to ~ N_0\,x ~ \to ~ NL_0^{\Box}\,y\, ~  
             \sim *(s(x)\,y),+(*(x\,y) *(y))$\,. & & \\
\end{tabular}\\

In this example let us define the elementary terms $t_i$, $i=0,1,2,...$,
by the recursion $t_0 := 0$ and $t_{i+1}:=s(t_i)$. 
Here the non-negative integer number $i$ is represented 
by the set $K_i$ of elementary terms $t$
for which $\sim t,t_i$ is R-derivable. For example,
$\sim *(+(\Box) s(s(0)\Box) 1),0$ is R-derivable, i.e.
$*(+(\Box) s(s(0)\Box) 1) \in K_0$. A member of $K_i$ may be 
an arbitrary complicate expression, but in principle a computing
machine will be able to decide whether any given elementary term
belongs to $K_i$ or not.\\
Axioms (4)-(8) represents the lists of non-negative integer numbers
including the empty list $\Box$, which are indicated by the predicate symbol
$NL_0^{\Box}$, and ensure that the empty list has no effect on the concatenation of lists.
If $\lambda$ represents a list {\it L} of integers, then $+(\lambda)$ 
in axioms (9)-(11) represents the sum of all integer numbers in $\it{L}$, 
whereas $*(\lambda)$ in axioms (12)-(14) 
stands for the product of all integer numbers in $\it{L}$.\\

\underline{Example 6:} ~ 
Representation of a language accepted by a finite automaton 

With $A_S := [\, 0 \,;\, 1 \,]$, $P_S := [\,A\,;\,B\,;\,C\,;\,D\,]$ and $x \in X$ 
we define a recursive system $S$ by the following complete list of basis R-axioms

\begin{tabular}{llll}
(1)   & $B\,1$ & \qquad \qquad (2) ~$D\,0$ & \\
(3)   & $\to ~ A\,x \,~ D\,x0 $ & \qquad \qquad (4) ~$\to ~ A\,x \,~ B\,x1 $ & \\
(5)   & $\to ~ B\,x \,~ C\,x0 $ & \qquad \qquad (6) ~$\to ~ B\,x \,~ A\,x1 $ & \\
\end{tabular}\\
\begin{tabular}{llll}
(7)   & $\to ~ C\,x \,~ B\,x0 $ & \qquad \qquad (8) ~$\to ~ C\,x \,~ D\,x1 $ & \\
(9)   & $\to ~ D\,x \,~ A\,x0 $ & \qquad \qquad (10) ~$\to ~ D\,x \,~ C\,x1 $ & \\ 
\end{tabular}

\begin{figure}[h] 
\unitlength1mm \hspace*{25.4mm} 
\epsfxsize=60,7mm \put(8,3){\epsffile{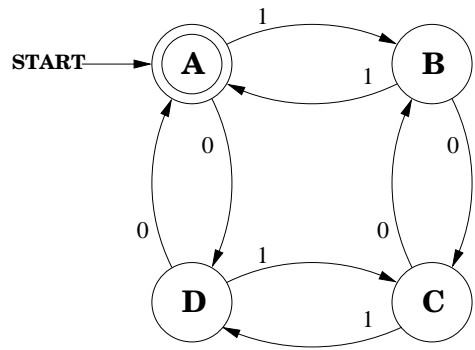}}
\caption{A finite automaton.}
\label{auto}
\end{figure}

We now consider the finite automaton with the states
$A$, $B$, $C$, $D$ depicted in Figure \ref{auto}, where $A$ is the initial
as well as the final state. A nonempty string $s=s_1...s_n$
of symbols $s_1,...,s_n \in \{\,0,\,1\}$ is called accepted 
by the finite automaton if we can follow a path of length $n$
in the graph of the automaton which starts and ends at the point $A$
and which follows a sequence of $n$ edges which are labeled by 
the symbols $s_1,...,s_n$ in the prescribed order.
The language accepted by the finite automaton consists on the set
of all strings accepted by this automaton, where we exclude for simplicity
the empty string.

For a general formal definition of a finite automaton
and the language accepted by this automaton 
see the textbook of Hopcroft and Ullman \cite{HU}.

The finite automaton accepts exactly the nonempty strings $\lambda$ over the alphabet 
$[\,0\,;\,1\,]$, for which the symbols 0 and 1 both occur an even number of times in 
$\lambda$. This set of strings is also generated in $S$ by the 1-ary predicate A. 
Here the states of the finite automaton
are the predicate symbols of the corresponding recursive system. 

The R-axioms (1)-(10) directly reflect the structure of the finite automaton. 
In the same way, any other regular language without the empty string
is accepted by a finite automaton, see \cite{HU}, and represented by a recursive
system $S$ such that all R-axioms of $B_S$ have the special form $A\,a$ or
$\to ~ B\,x \,~ C\,xb$ with $a,b \in A_S$, $A, B, C \in P_S$ and $x \in X$.\\

\underline{Example 7:} ~ 
Representation of a context-free language

For any finite alphabet or finite set of symbols $\Gamma$ let $\Gamma^*$ be the set of strings over $\Gamma$ including the empty string, whereas $\Gamma^+$
denotes the set of strings over $\Gamma$ without the empty string.

A \underline{context-free grammar} is a quadruple $G=(A,V,\Pi,v_0)$, where
\begin{itemize}
\item[(a)] $A$ is the finite set of terminal symbols\,,
\item[(b)] $V$ is the finite set of nonterminal symbols 
           with $A \cap V = \{\,\}$\,,
\item[(c)] $\Pi$ is a finite set of productions, which are strings of the form
           $v ~ \rightarrow ~ w_1\,...\,w_n$ with $v \in V$ and 
           $w_1,\,...\,w_n \in A \cup V$, $n \geq 1$. Here the symbol
           $\rightarrow$, which must not be confused with the implication arrow
           of a recursive system, neither occurs in $V$ nor in $A$.
\item[(d)] $v_0 \in V$ is a special symbol, called the starting symbol.                
\end{itemize}
The G-derivable strings $s \in (A \cup V)^+$ are defined recursively by 
\begin{itemize}
\item[(a)] $s=s_1\,...\,s_n$ is G-derivable for each production 
           $v_0 ~ \rightarrow ~ s_1\,...\,s_n \in \Pi$\,,
\item[(b)] if $s = \alpha\,v\,\beta$ with $v \in V$ and 
           $\alpha,\beta \in (A \cup V)^*$ is G-derivable and if
           $v \to w \in \Pi$, then $s' = \alpha\,w\,\beta$ is also G-derivable.
\end{itemize}
The \underline{context-free language} generated by G consists exactly 
of the G-derivable strings $s \in A^+$ without nonterminal symbols.

The standard definition also allows the derivation of empty strings, but
this would only cause technical complications in our case, whereas  
the main results about context-free grammars do not depend on this restriction.

Now we present an example. Define a context-free grammar $G$ by
$V=\{\,L\,\}$, i.e. $v_0 = L$, and $A=\{a\,;\,[\,;\,]\,;\,+\,;\,*\,\}$, 
and by the productions

\begin{tabular}{llll}
(1)   & $L ~ \to ~ a$ &&\\
(2)   & $L ~ \to ~ [L]$ &&\\
(3)   & $L ~ \to ~ L+L$ &&\\ 
(4)   & $L ~ \to ~ L*L$\,. &&\\
\end{tabular}\\

The context-free language generated by $G$ can be represented by the following
recursive system $S$: Choose $A_S=A=\{a\,;\,[\,;\,]\,;\,+\,;\,*\,\}$,
$P_S=\{L\}$, and let $x,y \in X$ be \underline{distinct} variables. 
The basis R-axioms of the recursive system $S$ are given by

\begin{tabular}{llll}
(1)   & $L~a$ &&\\ 
(2)   & $\to ~ L~x \,~ L~[x] $ && \\
\end{tabular}\\
\begin{tabular}{llll}
(3)   & $\to ~ L~x ~ \to ~ L~y ~\,L~x+y$ && \\
(4)   & $\to ~ L~x ~ \to ~ L~y ~\,L~x*y$\,. && \\
\end{tabular}\\

Here the 1-ary predicate $L$ represents the context-free language.

It is well known that every context-free language (without the empty string)
can be generated by a grammar in the normal form of Chomsky, 
where all the productions have of the special form $v \to a$ and $v_1 \to v_2\,v_3$
with $v, v_1, v_2, v_3 \in V$ and $a \in A$. 

One possible Chomsky-form of the grammar $G$ given before is\\
$G_N=(A_N,V_N,\Pi_N,L)$ with $A_N=\{a\,;\,[\,;\,]\,;\,+\,;\,*\,\}$,\\
$V_N=\{L\,;\,Bra\,;\,Ket\,;\,P\,;\,T\,;BraL\,;\,LP\,;\,LT\,\}$ 
and the productions $\Pi_N$

\begin{tabular}{llllllll}
(1)   & $L~\to~a$ & \qquad
(2) & $Bra~\to~[$ & \qquad
(3) & $Ket~\to~]$ && \\ 
\end{tabular}\\
\begin{tabular}{llllll}
(4) & $P~\to~+$ & ~~~
(5) & $T~\to~*$ && \\  
\end{tabular}\\
\begin{tabular}{llllll}
(6)   & $BraL~\to~Bra~\,L$ & \qquad
(7)   & $L~\to~BraL~\,Ket$ && \\
(8)   & $LP~\to~L~\,P$ & \qquad
(9)   & $L~\to~LP~\,L$ && \\
(10)   & $LT~\to~L~\,T$ & \qquad
(11)   & $L~\to~LT~\,L$\,. && \\
\end{tabular}\\

Then we can also replace the recursive system $S$ by another recursive system $S_N$
which is the counterpart of the grammar in Chomsky-form given before. 
In order to do this we choose the symbols and the basis R-axioms of $S_N$ as follows:

$A_{S_N}=\{a\,;\,[\,;\,]\,;\,+\,;\,*\,\}$,
$P_{S_N}=\{L\,;\,Bra\,;\,Ket\,;\,P\,;\,T\,;BraL\,;\,LP\,;\,LT\,\}$,\\ 

\begin{tabular}{lllllllllll}
(1)   & $L~a$ & \quad 
(2) & $Bra~[$ & \quad 
(3) & $Ket~]$ & \quad 
(4) & $P~+$ & \quad 
(5) & $T~*$ & \quad 
\end{tabular}\\
\begin{tabular}{llllll}
(6)   & $\to ~ Bra~x ~ \to ~ L~y \,~BraL~xy$ & \quad
(7)   & $\to ~ BraL~x ~ \to ~ Ket~y \,~L~xy$ && \\
(8)   & $\to ~L~x ~ \to ~ P~y \,~LP~xy$ &\quad
(9)   & $\to ~ LP~x ~ \to ~ L~y \,~L~xy$ && \\
(10)   & $\to ~L ~x ~ \to ~ T~y \,~LT~xy$ &\quad
(11)   & $\to ~ LT~x ~ \to ~ L~y \,~L~xy$\,. && \\
\end{tabular}\\

This example illustrates that every context-free language without the empty string
is represented by a recursive system $S$ where all the basis R-axioms have the
special form $A\,a$ and $\to ~ B\,x ~ \to ~ C\,y  ~\, D\,xy$ with $a \in A_S$
and $A, B, C, D \in P_S$ and distinct variables $x,y$. 

The restriction that $x,y \in X$ must be distinct is essential, 
which can be seen by representing the set of strings over the alphabet $[a]$ 
of length $2^n$, $n \geq 0$, with the two basis R-axioms
\[
L\,a \qquad \mbox{and} \qquad \to ~ L~x ~ \to ~ L~x \,~L~xx\,.
\] 
The language represented by the 1-ary $L$ is not context-free, as can be seen
by applying the pumping lemma for context-free languages, see Bar-Hillel,
Perles and Shamir \cite{BPS} (1961) and Wise \cite{Ws} (1976).

Note that by using a grammar or a recursive system the languages 
in our examples are generated in a quite nondeterministic way.

Finally we mention that for the context-free languages and an important
subclass, the deterministic context-free languages, one can define the so
called stack automata which are accepting these languages, 
see Chomsky \cite{Ch} (1962), Evey \cite{Ev} (1963) and \cite{HU}.\\

{\bf (1.13) Proposition, the avoidance of new symbols}~

Let $S = [A_S; P_S; B_S]$ be a recursive system and $A \supseteq A_S$
an extended set of symbols such that $S_A = [A; P_S; B_S]$ 
is also a recursive system. We suppose that $A_S$ is not empty and consider a mapping
$\gamma : A \to A_S$ with $\gamma(a)= a$ for all $a \in A_S$.
Then we can extend $\gamma$ to a function $\bar{\gamma}$, which assigns to each
R-list $\lambda$ and R-formula $F$ in $S_A$ a new R-list $\bar{\gamma}(\lambda)$
and a new R-formula $\bar{\gamma}(F)$ in S by replacing simultaneously in 
$\lambda$ and $F$ all the symbols $a \in A$ by $\gamma(a)$.

If $[\Lambda]=[F_1;...;F_l]$ is an R-derivation in $S_A$, then 
$[\Lambda]_{\bar{\gamma}}=[\bar{\gamma}(F_1);...;\bar{\gamma}(F_l)]$ 
is an R-derivation in $S$. Moreover, for all R-formulas
$F$ in $S$ there holds $F \in \Pi_R(S_A)$
if and only if $F \in \Pi_R(S)$.

\underline{Proof:} We first state the following properties of $\bar{\gamma}$,
which hold for all lists $\lambda$ and R-formulas $F$, $G$ in $S_A$ and for all $x \in X$\\
\begin{tabular}{llll}
(i)  & $\bar{\gamma}(F)=F\,,\mbox{~if~} F \in B_S\,,$ && \\
(ii)& $\bar{\gamma}(\to \,F\,G)~=~\to \,\bar{\gamma}(F)\,\bar{\gamma}(G)\,.$ && \\
(iii) & $\bar{\gamma}(F\frac{\lambda}{x})=\bar{\gamma}(F)\frac{\bar{\gamma}(\lambda)}{x}\,,$
&&\\
\end{tabular}\\
Then we employ induction with respect to the rules of inference in (1.11).
For Rule (a) we use (i),(ii),(iii), for Rule (b) we use (ii) and for Rule (c)
we use (iii). \dokend \\

{\bf (1.14) Theorem, the avoidance of equations}

Let $S = [A_S; P_S; B_S]$ be a recursive system and
$[\Lambda]$ an R-derivation in $S$. 
\begin{itemize}
\item[(a)] Suppose that the R-formulas of $B_S$ do not contain an equation
as an R-subformula. Let $[\hat{\Lambda}]$ result from $[\Lambda]$ by removing 
all the steps from $[\Lambda]$ which contain an equation as an R-subformula and 
by removing all the steps of the form $\to\,F\,F$ from $[\Lambda]$, where
$F$ is a prime R-formula. Then $[\hat{\Lambda}]$ is again an R-derivation in $S$. 
\item[(b)] Let $\sim^{*}$ be a new predicate symbol, which replaces the symbol $\sim$
and which is not present in the other set of symbols. 
Put $P_S^{*} = P_S \cup \{\sim^{*}\}$ and let $F^*$ result from
any R-formula $F$ by replacing everywhere in $F$ the symbol $\sim$ by $\sim^{*}$. 
Now we construct from $B_S$ another finite set $B_S^{*}$ of basis R-axioms 
without equations as R-subformulas as follows, where $x$, $y$, $s$, $t$
and $x_1,...,x_n, y_1,...,y_n$ are distinct variables, respectively.

\begin{tabular}{llll}
(i)   & $F^* \in B_S^*$ \, for all $F \in B_S$,&& \\
\end{tabular}\\
\begin{tabular}{llll}
(ii)  & $\sim^{*}\,x,x \, \in B_S^{*}$,&& \\
(iii) & $\to~\sim^{*}\,x,x~\to~\sim^{*}\,x,y~\sim^{*}\,y,x \, \in B_S^{*}$,&&\\
(iv)  & $\to~\sim^{*}\,x,y~\to~\sim^{*}\,y,s~\sim^{*}\,x,s \, \in B_S^{*}\,,$&&\\
\end{tabular}\\
\begin{tabular}{llll}
(v)  & $\to~\sim^{*}\,f(x),f(x) ~\to~\sim^{*}\,x,y ~\sim^{*}\,f(x),f(y)
\, \in B_S^{*}$ \,, $f \in A_S$,&&\\
(vi) & $\to~\sim^{*}\,xs,xs~\to~\sim^{*}\,s,t~\sim^{*}\,xs,xt \, \in B_S^{*}$,&&\\
(vii) & $\to~\sim^{*}\,xs,xt~\to~\sim^{*}\,x,y~\sim^{*}\,xs,yt \, \in B_S^{*}$,&&\\
\end{tabular}\\
\begin{tabular}{llll}
(viii) & $\to~\sim^{*}\,x_1,y_1~...\,\to~\sim^{*}\,x_n,y_n~
         \to\,p\,x_1,...,x_n \,p\,y_1,...,y_n \in B_S^{*}$&& \\
      & for all $p \in P_S$ and all $n \geq 1$ for which $p$ 
      occurs as a $n$-ary&&\\ 
      & prime R-subformula in $B_S$.&&\\
\end{tabular}
\end{itemize}
Let $S^* = [A_S;P_S^*;B_S^*]$ be the recursive system with the basis R-axioms
given in (i)-(viii), which do not contain any equation as an R-subformula. 
Let $n \geq 0$, $p \in P_S$ and let $\lambda$, $\mu$,
$\lambda_1$,...,$\lambda_n$ be any ($A_S$)-lists. 

Then $p\,\lambda_1,...,\lambda_n$ is R-derivable in $S$ if and only if it is
R-derivable in $S^*$, and $\sim\,\lambda,\mu$ is R-derivable in $S$ if and only if
\,$\sim^*\lambda,\mu$ is R-derivable in $S^*$.  

\underline{Remark:} $p\,\lambda_1,...,\lambda_n$ means $p$ for $n=0$.

\underline{Proof:} 
(a) Since the only R-axioms of $S$ which contain an equation as an R-subformula 
are given by (1.9), we conclude by a closer look at these R-axioms
that the only R-derivable equations must have the form 
$\sim \lambda,\lambda$. Therefore in addition to the R-formulas
containing equations the R-formulas
$\to~p\,\lambda_1,...,\lambda_n ~ p\,\lambda_1,...,\lambda_n$
occurring from (1.9)(c) after applying several times
the rules (1.11)(b,c) must be removed from an R-derivative in $S$.
A basis R-axiom of the form $\to ~ F ~F$, $F$ any prime R-formula,
is also superfluous and can be removed from an R-derivation in $S$.

(b) Let $[\Lambda^*]$ be an R-derivation of $p\,\lambda_1,...,\lambda_n$
or $\sim^*\,\lambda,\mu$ in $S^*$, respectively. Using (a) we can
suppose without loss of generality that $[\Lambda^*]$ does not contain the
symbol $\sim$. Then we can replace everywhere in $[\Lambda^*]$ the symbol 
$\sim^*$ by $\sim$ in order to obtain an R-derivation $[\Lambda]$ for 
$p\,\lambda_1,...,\lambda_n$ or $\sim\,\lambda,\mu$ in $S$, respectively.

Now let $[\Lambda]$ be any R-derivation of $p\,\lambda_1,...,\lambda_n$
or $\sim\,\lambda,\mu$ in $S$, respectively. First we cancel
all R-formulas $F$ in $[\Lambda]$ which contain any R-subformula
$q\,\lambda_1,...,\lambda_m$ with $q \in P_S$ for which $q$ 
does not occur $m$-ary as a prime R-subformula in $B_S$.
These R-formulas originating from the axioms (1.9)(c)
are clearly not prime R-formulas, so that
$p\,\lambda_1,...,\lambda_n$ and $\sim\,\lambda,\mu$ 
will not be canceled by this procedure, 
and we obtain a new R-derivation $[\hat{\Lambda}]$ in $S$.
Apart from two R-axioms of equality corresponding to (iii), (iv) we can suppose that 
the R-axiom (1.9)(b) is only used in $[\hat{\Lambda}]$ for the special cases
$\sim \lambda,\mu = \sim f(x), f(y)$, 
$\sim \lambda,\mu = xs,xt$ and $\sim \lambda,\mu = xs,yt$,
where $f \in A_S$ and $x,y,s,t \in X$ are distinct variables.
Replacing then everywhere in $[\hat{\Lambda}]$ the symbol $\sim$ by $\sim^*$
we obtain the corresponding R-derivation for
$p\,\lambda_1,...,\lambda_n$ or $\sim^*\,\lambda,\mu$ in $S^*$, respectively.
\dokend
\section{A universal recursive system}

In this section we construct a universal recursive system which simulates
any other recursive system. We prove a theorem which is due to
Smullyan \cite{Sm} and which turns out to be a version of G\"odels first 
Incompleteness Theorem. We derive a complexity result for a special type of
recursive systems and apply it to the universal recursive system.\\

{\bf (2.1) Encoding of the recursive systems}

Let $S = [A_S; P_S; B_S]$ be any recursive system. 
Here we suppose that $A_S$, $P_S$ and $X$ are \underline{lists} of symbols, 
i.e. they are ordered according to
\begin{enumerate}
\item[(a)]
$A_S = [\bf{a_1}\,;\,\bf{a_2}\,;...;\,\bf{a_k}\,]$ for the constants and operation symbols,
\item[(b)]
$P_S = [\bf{p_1}\,;\,\bf{p_2}\,;...;\,\bf{p_l}\,]$ for the predicate symbols,
\item[(c)]
$X = [\bf{x_1}\,;\,\bf{x_2}\,;\,\bf{x_3}\,;...\,]$ for the variable symbols.
\end{enumerate}

Next we define the alphabet 
$
A_{11} := [\,a\,; \, v \,;\, p \,;\, \Box \,;\, ' \,;\, * \,;
\, \us \, ; \, \uo \, ;\, \uc \, ; \, \uk \, ; \, \ui \,]
$\\
in order to encode the recursive system $S$ as follows, where $k,l$
are non-negative integer numbers which may be zero:
\begin{enumerate}
\item[(d)]
The symbols of $A_S$ are replaced by $a'\,;\,a''\,;\,a'''\,;\,...;\,a^{(k)}$.
\item[(e)]
The symbols of $P_S$ are replaced by $p'\,;\,p''\,;\,p'''\,;\,...;\,p^{(l)}$.
\item[(f)]
The variables of $X$ in $F$ are replaced by $v'\,;\,v''\,;\,v'''\,;\,...$, respectively.
\item[(g)]
The symbols of $E_R = [\, \sim \, ; \, ( \, ;\, ) \, ; \, , \, ; \, \to \,]$ 
in $F$ are replaced by \\
$ \us \, ; \, \uo \, ;\, \uc \, ; \, \uk \, ; \, \ui \,$\,, respectively.
\item[(h)]
Let $A_{11}^+$ be the set of all finite nonempty strings with 
respect to the alphabet $A_{11}$. Then to every R-formula $F$ of $S$ 
there corresponds exactly one string ${\tilde F} \in A_{11}^+$
which results from $F$ if the symbols in $F$ are replaced according to (d)-(g).
Therefore we only need the finite alphabet $A_{11}$ of symbols in order to 
encode all R-formulas of any recursive system $S$. 
\item[(i)]
We suppose that the basis R-axioms in $B_S$ are ordered according to
$B_S=[F_1\,;\,F_2\,;\,...;\,F_m]$\,, where $m$ may be zero. 
We encode the recursive system $S$ by defining the corresponding 
\underline{R-basis string} ${\tilde S}$ according to 
\[
{\tilde S} = u*w*{\tilde F}_1*{\tilde F}_2*...*{\tilde F}_m*
\]
If $m=0$, then 
$
{\tilde S} = u*w*\Box*\,.
$
Here $u,w \in \{\,\Box \,;\, ' \,;\, '' \,;... \,\}$ are strings which recover 
the finite alphabets $A_S$ and $P_S$. If $k=0$, i.e. $A_S$ is empty,
then $u = \Box$, otherwise $u$ consists on a string of $k=|A_S|$ accents.
If $l=0$, i.e. $P_S$ is empty, then $w = \Box$, otherwise $w$ consists 
on a string of $l=|P_S|$ accents. Note that the knowledge of ${\tilde S} \in A_{11}^+$ 
allows a complete reconstruction of the original recursive system $S$.
\end{enumerate}

\underline{Example 1:}~We define the recursive system $S$ by
$A_S = [\,a,b\,]$, $P_S = [\,p,q\,]$ and the three basis R-axioms
for distinct variables $x$, $y$

(1) $p\,a,ab$ \qquad 
(2) $\to~ p\,x,y ~ p\,xa,yab$ \qquad 
(3) $\to~ p\,x,y ~ q\,y$\,.

If we put $x = {\bf x_1}$ and $y = {\bf x_2}$, then 
the encoding of $S$ gives the R-basis string

${\tilde S} = 
''*''*p'a'\uk a'a''*\ui p'v'\uk v''p'v'a'\uk v''a'a''*\ui p'v'\uk v''p''v''*$\,.\\

{\bf (2.2)  The universal recursive system $S_{11}$}

The constants and operation symbols of $S_{11}$ are given by the 
alphabet $A_{11}$. The symbols $x,y,u,w,z,r,t,s$ denote distinct variables
and the predicate symbols of $S_{11}$ are included in the list of basis R-axioms 
of $S_{11}$ given by

\begin{tabular}{llll}
(1a)~   & $Acc~'$ & & \\
(1b)~   & $\to ~ Acc\,x \,~ Acc\,x'$ & & \\
(2a)~   & $N_0\,\Box$ & & \\
(2b)~   & $\to ~ Acc\,x \,~ N_0\,x$ & & \\
\end{tabular}

$Acc\,x$ means that $x$ is a nonempty string consisting only of accents,
whereas $N_0\,x$ means that $x \in \{\Box\,;\,'\,;\,''\,;... \}$ represents 
a non-negative integer number.

\begin{tabular}{llll}
(3)~  & $\to ~ Acc\,x ~ \to ~ Acc\,y\,~ <\,x,xy$&&\\
(4a)~  & $\to ~ Acc\,x ~ \leq\,x,x$&&\\
(4b)~  & $\to ~ <\,x,y ~ \leq\,x,y$&&\\
\end{tabular}\\
\begin{tabular}{llll}
(5)~~  & $\to ~ \leq\,x,u ~ A_s\,ax,u$&&\\
(6)~~  & $\to ~ \leq\,x,w ~ P_s\,px,w$&&\\
(7)~~  & $\to ~ Acc\,x ~V\,vx$&&\\
\end{tabular}

From now on $u$ and $w$ represent the non-negative integer numbers
$|A_S| \geq 0$ and $|P_S| \geq 0$, respectively. 
$A_s\,ax,u$ means that $ax$ represents a constant symbol in $A_S$ 
and $P_s\,px,w$ that $px$ represents a predicate symbol in $P_S$\,.  
$V\,vx$ means that $vx$ represents the variable symbol $\bf{x_i}$, where 
$x$ consists on $i \geq 1$ accents.

\begin{tabular}{llll}
(8a)~  & $\to ~ A_s\,x,u ~ L\,x,u $&&\\
(8b)~  & $\to ~ V\,x     ~ \to ~ N_0\,u     ~ L\,x,u $&&\\
(8c)~  & $\to ~ A_s\,x,u ~ \to ~ L\,y,u ~ L\,x \uo y \uc,u $&&\\
(8d)~  & $\to ~ L\,x,u ~ \to ~ L\,y,u ~ L\,xy,u $&&\\
\end{tabular}

$L\,x,u$ means that $x$ represents a list (with respect to $A_S$).

\begin{tabular}{llll}
(9a)~  & $\to ~ A_s\,x,u ~ EL\,x,u $&&\\
(9b)~  & $\to ~ A_s\,x,u ~ \to ~ EL\,y,u ~ EL\,x \uo y \uc,u $&&\\
(9c)~  & $\to ~ EL\,x,u ~ \to ~ EL\,y,u ~ EL\,xy,u $&&\\
\end{tabular}

$EL\,x,u$ means that $x$ represents an elementary list (with respect to $A_S$).

\begin{tabular}{llll}
(10a)  & $\to ~ L\,x,u ~ LL\,x,u$&&\\
(10b)  & $\to ~ LL\,x,u ~ \to ~ L\,y,u ~ LL\,x \uk y,u$&&\\
\end{tabular}

$LL\,x,u$ means that $x$ represents a finite sequence of lists which are separated
by the underlined comma.

\begin{tabular}{llll}
(11)   & $\to ~ N_0\,w ~ \to ~ L\,x,u ~ \to ~ L\,y,u ~ Eq\,\us x \uk y,u,w $&&\\
(12a)  & $\to ~ Eq\,x,u,w ~ PRF\,x,u,w $&&\\
(12b)  & $\to ~ P_s\,x,w ~ \to ~ N_0\,u ~ PRF\,x,u,w $&&\\
(12c)  & $\to ~ P_s\,x,w ~ \to ~ LL\,y,u ~ PRF\,xy,u,w $&&\\
\end{tabular}

$Eq\,x,u,w$ and $PRF\,x,u,w$ means that $x$ represents an equation and a
prime R-formula, respectively.

\begin{tabular}{llll}
(13a)  & $\to ~ EL\,x,u ~ ELL\,x,u$&&\\
(13b)  & $\to ~ ELL\,x,u ~ \to ~ EL\,y,u ~ ELL\,x \uk y,u$&&\\
\end{tabular}

$ELL\,x,u$ means that $x$ represents a finite sequence of elementary
lists which are separated by the underlined comma.

\begin{tabular}{llll}
(14a)  & $\to ~ N_0\,w ~ \to ~ EL\,x,u ~ \to ~ EL\,y,u ~ EPRF\,\us x \uk y,u,w $&&\\
(14b)  & $\to ~ P_s\,x,w ~ \to ~ N_0\,u ~ EPRF\,x,u,w $&&\\
(14c)  & $\to ~ P_s\,x,w ~ \to ~ ELL\,y,u ~ EPRF\,xy,u,w $&&\\
\end{tabular}

$EPRF\,x,u,w$ means that $x$ represents an elementary prime R-formula.

\begin{tabular}{llll}
(15a)  & $\to ~ PRF\,x,u,w ~ RF\,x,u,w$&&\\
(15b)  & $\to ~ PRF\,x,u,w ~ \to ~ RF\,y,u,w ~ RF\,\ui xy,u,w $&&\\
\end{tabular}

$RF\,x,u,w$ means that $x$ represents an R-formula.

\begin{tabular}{llll}
(16a)  & $\to ~ <\,x,y ~ VV\,vx,vy $&&\\
(16b)  & $\to ~ <\,x,y ~ VV\,vy,vx $&&\\
\end{tabular}

$VV\,x,y$ means that $x$ and $y$ represent two different variables. 

\begin{tabular}{llll}
(17a)  & $\to ~ A_s\,x,u ~ \to ~ V\,z ~ \to L\,r,u ~ SbL\,x,r,z,x,u$&&\\
(17b)  & $\to ~ V\,x ~ \to L\,r,u ~ SbL\,x,r,x,r,u$&&\\
(17c)  & $\to ~ VV\,x,z ~ \to L\,r,u ~ SbL\,x,r,z,x,u$&&\\
(17d)  & $\to ~ A_s\,x,u ~\to ~ SbL\,y,r,z,t,u~ SbL\,x \uo y \uc ,r,z,x \uo t \uc,u$&&\\
(17e)  & $\to ~ SbL\,x,r,z,s,u ~\to ~ SbL\,y,r,z,t,u~ 
               SbL\,xy,r,z,st,u$&&\\
\end{tabular}

$SbL\,x,r,z,s,u$ means that $s$ represents the list which results from the list
represented by $x$ after the substitution of the variable represented by $z$
by the list represented by $r$.

\begin{tabular}{llll}
(18a)  & $\to ~ SbL\,x,r,z,s,u ~ SbLL\,x,r,z,s,u$&&\\
(18b)  & $\to ~ SbLL\,x,r,z,s,u ~ \to ~ SbL\,y,r,z,t,u ~ 
                SbLL\,x \uk y,r,z,s \uk t,u$&&\\
\end{tabular}

$SbLL\,x,r,z,s,u$ is the generalization of $SbL\,x,r,z,s,u$ for finite sequences 
of lists separated by the underlined comma, which are represented here by $x$ and $s$, whereas
$r$ represents a list as before.

\begin{tabular}{llll}                
(19a)  & $\to ~ N_0\,w ~ \to ~ SbL\,x,r,z,s,u ~ \to ~ SbL\,y,r,z,t,u$&&\\ 
       & $ \quad ~~ SbPRF\,\us x \uk y,r,z,\us s \uk t,u,w$&&\\
(19b)  & $\to ~ P_s\,x,w ~ \to ~ V\,z ~ \to ~ L\,r,u ~ SbPRF\,x,r,z,x,u,w$&&\\
(19c)  & $\to ~ P_s\,x,w ~ \to ~ SbLL\,y,r,z,t,u ~ SbPRF\,xy,r,z,xt,u,w$&&\\
\end{tabular}

$SbPRF\,x,r,z,s,u,w$ means that $s$ represents the prime R-formula which results 
from the prime R-formula represented by $x$ after the substitution of the variable 
represented by $z$ by the list represented by $r$.

\begin{tabular}{llll}
(20a)  & $\to ~ SbPRF\,x,r,z,s,u,w ~ SbRF\,x,r,z,s,u,w $&&\\
(20b)  & $\to ~ SbPRF\,x,r,z,s,u,w ~ \to ~ SbRF\,y,r,z,t,u,w$&&\\
       & $\quad ~~ SbRF\,\ui xy,r,z,\ui st,u,w$&&\\
\end{tabular}

$SbRF\,x,r,z,s,u,w$ means that $s$ represents the R-formula which results 
from the R-formula represented by $x$ after the substitution of the variable 
represented by $z$ by the list represented by $r$.

\begin{tabular}{llll}
(21a)  & $\to ~ RF\,x,u,w ~ SbRF\,x,x,u,w $&&\\
(21b)  & $\to ~ SbRF\,xz,s,z,rs,u,w ~ SbRF\,xz,rs,u,w $&&\\
(21c)  & $\to ~ SbRF\,xzy,s,z,rst,u,w ~ SbRF\,xzy,rst,u,w $&&\\
\end{tabular}

$SbRF\,x,s,u,w$ means that there is a variable represented by $z$ and 
a list represented by $r$ such that $SbRF\,x,r,z,s,u,w$ is R-derivable.

\begin{tabular}{llll}
(22a)  & $\to ~ V\,x ~ \to ~ V\,y ~ AP\,\ui \us x \uk y,x,y$&&\\
(22b)  & $\to ~ V\,x ~ \to ~ V\,y ~ \to ~ AP\,r,s,t ~\, 
          AP\,\ui \us x \uk y r,x \uk s,y \uk t$&&\\
\end{tabular}

$AP$ is an auxiliary predicate needed for the
representation of the equality axioms of the form (1.9)(c).
 
\begin{tabular}{llll}
(23a)  & $\to ~  N_0\,u ~ \to ~  N_0\,w ~ \to ~ V\,x ~\, 
         EqA \, \us x \uk x,u,w$&&\\
(23b)  & $ \to ~ V\,x ~ \to ~ V\,y ~ \to ~ Eq\,z,u,w
           ~ \to ~ SbPRF\,z,x,y,s,u,w$&&\\
& $\quad~\, EqA \, \ui \,s ~ \ui \, \us x \uk y \,z,u,w$&&\\

(23c)   & $\to ~ AP\,r,s,t ~ \to ~ P_s\,z,w ~ \to ~ N_0\,u ~\,
         EqA \,r \ui zs zt,u,w$&&\\
\end{tabular}

$EqA~x,u,w$ means that $x$ represents an axiom of equality.

\begin{tabular}{llll}
(24a)  & $\to ~ RF\,x,u,w ~\, RBasis^+\,u*w*x*$&&\\
(24b)  & $\to ~ RF\,x,u,w ~ \to ~  RBasis^+\,u*w*s* ~\, RBasis^+\,u*w*s*x*$&&\\
\end{tabular}

$RBasis^+\,x$ means that $x$ is an R-basis string with $|B_S| \geq 1$.

\begin{tabular}{llll}
(25a)  & $\to ~  N_0\,u ~ \to ~  N_0\,w ~ RBasis\,u*w*\Box*$&&\\
(25b)  & $\to ~ RBasis^+\,x ~\, RBasis\,x$&&\\
\end{tabular}

$RBasis\,x$ means that $x$ is an R-basis string including $|B_S| = 0$.

\begin{tabular}{llll}
(26a)  & $\to ~ RF\,x,u,w ~\, BRA\,u*w*x*,x$&&\\
(26b)  & $\to ~ RF\,x,u,w ~ \to ~ RBasis^+\,u*w*s* ~\, BRA\,u*w*s*x*,x$&&\\
(26c)  & $\to ~ RF\,x,u,w ~ \to ~ BRA\,u*w*s*,y ~\, BRA\,u*w*s*x*,y$&&\\
\end{tabular}

$BRA\,x,y$ means that $x$ is an R-basis string and that $y$ represents a
basis R-axiom of the recursive system determined by $x$. Then $|B_S| \geq 1$.

\begin{tabular}{llll}
(27a)  & $\to ~ EqA \,x,u,w ~ \to ~ RBasis\,u*w*s* ~\, RA\,u*w*s*,x$&&\\
(27b)  & $\to ~ BRA\,x,y ~\, RA\,x,y$&&\\
\end{tabular}

$RA\,x,y$ means that $x$ is an R-basis string and that $y$ represents an
R-axiom of the recursive system determined by $x$. 

\begin{tabular}{llll}
(28)  & $\to ~ PRF\,x,u,w ~ \to ~ BRA\,u*w*s*,x ~\, PBRA\,u*w*s*,x$&&\\
\end{tabular}

$PBRA\,x,y$ means that $x$ is an R-basis string and that $y$ represents 
a prime basis R-axiom of the recursive system determined by $x$.

\begin{tabular}{llll}
(29a)  & $\to ~  N_0\,u ~ \to ~  N_0\,w ~ \to ~ RA\,u*w*s*,x ~\,
       D_s^+\, u*w*x*, u*w*s*$&&\\
(29b)  & $\to ~  D_s^+\, x,y ~ \to ~ RA\,y,z ~\, D_s^+\, xz*,y$&&\\
\end{tabular}
\begin{tabular}{llll}
(29c)  & $\to ~  D_s^+\, x\ui rsz,t ~ \to ~ BRA\,x\ui rsz,\ui rs~ \to ~ PBRA\,x\ui rsz,r$&&\\
       & $\to ~ RBasis^+\, x\ui rszs* ~ D_s^+\, x\ui rszs*,t$&&\\
(29d)  & $\to ~  D_s^+\, u*wxyz, t ~ \to ~ BRA\,u*wxyz,y~ \to ~ SbRF\,y,s,u,w$&&\\
       & $\quad~\, D_s^+\, u*wxyzs*, t$&&\\
\end{tabular}

$D_s^+\,x,y$ means that $x$ represents a nonempty R-derivation in the recursive system
given by the R-basis string $y$. The premise $RBasis^+\, x\ui rszs*$ in (29c) guarantees that
$s$ represents an R-formula.

\begin{tabular}{llll}
(30a)  & $\to ~  N_0\,u ~ \to ~  N_0\,w ~ 
          \to ~ RBasis\,u*w*s* ~\, D_s\, u*w*\Box*, u*w*s*$&&\\
(30b)  & $\to ~ D_s^+\, x,y  ~\, D_s\, x, y$&&\\
\end{tabular}

$D_s\,x,y$ means that $x$ represents an R-derivation (which may be empty) 
in the recursive system given by the R-basis string $y$.

\begin{tabular}{llll}
(31)  & $\to ~ EPRF\, x,u,w  ~ \to ~ BRA~t,x  ~ \to ~ D_s\,t,y ~\, \Omega_s\,yx$&&\\
\end{tabular}

In this context $\Omega_s\,yx$ means that $x$ represents an elementary prime R-formula
which is R-derivable in the recursive system given by the R-basis string $y$.\\

{\bf (2.3)  Definition of $S_{11}$-statements and $S_{11}$-theorems}

$z=yx$ with $x, y \in A_{11}^+$ is called 
\underline{$S_{11}$-statement} if and only if
$y$ is an R-basis string which represents a recursive system $S$
and $x$ represents an elementary prime R-formula, not necessary in $S$. 
If in addition $\Omega_s\,yx$ is R-derivable in $S_{11}$, 
then $z$ is called \underline{$S_{11}$-theorem}.

Note that $z$ is not an R-basis string since the last symbol in $z$ 
is not the ``*". The $S_{11}$-statement $z=yx$ is called $n$-ary, $n \geq 0$, 
if the elementary prime R-formula represented by $x$ is $n$-ary.\\

{\bf (2.4)  Definition of $S_{11}$-predicates}

If $y$ is an R-basis string which represents a recursive system $S$
and if $q$ represents a predicate symbol, not necessary in $S$, 
then $P=yq$ is called  \underline{$S_{11}$-predicate}.
If \, $ELL\,s,u$ \, is R-derivable in $S_{11}$ for some $s,u \in A_{11}^+$, 
then it is easy to check that $Ps=yqs$ is an $S_{11}$-statement. 
We say that $s$ \underline{satisfies} 
the $S_{11}$-predicate $P$ if in addition $\Omega_s\,Ps$ is R-derivable in $S_{11}$. 
In this case $Ps=yqs$ is an $S_{11}$-theorem.\\

\underline{Example 2:}~The $A_{11}$-string 

$P = \,''*'*p'a' \uk a' \uo a'' \uc *p'a'a''* \ui p'v' p'v'a'a''*p'$ 

is an $S_{11}$-predicate
which is satisfied by the elementary lists 

$a'a''$,~ $a'a''a'a''$,~ $a'a''a'a''a'a''$,~ ... and so on, 

and therefore we obtain the following 1-ary $S_{11}$-theorems:

$P\,a'a''$,~ $P\,a'a''a'a''$,~ $P\,a'a''a'a''a'a''$,~ ...\,.

Moreover the string $s = a' \uk a' \uo a'' \uc$ satisfies the $S_{11}$-predicate $P$
and gives the 2-ary $S_{11}$-theorem

$''*'*p'a' \uk a' \uo a'' \uc *p'a'a''* \ui p'v' p'v'a'a''*p'a' \uk a' \uo a'' \uc$ .

On the other hand, for $s = a'''$ the $A_{11}$-string $Ps$ is an 
$S_{11}$-statement, but not an $S_{11}$-theorem.\\

{\bf (2.5)  The diagonalization of $S_{11}$-predicates}

There is a very simple method in order to generate a so called 
\underline{self-referential} $S_{11}$-statement. We first define the mapping 
$g_{11}:A_{11}^+ \to A_{11}^+$ by

\begin{tabular}{llll}
$g_{11}(a) = a'~$   & $g_{11}(v) = a''~$ & $g_{11}(p) = a'''~$ \\
$g_{11}(\Box) = a''''~$ & $g_{11}(') = a'''''~$   & $g_{11}(*) = a''''''~$\\
$g_{11}(\us) = a'''''''~$ & $g_{11}(\uo) = a''''''''~$ &
$g_{11}(\uc) = a'''''''''~$\\   
$g_{11}(\uk) = a''''''''''~$ & $g_{11}(\ui) = a'''''''''''~$ & 
$g_{11}(xy) = g_{11}(x)g_{11}(y)~$. & \\
\end{tabular}

Then the \underline{diagonalization} of
any $S_{11}$-predicate $P$ is given by 

$\mbox{Diag}(P) = P\,g_{11}(P)$.
Note that $\mbox{Diag}(P)$ is an $S_{11}$-statement.\\ 

\underline{Example 3:}~The $A_{11}$-string 
$P = \,'''''''''''*'*p'v'*p'$ 
is an $S_{11}$-predicate which represents any list over an alphabet
consisting on 11 symbols,
and therefore we conclude that its diagonalization is an $S_{11}$-theorem:

\begin{tabular}{llll}
$\mbox{Diag}(P) =$ &$\, '''''''''''*'*p'v'*p'a'''''a'''''a'''''a'''''
a'''''a'''''a'''''a'''''a'''''a'''''a'''''$&&\\
&$a''''''a'''''a''''''a'''a'''''a''a'''''a''''''a'''a'''''$\,.&&\\
\end{tabular}\\

{\bf (2.6)  A version of G\"odel's First Incompleteness Theorem}

Let $B_s^{(1)}$ be the set of all 1-ary $S_{11}$-statements and 
$\Omega_s^{(1)}$ the set of all 1-ary $S_{11}$-theorems. Then
$B_s^{(1)}$ and $\Omega_s^{(1)}$ are recursively enumerable, but not
the complement ${\overline \Omega}_s^{(1)}=B_s^{(1)} \setminus \Omega_s^{(1)}$.

\underline{Proof:}

In order to see that $B_s^{(1)}$ and $\Omega_s^{(1)}$ are recursively enumerable,
we extend $S_{11}$ by the predicate symbols ``$B_s^{(1)}$" and ``$\Omega_s^{(1)}$"
and add the two basis R-axioms

\begin{tabular}{llll}
 & $\to ~ RBasis~x  ~ \to ~ P_s\,y,w  ~ \to ~ EL\,z,u ~\, B_s^{(1)}\,xyz$&&\\
 & $\to ~ B_s^{(1)}\,x ~ \to ~ \Omega_s\, x ~\, \Omega_s^{(1)}\, x$&&\\
\end{tabular}

with distinct variables $x,y,z,u,w \in X$.
To use the same notation for the sets and the corresponding predicate symbols 
will not lead to confusions.

We assume that the set ${\overline \Omega}_s^{(1)}$ is recursively enumerable.
Then due to Theorem (1.14) there is a recursive system $S=[A_S;P_S;B_S]$ 
which represents ${\overline \Omega}_s^{(1)}$ such that no equation is 
involved in $B_S$. Let the members of $A_S$, $P_S$ and $B_S$ be given
in a fixed order. Since $A_S$ must contain the symbols of $A_{11}$
due to our assumption that ${\overline \Omega}_s^{(1)}$ is represented in $S$,
we can suppose without loss of generality that $A_S$ starts with the alphabet $A_{11}$
in the prescribed order given for $A_{11}$. 
We can suppose that the predicates $RBasis$ and $P_s$ from $S_{11}$ are also represented
in $S$, using the predicate symbols $``RBasis"$ and $``P_s"$.
The predicate symbol representing ${\overline \Omega}_s^{(1)}$ in $S$ 
may also be denoted by ``${\overline \Omega}_s^{(1)}$". Moreover we extend $P_S$
by the two \underline{new} symbols ``$G_{11}$" and ``${\overline \Omega}_s^{(1)\#}$"
to a new alphabet $P_S^{\#}$ and we extend $B_S$ to a new list $B_S^{\#}$
of basis R-axioms by adding the following R-axioms to $B_S$ for distinct
$x,y,r,s,w \in X$

\begin{tabular}{llll}
 (1) \quad $G_{11}\,a,a'~$   
&(2) \quad $G_{11}\,v,a''~$ 
&(3) \quad $G_{11}\,p,a'''~$\\ 
 (4) \quad $G_{11}\,\Box,a''''~$
&(5) \quad $G_{11}\,',a'''''~$   
&(6) \quad $G_{11}\,*,a''''''~$\\
 (7) \quad $G_{11}\,\us,a'''''''~$ 
&(8) \quad $G_{11}\,\uo,a''''''''~$ 
&(9) \quad $G_{11}\,\uc,a'''''''''~$\\  
(10) \quad $G_{11}\,\uk,a''''''''''~$ 
&(11) \quad $G_{11}\,\ui,a'''''''''''~$ &\\ 
\end{tabular}\\ 
\begin{tabular}{llll}
(12) \, $\to \, G_{11}\,x,r \, \to \, G_{11}\,y,s ~ G_{11}\,xy,rs$&& \\

(13) \, $\to \, RBasis \, x \to \, P_s \, y,w \to \, G_{11}\,x,r \, 
\, \to G_{11}\,y,s \, \to \, {\overline \Omega}_s^{(1)}\,xyrs ~ 
{\overline \Omega}_s^{(1)\#}\,xy$&& \\
\end{tabular}

There results an extended recursive system $S^{\#}=[A_S;P_S^{\#};B_S^{\#}]$.

The relation $G_{11}\,\lambda,\mu$ generated by the R-axioms (1)-(12)
is satisfied if and only if there hold $\lambda,\mu \in A_{11}^+$
and $\mu = g_{11}(\lambda)$. Moreover, ${\overline \Omega}_s^{(1)\#}\,\lambda$ 
is R-derivable in $S^{\#}$ if and only if \, $\lambda$ is an $S_{11}$-predicate
and $\mbox{Diag}(\lambda)=\,\lambda g_{11}(\lambda) \in {\overline \Omega}_s^{(1)}$.
We write $\lambda \in {\overline \Omega}_s^{(1)\#}$ in order to express that 
${\overline \Omega}_s^{(1)\#}\,\lambda$ is R-derivable in $S^{\#}$.

These representation properties are guaranteed since the equations
are excluded from $B_S^{\#}$ and since the symbols 
$G_{11},{\overline \Omega}_s^{(1)\#} \in P_S^{\#}$ are not in $P_S$.

By forming the R-basis string for the recursive system $S^{\#}$
we can construct the $S_{11}$-predicate $P$ corresponding to the set
${\overline \Omega}_s^{(1)\#}$ represented in $S^{\#}$.
Since the alphabet $A_S$ of $S^{\#}$ starts with the alphabet $A_{11}$,
we obtain for all $\lambda \in A_{11}^+$
\[
\lambda \in {\overline \Omega}_s^{(1)\#} \quad  \Leftrightarrow \quad 
P\,g_{11}(\lambda) \in \Omega_s^{(1)}\,.
\]
If we put $\lambda = P$, then 
\[
P \in {\overline \Omega}_s^{(1)\#} \quad  \Leftrightarrow \quad 
P\,g_{11}(P) \in \Omega_s^{(1)}\,.
\]
This equivalence contradicts the construction of 
the set ${\overline \Omega}_s^{(1)\#}$, which requires that 
the $S_{11}$-predicate $P$ should satisfy the equivalence
\[P \in {\overline \Omega}_s^{(1)\#} \quad  \Leftrightarrow \quad 
P\,g_{11}(P) \in {\overline \Omega}_s^{(1)}\,.\]
Thus we have proven Theorem (2.6). \dokend \\

In Section 5, Theorem (5.6) we will explain
in what sense this result may be regarded as a version of 
G\"odels First Incompleteness Theorem.\\
\vspace*{1cm}\\

\underline{Remarks:}
\begin{itemize}
\item[(i)] The recursive systems considered in Smullyan \cite{Sm}
are called elementary formal systems there. Like the recursive system $S_{11}$,
they do not contain the equations and the operation terms, 
but this is of course not a principle restriction 
for the construction of recursively enumerable relations.
\item[(ii)] The construction of $S_{11}$ was only needed in order to prove that
${\Omega}_s^{(1)}$ and $B_s^{(1)}$ are recursively enumerable. In order to
prove that $\overline{\Omega}_s^{(1)}$ is not recursively enumerable we can
directely define all the necessary ingredients like $S_{11}$-statements, 
$S_{11}$-theorems and $S_{11}$-predicates by using the encoding (2.1) for
the recursive systems.
\item[(iii)] Due to Church's thesis and Theorem (2.6) we conclude that
there is no algorithm which enables us to decide whether a given R-formula
of the recursive system $S_{11}$ is R-derivable or not. The reason for this
is the fact that the 1-ary predicate ${\Omega}_s^{(1)}$ is not decidable.
But the other predicates of $S_{11}$ generated by (2.2) (1a)-(30b) are decidable,
since they form a recursive subsystem which satisfies the following
\end{itemize}

{\bf (2.7) Definition of special recursive systems and predicates}

We consider a recursive system $S=[A_S; P_S; B_S]$.
Then $S$ and the predicates represented in $S$ are called 
\underline{special recursive} if
\begin{itemize}
\item there is no equation involved in $B_S$,
\item every argument list occurring in the R-premises of any R-axiom $F$ 
also occurs as a sublist in an argument list of the R-conclusion of $F$.
\end{itemize}

In order to estimate the complexity of an algorithm looking for
an R-derivation of an elementary prime R-formula $p\,\lambda_1,...,\lambda_i$
in a given special recursive system $S$ we need two Lemmata.
We shall prove that resolution strategies for special recursive predi\-cates
will only require polynomial effort with respect to the length
of the ``input formula" $p\,\lambda_1,...,\lambda_i$.
As a consequence, special recursive predicates are decidable.\\

{\bf (2.8) Lemma} 

Let $\lambda$ be any $A_S$-list consisting on $|\lambda| = n$ symbols.
Then the number of sublists in $\lambda$ is less or equal to
$
\frac{n(n+1)}{2}\,.
$

\underline{Proof}: Induction with respect to $n$.\\

{\bf (2.9) Lemma}

Let $\mu$ be any $A_S$-list and let $\lambda$ be any elementary 
$A_S$-list consisting on $|\lambda|= n$ not necessary 
distinct symbols. Let $x_1,...,x_k$ with $1 \leq k \leq n$ be
the list of distinct variables occurring in $\mu$, 
ordered according to their first appearance. 
By $\mbox{Inst}(\mu,\lambda)$ we denote the set of
all mappings which assign to each variable $x_j$ in $\mu$
an elementary $A_S$-list $\kappa_j$ such that
$
\lambda = \mu\,\frac{\kappa_1}{x_1}...\frac{\kappa_k}{x_k}
$.
Then 
$$
|\mbox{Inst}(\mu,\lambda)| \leq \begin{pmatrix}
n-1\\k-1
\end{pmatrix}\,.
$$

\underline{Proof}: Induction with respect to $k$.\\

For $1 \leq k \leq n$ we put $\Gamma(n,k)=\max \limits_{j=1}^{k}
\begin{pmatrix}
n-1\\j-1
\end{pmatrix}$. We define  $\Gamma(n,0)= 0$ and $\Gamma(n,k)=\Gamma(n,n)$
for $k > n$.

{\bf (2.10) Theorem}

Let $S=[A_S; P_S; B_S]$ be a special recursive system and let
$p\,\lambda_1,...,\lambda_i$ be an elementary prime R-formula
which is R-derivable in $S$. 
Let $n$ be the maximal number of not necessary distinct symbols
occurring in one of the lists $\lambda_1,...,\lambda_i$, i.e.
$n = \max \limits_{j=1}^{i}|\lambda_j|$.
We introduce the following numbers
which describe certain complexity properties of the special recursive system $S$:
\begin{itemize}
\item $k$ is the maximal number of distinct variables occurring
in an argument list of any $F \in B_S$, 
\item $\alpha$ is the maximal number of argument lists occurring
in a prime R-formula which is subformula of any $F \in B_S$,
\item $\rho$ is the maximal number of prime R-formulas occurring
in any $F \in B_S$.
\end{itemize}
Then there is an R-derivation $[\Lambda]$ of $p\,\lambda_1,...,\lambda_i$
with a number of steps  $|[\Lambda]|$ such that
$$
|[\Lambda]| \leq 
|B_S| \, \rho\,
\left( 1+ \alpha k 
\left(
\alpha\frac{n(n+1)}{2}  
\Gamma(n,k)
\right)^{\alpha}
\, \right)
\,.
$$

\underline{Remark:}~
This Theorem implies that for each 
$p\,\lambda_1,...,\lambda_i \in \Pi_R(S)$ 
there is an R-derivation $[\Lambda]$ of polynomial length
with respect to $n = \max \limits_{j=1}^{i}|\lambda_j|$.
We conclude that special recursive predicates are decidable.

\underline{Proof:}~An R-derivation $[\Lambda]$ of $p\,\lambda_1,...,\lambda_i$ 
can be chosen with the following properties:

1) All the R-formulas in $[\Lambda]$ are distinct.

2) $[\Lambda]$ starts with $[\Lambda_1] = B_S$, where
the R-axioms in $B_S$ are given in a fixed order (we may suppose that the formulas in $B_S$ are distinct).
 
3) Any application of the Substitution Rule
is restricted to the basis R-axioms, where each variable
is only replaced by elementary $A_S$-lists.

4) For all argument lists $\mu$ in $[\Lambda]$ 
with at least one variable occurring beyond $[\Lambda_1]$
there is an elementary list $\lambda$ which occurs as a sublist
in $p\,\lambda_1,...,\lambda_i$ such that $\mbox{Inst}(\mu,\lambda)$
is not empty. If $\mu$ is an elementary argument list in $[\Lambda] \setminus [\Lambda_1]$, 
then it must occur as a sublist in $p\,\lambda_1,...,\lambda_i$.

5) The Modus Ponens Rule is only applied if
all possible substitutions are done.

We extend $[\Lambda_1] = B_S$ given in 2) to a new R-derivation $[\Lambda_2]$
by applying the Substitution Rule on $[\Lambda_1]$ due to 3). 
In order to do this, we fix a given R-axiom $F \in B_S$ 
with R-conclusion $F_c = q\,\mu_1,...,\mu_l$, where $F$ may or may not have R-premises. 
We suppose that all argument lists $\mu$ in $F$ satisfy condition 4).
Due to $\mbox{var}(F) = \mbox{var}(F_c)$ it is sufficient 
to assign elementary $A_S$-lists to all variables in $F_c$
in order to get all possible substitutions which reduce $F$
to an elementary R-formula $F'$. Let 
$F'_c = q\,\mu'_1,...,\mu'_l$ result from $F_c$ by replacing
all the variables in $F_c$ by elementary $A_S$-lists.
Due to 4) we will only permit substitutions leading to
elementary $A_S$-lists $\mu'_1,...,\mu'_l$ which are sublists 
of the elementary $A_S$-lists $\lambda_1,...,\lambda_i$.
Due to Lemma (2.8) we have at most $\alpha \frac{n(n+1)}{2}$
possibilities to choose $\mu'_{\kappa}$ for any fixed $\kappa$. 
Due to Lemma (2.9) we have at most 
$\Gamma(|\mu'_{\kappa}|,k) \leq \Gamma(n,k)$ 
possibilities to assign elementary $A_S$-lists to all variables in 
$\mu_{\kappa}$ to obtain $\mu'_{\kappa}$. If we do these substitutions
for all $A_S$-lists $\mu_1,...,\mu_l$, we obtain at most
$$
\left(
\alpha\frac{n(n+1)}{2}  
\Gamma(n,k)
\right)^{l} \leq 
\left(
\alpha\frac{n(n+1)}{2}  
\Gamma(n,k)
\right)^{\alpha}
$$
elementary R-formulas $F'$ resulting from the substitutions 
of all variables in $F$. Since the total number of distinct variables
in $F$ or $F_c$ is bounded by $\alpha k$, we obtain the upper bound
$$
\alpha k\,
\left(
\alpha\frac{n(n+1)}{2}  
\Gamma(n,k)
\right)^{\alpha}
$$
of possible substitution steps, applied on the fixed R-axiom $F \in B_S$.
But $F$ is also part of $[\Lambda]$, and therefore we obtain the upper bound
$$
|[\Lambda_2]| \leq 
|B_S| \,
\left( 1+ \alpha k 
\left(
\alpha\frac{n(n+1)}{2}  
\Gamma(n,k)
\right)^{\alpha}
\, \right)
$$
for the number of steps of an R-derivation $[\Lambda_2]$, where 
$[\Lambda_2]$ is the part of $[\Lambda]$ which extends $[\Lambda_1] = B_S$ 
by applying the Substitution Rule.
This is possible due to the fifth property imposed on $[\Lambda]$.
The possible applications of the Modus Ponens Rule on $[\Lambda_2]$
yields $[\Lambda]$ with
$$
|[\Lambda]| \leq \rho\,|[\Lambda_2]|\,.
$$
From the last two inequalities we obtain Theorem (2.10). \dokend\\

\underline{Remark:} The proof of Theorem (2.10) enables the construction
of a determi\-nis\-tic resolution strategy which decides with polynomial effort
whether or not an 
elementary prime R-formula $p\,\lambda_1,...,\lambda_i$ is R-derivable 
in a special recursive system $S$. 
If $p\,\lambda_1,...,\lambda_i$ is R-derivable, then the algorithm constructs
an R-derivation $[\Lambda]$ obeying the restrictions 1)-5) 
in the proof of the Theorem.\\
 
Finally we mention that there are many other formalisms in order to generate 
recursively enumerable relations. One possible way is the definition of
recursive (or computable) functions for the non-negative integer numbers,
which can be formalized immediately in appropriate recursive systems,
or the use of Turing machines. Other approaches are given by Semi-Thue systems,
see Thue \cite{Th1}, \cite{Th2} and Jantzen \cite{Jn}, which are the
foundation for the use of grammars, see Hopcroft-Ullman \cite{HU},
and by logic programming, see Lloyd \cite{Ll}.
One very impressive result for the characterization of recursively enumerable
sets of positive integer numbers was finally solved by Matijasevi\u{c}
\cite{Mt1}, \cite{Mt2}, see also the extensive study of Davis \cite{Dv}:\\

{\bf (2.11) Theorem} (Matijasevi\u{c}, Robinson, Davis, Putnam)

One can construct a polynomial $M(y_1,...,y_n,z)$ with integer coefficients
such that for every recursively enumerable relation $R=R(x)$ of positive integer 
numbers $x$ there is a positive integer number $g_R$ with
\begin{align*}
R(x) ~ \Leftrightarrow &\mbox{\quad there are positive integer numbers~} k_1,...,k_n
\mbox{~such that~}\\ &\quad x=M(k_1,...,k_n,g_R)>0\,.
\end{align*}
 
This Theorem implies that the recursively enumerable sets of positive integer 
numbers are exactly the Diophantine sets. As a consequence, Hilbert's tenth
problem is unsolvable, i.e. there is no computing algorithm which will tell
of a given polynomial Diophantine equation with integer coefficients
whether or not it has a solution in integers.

\section{Embedding of a recursive system in a\\ mathematical system}
In this section we define a formal mathematical system which 
includes the predicate calculus and the structural induction with respect
to the recursively enumerable relations generated by an underlying recursive system
denoted by $S = [A_S;P_S;B_S]$. We will also define mathematical systems
with restricted argument lists.\\

{\bf (3.1) The symbols of the mathematical system}

Given are the following pairwise disjoint sets of symbols
\begin{enumerate}
\item[(a)]
A set $A_M \supseteq A_S$ of {\bf constant symbols} or {\bf operation symbols},
which must not be finite.
\item[(b)]
A set $P_M \supseteq P_S$ of {\bf predicate symbols}, which must not be finite.
\item[(c)]
The infinite alphabet $X$ of {\bf variable symbols} is the same
as in (1.1)(c). 
\item[(d)]
We define the following extension of the alphabet $E_R$ in (1.1)(d):
\[
E := [\, \sim \, ; \, ( \, ;\, ) \, ; \, , \, ; \, \to \, ;
\neg \, ; \leftrightarrow \, ; \& \, ;\, \vee \, ;
\, \forall \, ;\, \exists \,]\,.
\]
\end{enumerate}
If $A_M$ and $P_M$ are finite or denumerable then we may also assume that their symbols
are arranged in a fixed order and that $A_M$ and $P_M$ are extensions of finite alphabets $A_S$ and $P_S$, respectively.\\

{\bf (3.2) The basic structures of the mathematical system}
 
are the ($A_M$-)lists, ($A_M$-)sublists, ($A_M$-)terms and 
the elementary ($A_M$-)lists and ($A_M$-)terms, 
which are defined as in (1.1)-(1.4), 
but for the extended set $A_M$ instead of $A_S$.
The prime formulas and the elementary prime formulas are defined 
in the same way
as in (1.5) and (1.6), but with respect to the set $P_M$
of extended predicate symbols. 
Note that every prime R-formula is also a prime formula.

\newpage
{\bf (3.3) The formulas of the mathematical system}
\begin{enumerate}
\item[(a)]
Every prime formula is a formula.
\item[(b)]
Let $F$, $G$ be formulas and $x \in X$ be any variable.
Then the following expressions are formulas with the subformulas $F$, $G$, respectively.
\[
\neg \, F \, ; ~ \to F G \, ; ~ \leftrightarrow F G \, ; ~
\& F G \, ; ~ \vee F G \, ; ~ \forall \, x F \, ;~ \exists \, x F 
\]
\end{enumerate}
For example, if $f \in A_M$, $B \in P_M$ and $x,y \in X$, 
then the following expression is a formula of the mathematical system:
\[
\exists x ~~ \& ~ \forall x \, B \, x,y,f(xy) ~ \, \neg \sim x,y \,.
\]

The generalization of (1.7)(c) to subformulas is obvious.
A maximal sublist which occurs in a formula $F$ and which is not immediately following
$\forall$ or $\exists$ is also called an \underline{argument list} of $F$.
Finally we note that every R-formula is also a formula of the mathematical system.

{\bf (3.4) Variables in lists and formulas, free variables}
\begin{enumerate}
\item[(a)]
$\mbox{var}(\lambda)$ denotes the set of all variables
which occur in the list $\lambda$.
\item[(b)]
$\mbox{var}(F)$ denotes the set of all variables
occurring in a formula $F$.
\item[(c)] Recursive definition of $\mbox{free}(F)$, where
$F$, $G$ are formulas and  $x \in X$:\\
(i) ~ $\mbox{free}(F) = \mbox{var}(F)$ ~ for any prime formula $F$,\\
(ii) ~ $\mbox{free}(\neg \, F) = \mbox{free}(F)$,\\
(iii) ~ $\mbox{free}(J F G) = \mbox{free}(F) \cup \mbox{free}(G)$
for $J \in \{\to\, ; ~ \leftrightarrow\, ; ~
\& \, ; ~ \vee \, \}$.\\
(iv) ~ $\mbox{free}(\forall \, x F) = \mbox{free}(\exists \, x F)
= \mbox{free}(F) \setminus \{ x\}$.\\
\end{enumerate}

{\bf (3.5) The substitution of variables in lists (SbL)}

The expression $\mbox{SbL}(\lambda; \mu; x)\,=\,\lambda\,\frac{\mu}{x}$ 
describes the substitution of the variable $x$ in a list $\lambda$ by the list $\mu$. 
The following recursive definition of SbL holds for all lists $\lambda$, $\mu$, $\nu$,
for all $x,y \in X$ and $a, f \in A_M$
\begin{enumerate}
\item[(a)]
$
\mbox{SbL}(a;\mu;x)=a\,, \qquad
\mbox{SbL}(y;\mu;x)=\left\{
\begin{array}{r@{\quad,\quad}l}
y &  x \neq y\\ \mu & x = y\,.
\end{array}
\right.
$
\item[(b)]
$
\mbox{SbL}(f(\lambda);\mu;x)=f(\mbox{SbL}(\lambda;\mu;x))\,.
$
\item[(c)]
$
\mbox{SbL}(\lambda\,\mu;\nu;x)=
\mbox{SbL}(\lambda;\nu;x)\,\mbox{SbL}(\mu;\nu;x)\,.
$
\end{enumerate}

{\bf (3.6) The substitution of variables in formulas (SbF)}

The expression $\mbox{SbF}(F; \mu; x)=F\,\frac{\mu}{x}$ describes the substitution 
for each free occurrence of the variable $x$ in a formula $F$ 
by the list $\mu$. The recursive definition of SbF holds for all lists $\mu$, 
$\lambda_1$, $\lambda_2$,\,...\,, $\lambda_m$ $(m \geq 2)$, any $p \in P_M$,
$x,y \in X$\,, for all formulas $F$, $G$ and for\\ 
$J \in \{\to\, ; ~ \leftrightarrow\, ; ~\& \, ; ~ \vee \, \}$,
$Q \in \{\, \forall \, ; \, \exists\, \}$:
\begin{enumerate}
\item[(a)] let $\lambda_j' := \mbox{SbL}(\lambda_j ; \mu ; x)$ 
for $j=1,...,m$:
\[
\begin{array}{l@{\quad,\quad}l}
\mbox{SbF}(p \,;\,\mu \,;\, x) = p &
\mbox{SbF}(\sim \lambda_1 , \lambda_2 \,;\, \mu \,;\, x) = \,
\sim \lambda_1' , \lambda_2' \,, ~ \nonumber \\
\mbox{SbF}(p \, \lambda_1 \,;\, \mu \,;\, x) = p \, \lambda_1' &
\mbox{SbF}(p \, \lambda_1,\,...\,,\lambda_m \,;\, \mu \,;\, x) = 
p \, \lambda_1',\,...\,,\lambda_m' \,.\nonumber
\end{array}
\]
\item[(b)]
$\mbox{SbF}(\neg \, F; \mu; x) = \neg \, \mbox{SbF}(F; \mu; x)$ \,.
\item[(c)]
$\mbox{SbF}(J \, F G; \mu; x) = J ~
\mbox{SbF}(F; \mu; x) \, \mbox{SbF}(G; \mu; x)$ \,.
\item[(d)]
$
\mbox{SbF}(Q \, y \, F ; \mu ; x)=\left\{
\begin{array}{r@{\quad,\quad}l}
Q \, y \, F & x=y\\
Q \, y \, \mbox{SbF}(F ; \mu ; x) &  x \neq y\,. 
\end{array}
\right.
$
\end{enumerate}

{\bf (3.7) Avoiding collisions for the substitution SbF}

In order to ensure that the SbF-substitution of the variable $x$ by the list
$\mu$ in the formula $F$ is well defined we introduce the metamathematical 
predicate $\mbox{CF}(F;\mu; x)$, which means that $F$ and $\mu$ are
collision-free with respect to $x$. The recursive definition
of CF holds for all lists $\mu$, for all $x,y \in X$\,, 
for all formulas $F$, $G$ and for any
$J \in \{\to\, ; ~ \leftrightarrow\, ; ~\& \, ; ~ \vee \, \}$ \,,
$Q \in \{\, \forall \, ; \, \exists\, \}$:
\begin{enumerate}
\item[(a)] $\mbox{CF}(F;\mu; x)$ holds for any prime formula $F$.
\item[(b)]
$\mbox{CF}(\neg \, F; \mu; x)$ holds if and only if
$\mbox{CF}(F;\mu; x)$ holds.
\item[(c)]
$\mbox{CF}(J \, F G; \mu; x)$ holds if and only if
$\mbox{CF}(F; \mu; x)$ and $\mbox{CF}(G; \mu; x)$\\ are both satisfied.
\item[(d)]
$\mbox{CF}(Q \, y \, F ; \mu ; x)$ is satisfied if and only if:\\
i) ~ $x \not\in \mbox{free}(F) \setminus\{y\}$ ~\, or ~\,
ii)\, $y \not\in \mbox{var}(\mu)$ and $\mbox{CF}(F; \mu; x)$.
\end{enumerate}
\underline{Remarks:} The CF-condition is necessary in order to exclude 
undesired subs\-titutions like 
$\mbox{SbF}(\exists y \, \neg \,\sim x,y \,; y ; x)=
\exists y \, \neg \,\sim y,y$ with $x \neq y$. \\

It is also important to note that $x \not\in \mbox{free}(F)$ implies 
$\mbox{CF}(F;\mu; x)$ as well as $\mbox{SbF}(F;\mu; x)=F$.\\

{\bf (3.8) Propositional functions and truth values}

Let $\xi_1,...,\xi_j$ $(j \geq 1)$ be new distinct symbols, which are not occurring
in the given sets $A_M$, $P_M$, $X$, $E$
and \underline{not} part of the formal system. We call them propositional variables.
The propositional functions (of $\xi_1,...,\xi_j$) are defined as follows, where
$J \in \{\to\, ; ~ \leftrightarrow\, ; ~\& \, ; ~ \vee \, \}$:

\begin{tabular}{llll}
(a)   & $\xi_i$ & is a propositional function for $1 \leq i \leq j$. & \\
(b)   & $\neg \, \alpha$ & is a propositional function 
        if $\alpha$ is a propositional function. & \\
(c)   & $J \alpha \beta$ & is a propositional function if $\alpha$ and $\beta$ 
        are propositional functions. & \\
\end{tabular}\\

Let $\Psi : \{ \xi_1,...,\xi_j\} \to \{\,\top\,,\bot\,\}$ be any mapping 
which assigns to each propositional variable a truth value $\top$ for true 
or $\bot$ for false. Then we can canonically extend $\Psi$ to a function
$\bar{\Psi}$, which assigns to each propositional function 
of $\xi_1,...,\xi_j$
either the value $\top$ or $\bot$ according to 

\begin{tabular}{llll}
(d)   & $\bar{\Psi}(\neg \, \alpha)=\top$  &$\Leftrightarrow$&
         $\bar{\Psi}(\alpha)=\bot\,,$ \\
(e)   & $\bar{\Psi}(\to \, \alpha \, \beta)=\top$&$\Leftrightarrow$&
         $\bar{\Psi}(\alpha)=\bot$ or $\bar{\Psi}(\beta)=\top$,\\
(f)   & $\bar{\Psi}(\leftrightarrow \, \alpha \, \beta)=\top$ &$\Leftrightarrow$&
         $\bar{\Psi}(\alpha)=\top$ if and only if $\bar{\Psi}(\beta)=\top$,\\
(g)   & $\bar{\Psi}(\& \, \alpha \, \beta)=\top$ &$\Leftrightarrow$&
         $\bar{\Psi}(\alpha)=\top$ and $\bar{\Psi}(\beta)=\top$,\\
(h)   & $\bar{\Psi}(\vee \, \alpha \, \beta)=\top$ &$\Leftrightarrow$&
         $\bar{\Psi}(\alpha)=\top$ or $\bar{\Psi}(\beta)=\top$.\\
\end{tabular}

Here $\bar{\Psi}(\alpha)=\bot$ $\Leftrightarrow$ not $\bar{\Psi}(\alpha)=\top$
holds for all propositional functions $\alpha$.

A propositional function $\alpha=\alpha(\xi_1,...,\xi_j)$ is called 
\underline{identically true}, if there holds $\bar{\Psi}(\alpha) = \top$ for 
\underline{every} mapping $\Psi : \{ \xi_1,...,\xi_j\} \to \{\,\top\,,\bot\,\}$.\\

{\bf (3.9)  The axioms of the propositional calculus} 

Let $\alpha=\alpha(\xi_1,...,\xi_j)$ be a propositional function
of the distinct propositional variables $\xi_1,...,\xi_j$, $j \geq 1$. 
Let $F_1$,...,$F_j$ be formulas and suppose that $\alpha$ is identically true.
Then the formula $F := \alpha(F_1,...,F_j)$ is an axiom of the propositional calculus.\\

{\bf (3.10) The axioms of equality} 

Let $x,y \in X$ and let $\lambda, \mu$ be $A_M$-lists.
Then the following formulas are axioms of equality

\begin{tabular}{llll}
(a)  & $ \sim x,x~.$ &&\\
(b)  & $ \to ~\mbox{SbF}(\,\sim \lambda,\mu\,;\,x\,;\,y\,) ~            
         \to ~ \sim x,y ~
         \sim \lambda,\mu\,.$ &&\\
\end{tabular}

Let $p \in P_M$, $n \geq 1$ and $x_1,y_1,...,x_n,y_n \in X$.
Then the following formula is an axiom of equality

\begin{tabular}{llll}
(c)  & $\to ~ \sim x_1,y_1 ~...~ \to ~ \sim x_n,y_n ~ \to ~ 
p\,x_1,...,x_n~p\,y_1,...,y_n.$ &&\\
\end{tabular}\\

{\bf (3.11) The quantifier axioms}

Let $F$, $G$ be formulas and $x \in X$. Then we define the quantifier axioms
 
\begin{tabular}{llll}
(a)   & $\to ~ \forall \, x F  ~~ F$ & & \\
(b)   & $\to ~ \forall \, x ~ \to F G ~~
        \to ~ F ~ \forall \, x \,G\,,$ & if $x \not\in \mbox{free}(F)$& \\
(c)   & $\leftrightarrow ~ \neg\,\forall \, x \, \neg\,F ~~\exists \, x \,F$\,. & &\\
\end{tabular}\\

{\bf (3.12) The mathematical system M} is given for fixed $X$ and $E$ 
\begin{enumerate}
\item[(i)] by the recursive system $S=[A_S;P_S;B_S]$ defined in (1.10)\,,
\item[(ii)] by the sets $A_M \supseteq A_S$ and $P_M \supseteq P_S$
and by a set $B_M \supseteq B_S$ of formulas in $M$. 
The formulas of $B_M$ are called the basis axioms of the mathematical system $M$.
Recall that $A_M$, $P_M$, $X$ and $E$ are pairwise disjoint.
Often we have that $A_M$ and $P_M$ are countable or even finite sets, or that
$B_M$ is recursively solvable, i.e. decidable, 
but this must not be required in the general case.
\end{enumerate}
The mathematical system may be denoted by $M = [S; A_M; P_M; B_M]$.\\
The \underline{axioms of M} are the axioms 
of the propositional calculus, the axi\-oms of equality, 
the quantifier axioms and the formulas in $B_M$.\\

{\bf (3.13)  Rules of inference and (formal) proofs in M} 

A (formal) proof in $M$ is a
list $[\Lambda]:=[F_1;...;F_l]$ of formulas $F_1,...,F_l$ 
including the empty list $[\,]$.
The formulas $F_1,...,F_l$ are the steps of the proof, 
which is generated by the rules of inference
\begin{enumerate}
\item[(a)] \underline{Axiom Rule}:~ 
The empty list $[\,]$ is a proof in the mathematical system $M$.
If $[\Lambda]$ is a proof and $F$ an axiom, 
then $[\Lambda;~F]$ is also a proof.
\item[(b)] \underline{Modus Ponens Rule}:~
Let $F$, $G$ be two formulas and $F$, $\to ~ F ~ G$ both steps 
of the proof $[\Lambda]$.
The $[\Lambda;~G]$ is also a proof.
\item[(c)] \underline{Substitution Rule}:~ 
Let $F$ be a step of the proof $[\Lambda]$, $x \in X$ and $\lambda$ a list.
If $\mbox{CF}(F;\lambda;x)$ holds,
then $[\Lambda;~\mbox{SbF}(F;\lambda;x)]$ is also a proof.
\item[(d)] \underline{Generalization Rule}:~
Let $F$ be a step of the proof $[\Lambda]$, $x \in X$.
Then $[\Lambda;~\forall \, x F]$ is also a proof.
Here it is not required that $x$ occurs in $F$.
\item[(e)] \underline{Induction Rule}:~
In the following we fix a predicate symbol $p \in P_S$,
a list $x_1,...,x_i$ of $i\geq 0$ distinct variables
and a formula $G$ in $M$. 
We suppose that $x_1,...,x_i$ and the variables of $G$
are not involved in $B_S$.

Then to every R-formula $F$ of $B_S$ there corresponds exactly one formula $F'$ 
of the mathematical system, which is obtained if we replace in $F$ each 
$i-$ary subformula $p \, \lambda_1,...,\lambda_i$, 
where $\lambda_1,...,\lambda_i$ are lists, by the formula 
$G ~ \frac{\lambda_1}{x_1}...\frac{\lambda_i}{x_i}$.
 
If $F'$ is a step of a proof $[\Lambda]$ for all R-formulas $F \in B_S$ 
for which $p$ occurs $i-$ary in the R-conclusion of $F$, 
then $[\Lambda; ~ \to ~p\,x_1,...,x_i ~G]$ is also a proof.
\end{enumerate}

\underline{Remarks on the rules of inference:} 

Any R-derivation in the recursive system $S=[A_S;P_S;B_S]$ is also a proof
in the mathematical system $M=[S;A_M;P_M;B_M]$ due to the first three Rules (a)-(c),
due to $A_M \supseteq A_S$, $P_M \supseteq P_S$, $B_M \supseteq B_S$ 
and due to the fact that every R-axiom of equality
is also an axiom of equality in the mathematical system $M$. 
Rule (e) enables the structural induction with respect to the recursively
enumerable relations represented in S. If we put $P_S=[\,]$,
then the use of the Induction Rule (e) is suppressed.

The axioms of the propositional calculus can also be reduced to 
axiom schemes resulting from a small list of 
identically true propositional functions, 
which requires an own study of the propositional calculus.\\

{\bf (3.14) Provable formulas}

The steps of a proof $[\Lambda]$ are called \underline{provable} formulas. 
By $\Pi(M)$ we denote the set of all provable formulas $F$ in $M$.\\

\underline{Example 1:}~Let $A_S:=[\,0\,;\,'\,]$, $B,C,D \in P_S$
and $x,y,z \in X$ be distinct symbols. 
We consider the recursive system $S=[A_S;P_S;B_S]$ with the complete list of basis R-axioms given by

\begin{tabular}{llll}
(1)   & $B\,0$ & 
\qquad \qquad (2) ~$\to ~ B\,x \,~  B\,x'$ & \\
(3)   & $\to ~ B\,x ~ C\,x$ & 
\qquad \qquad (4) ~$\to ~ B\,x \,~ \to ~ C\,y \,~ C\,xy$ & \\
(5)   & $\to ~ B\,x ~ D\,x$ & 
\qquad \qquad (6) ~$\to ~ B\,x \,~ \to ~ C\,y \,~ D\,xy\,.$ & \\
\end{tabular}

The basis R-axioms (1)-(6) of $B_S$ form a proof in 
any mathe\-matical system $M=[S;A_M;P_M;B_M]$ which can be extended as follows 

\begin{tabular}{llll}
(7)   & $\to ~ D\,z ~ C\,z$ & 
\mbox{with Rule (e) and (3)-(6) for $p\,z=D\,z$, $G=C\,z$}& \\
(8)   & $\to ~ D\,y ~ C\,y$ & \mbox{with Rule (c) and (7)}& \\
\end{tabular}\\
\begin{tabular}{llll}
(9)   & $\to ~~ \to \,D\,y \, C\,y ~~ \to ~~ \to\,B\,x\,\to\,C\,y\,D\,xy ~~~
\to\,B\,x\to\,D\,y\,D\,xy$ & &\\
  & with (3.9) and the identically true propositional function & &\\
  & 
  $
    \alpha(\xi_1,\xi_2,\xi_3,\xi_4)~:= ~
   \to ~\, \to \,\xi_1 \, \xi_2 ~\, \to ~\, \to\,\xi_3\,\to\,\xi_2\,\xi_4 ~\,
   \to\,\xi_3\to\,\xi_1\,\xi_4
  $ & &\\
\end{tabular}\\
\begin{tabular}{llll}
(10)   & $\to ~~ \to\,B\,x\,\to\,C\,y\,D\,xy ~~\to\,B\,x\to\,D\,y\,D\,xy$ & &\\
       & with Rule (b) and (8), (9)& &\\
\end{tabular}\\
\begin{tabular}{llll}
(11)   & $\to\,B\,x\to\,D\,y\,D\,xy$ & with Rule (b) and (6), (10)&\\
\end{tabular}\\
\begin{tabular}{llll}
(12)   & $\to ~ C\,z ~ D\,z$ & 
\mbox{with Rule (e) and (3),(4),(5),(11) for $p\,z=C\,z$, $G=D\,z$}& \\
\end{tabular}\\
\begin{tabular}{llll}
(13)   & $\to ~~ \to ~ C\,z ~ D\,z ~~\to~~ \to ~ D\,z ~ C\,z ~~
\leftrightarrow ~ C\,z ~ D\,z$ & with (3.9) &\\
\end{tabular}\\
\begin{tabular}{llll}
(14)   & $\to ~~ \to ~ D\,z ~ C\,z ~~\leftrightarrow ~ C\,z ~ D\,z$ 
& with Rule (b) and (12), (13)&\\
(15)   & $\leftrightarrow ~ C\,z ~ D\,z$ 
& with Rule (b) and (7), (14)&\\
(16)   & $\forall \, z \,\leftrightarrow ~ C\,z ~ D\,z$ 
& with Rule (d) and (15)\,.&\\
\end{tabular}\\

\underline{Example 2:}~We consider a mathematical system $M = [S;\,A_M;\,P_M;\,B_M]$, 
fix a predicate symbol $p \in P_S$ of the recursive system
and a non-negative integer number $i \geq 0$. 
We suppose that there is no $i$-ary R-conclusion of the form $p\,\lambda_1,...,\lambda_i$ 
in the R-formulas of $B_S$. We consider a list of distinct new variables $x_1,...,x_i$ 
and obtain the following proof $[\Lambda]$ of $\neg\,p\,x_1,...,x_i$ in $M$
due to the Induction Rule (e):

\begin{tabular}{llll}
 $[\Lambda] = $& $[~\to ~ p\,x_1,...,x_i ~ \neg\,p\,x_1,...,x_i\,;$ &&\\
& $~\to ~ \to ~ p\,x_1,...,x_i ~ \neg\,p\,x_1,...,x_i ~ \neg\,p\,x_1,...,x_i\,;$ &&\\
& $~\neg\,p\,x_1,...,x_i\,]$\,. &&\\
\end{tabular}\\

\underline{Example 3:}~Let $A_S:=[\,a\,]$, $N,< \,\in P_S$ and
$x,y,u,v \in X$ be distinct symbols. We consider the complete list 
of basis R-axioms given by

\begin{tabular}{llll}
(1) & $N\,a$ & &\\
\end{tabular}\\
\begin{tabular}{llll}
(2) & $\to ~ N\,x \,~  N\,xa$ & &\\
\end{tabular}\\
\begin{tabular}{llll}
(3) & $\to ~ N\,x ~ \to ~ N\,y ~<\,x,xy$ & & \\
\end{tabular}

As in the first example they form a proof which will be extended by

\begin{tabular}{llll}
(4) & $\to ~ N\,u \,~  N\,ua$ & &\\
\end{tabular}\\
\begin{tabular}{llll}
(5) & $\to ~ N\,v \,~  N\,va$ & &\\
\end{tabular}\\
\begin{tabular}{llll}
(6) & $\to ~ N\,uv \,~  N\,uva$ & &\\
\end{tabular}\\
\begin{tabular}{llll}
(7) &$\to ~ \to ~ N\,v \,~  N\,va ~~ \to ~ \to ~ N\,uv \,~  N\,uva$ &&\\
& $\qquad \to ~ \to ~ N\,u \,~\&  N\,v \, N\,uv~~ \to ~N\,u \,~\&  
N\,va \, N\,uva $ &&\\
\end{tabular}\\
\begin{tabular}{llll}
(8) &$\to ~ \to ~ N\,uv \,~  N\,uva~
\to ~ \to ~ N\,u \,~\&  N\,v \, N\,uv~~ \to ~N\,u \,~\&  N\,va \, N\,uva $ &&\\
\end{tabular}\\
\begin{tabular}{llll}
(9) &$\to ~ \to ~ N\,u \,~\&  N\,v \, N\,uv~~ \to ~N\,u \,~\&  N\,va \, N\,uva $ &&\\
\end{tabular}

In (7) we have used the identically true propositional function

$\to ~ \to \, \xi_1\,\xi_2~\to~ ~\to\,\xi_3 \,\xi_4~\to~ ~\to\,
    \xi_5 \,\&\, \xi_1\,\xi_3~\to\,\xi_5 \,\&\, \xi_2\,\xi_4$.

\begin{tabular}{llll}
(10) &$\to ~ N\,a $ &&\\
 & $\to~ \to ~ N\,u \, N\,ua~  $ &&\\
 & $\qquad \to ~ N\,u \,~\&  N\,a ~ N\,ua$ &&\\
\end{tabular}\\
\begin{tabular}{llll}
(11) 
 & $\to~ \to ~ N\,u \, N\,ua~  $ &&\\
 & $\qquad \to ~ N\,u \,~\&  N\,a ~ N\,ua$ &&\\
\end{tabular}\\
\begin{tabular}{llll}
(12)  & $\qquad \to ~ N\,u \,~\&  N\,a ~ N\,ua$ &&\\
\end{tabular}\\
\begin{tabular}{llll}
(13) &$\to ~ \to ~ N\,u \,~\&  N\,x \, N\,ux~~ \to ~N\,u \,~\&  N\,xa \, N\,uxa $ &&\\
\end{tabular}\\
\begin{tabular}{llll}
(14) &$\to ~N\,v ~\to ~ N\,u \,~\&  N\,v \, N\,uv$ &&\\
\end{tabular}

In (10)-(13) we have prepared the 
first application of the Induction Rule.
For (10) we use
$\to ~ \xi_1 ~ \to ~\to\,\xi_2\,\xi_3
~ \to\,\xi_2 \,\&\, \xi_1\,\xi_3$ as an identically true propositional function.
Formula (13) results from (9) and Rule (c)
and formula (14) from (12), (13), (1), (2) and Rule (e).

Finally we listen the remaining steps of the proof

\begin{tabular}{llll}
(15) &$\to ~~\to ~N\,v ~\to ~ N\,u \,~\&  N\,v \, N\,uv$ &&\\
      &$\qquad      ~\to ~N\,u ~\to ~ N\,v \,~ N\,uv$ &Rule (a)&\\
\end{tabular}\\
\begin{tabular}{llll}
(16) &$\to ~N\,u ~\to ~ N\,v \,~ N\,uv$ &Rule (b)&\\
\end{tabular}\\
\begin{tabular}{llll}
(17)  
 &$\to  ~ N\,x ~\to~ N\,v \, N\,xv$ &Rule (c)&\\ 
\end{tabular}\\
\begin{tabular}{llll}
(18)  
 &$\to  ~ N\,x ~\to~ N\,y \, N\,xy$ &Rule (c)&\\ 
\end{tabular}\\
\begin{tabular}{llll}
(19)  
 &$\to  ~ N\,x ~\to~ N\,y \, N\,x$ &Rule (a)&\\ 
\end{tabular}\\
\begin{tabular}{llll}
(20)  
 &$\to  ~ <\,u,v \, N\,v$ &\qquad ~Rule (e)&\\ 
\end{tabular}\\
\begin{tabular}{llll}
(21)  
 &$\to  ~ <\,u,v \, N\,u$ &\qquad ~Rule (e)&\\ 
\end{tabular}\\
\begin{tabular}{llll}
(22)  
 &$\to  ~ \to  ~ <\,u,v \, N\,v$ &&\\
 &$\to  ~ \to  ~ <\,u,v \, N\,u$ &&\\  
 &$\qquad \to  ~ <\,u,v ~ \& \, N\,u \, N\,v $ &Rule (a)&\\  
\end{tabular}\\
\begin{tabular}{llll}
(23)  
 &$\to  ~ \to  ~ <\,u,v \, N\,u$ &&\\  
 &$\qquad \to  ~ <\,u,v ~ \& \, N\,u \, N\,v $ &Rule (b)&\\  
\end{tabular}\\
\begin{tabular}{llll}
(24)  
 &$\to  ~ <\,u,v ~ \& \, N\,u \, N\,v $ &\qquad Rule (b)&\\  
\end{tabular}\\
\begin{tabular}{llll}
(25)  
 &$\forall \, v ~\to  ~ N\,u ~\to~ N\,v \, N\,uv$ &Rule (d)&\\ 
\end{tabular}\\
\begin{tabular}{llll}
(26)  
 &$\forall \, v ~\to  ~ <\,u,v ~ \& \, N\,u \, N\,v $ &Rule (d)&\\  
\end{tabular}

We finally end up with the two formulas, using again Rule (d) 

\begin{tabular}{llll}
(27)  
 &$\forall \, u \, \forall \, v ~\to  ~ N\,u ~\to~ N\,v \, N\,uv$ &&\\ 
\end{tabular}\\
\begin{tabular}{llll}
(28)  
 &$\forall \, u\,\forall \, v ~ \to  ~ <\,u,v ~ \& \, N\,u \, N\,v $ &&\\  
\end{tabular}\\

Now we consider mathematical systems with
given restrictions for the argument lists in their formulas. This is
important since we are often concerned with the representation of
functions with a given number of arguments or with special lists of terms.

The restriction of the argument lists is described by a subset of lists
which contains the variables and which is invariant with respect 
to substitutions. This is described in the next definition.

\newpage
{\bf (3.15) Mathematical systems with restricted argument lists}\\
Let $M = [S; A_M; P_M; B_M]$ be a mathematical system and ${\cal L}$
a given subset of $A_M$-lists with the properties 
\begin{itemize}
\item[(i)] $X \subseteq {\cal L} $\,,
\item[(ii)] $\lambda \frac{\mu}{x} \in {\cal L}$ ~
for all $\lambda, \mu \in {\cal L}$, $x \in X$\,,
\item[(iii)] all formulas in $B_M$ contain only argument lists in ${\cal L}$\,.
\end{itemize}
Then $[M;{\cal L}]$ is called a
mathematical system with restricted argument lists.
A formula in $[M;{\cal L}]$ is a formula in $M$ 
which has only argument lists in ${\cal L}$\,.
A proof $[\Lambda]$ in $[M;{\cal L}]$ is a proof in $M$
with the restrictions 
\begin{itemize}
\item[(iv)] the formulas in $[\Lambda]$ and 
the formulas $F$ and $G$ in (3.13)(a)-(e)
contain only argument lists in ${\cal L}$\,,
\item[(v)] there holds $\lambda \in {\cal L}$ 
for the list $\lambda$ in (3.13)(c)\,.
\end{itemize}
By $\Pi(M;{\cal L})$ we denote the set of provable formulas 
in $[M;{\cal L}]$\,. 


\underline{Example 4:}~ The Peano arithmetic PA

Let $\tilde{S}$ be the recursive system $\tilde{S}=[\tilde{A};\tilde{P};\tilde{B}]$ 
where $\tilde{A}$,$\tilde{P}$ and $\tilde{B}$ are empty, and put $A_{PA} = [\,0;\,s\,;\,+\,;\,*\,]$, \
$P_{PA} = [\,]$.\\ Next we define the set ${\cal L}$ of
\underline{numeral terms} by the recursive definition

\begin{tabular}{llll}
(i)& $0$ and $x$ are numeral terms for any $x \in X$.&&\\
(ii)& If $\vartheta$ is a numeral term, then also $s(\vartheta)$.&&\\
(iii)& If $\vartheta_1$, $\vartheta_2$ are numeral terms,
then also $+(\vartheta_1 \vartheta_2)$ and $*(\vartheta_1 \vartheta_2)$.&&\\
\end{tabular}\\

We define the mathematical system $M'=[\tilde{S};A_{PA};P_{PA};B_{PA}]$ 
by giving the following basis axioms for $B_{PA}$ with distinct variables $x,y$

\begin{tabular}{llll}
(1)  & $\forall\,x ~ \sim +(0x),x$&&\\
(2)  & $\forall\,x \, \forall\,y ~\sim +(s(x)y),s(+(xy))$&&\\
\end{tabular}\\
\begin{tabular}{llll}
(3)  & $\forall\,x  ~ \sim *(0x),0$&&\\
(4)  & $\forall\,x \, \forall\,y ~ \sim *(s(x)y),+(*(xy)y)$&&\\
\end{tabular}\\
\begin{tabular}{llll}
(5) & $\forall\,x \, \forall\,y ~ \to ~ \sim s(x),s(y)\,~ \sim x,y$&&\\
(6) & $\forall x  ~\neg \sim s(x),0\,.$&&\\
\end{tabular}

Moreover, for all formulas $F$ (with respect to $A_{PA}$ and $P_{PA}$)
which have only numeral argument lists, the following formulas
belong to $B_{PA}$ according to the Induction Scheme

\begin{tabular}{llll}
(IS) & $\to ~~ \forall\,x\,~ \& ~ \mbox{SbF}(F;0\,;x) ~\, 
   \to~ F\, \mbox{SbF}(F;s(x)\,;x)  ~~
   \forall\,x\,F$ \,.&&\\ 
\end{tabular}

The system PA of \underline{Peano arithmetic} is given by
PA = $[M';{\cal L}]$, i.e. the argument lists of PA are restricted 
to the set ${\cal L}$ of numerals.
The Induction Rule (3.13)(e) is not used in PA since 
$\tilde{A}$,$\tilde{P}$ and $\tilde{B}$ are empty here
and since we are using the Induction Scheme (IS).

The following formulas are provable in PA for all $x,y,z \in X$:

\begin{tabular}{llll}
& $\forall\,x \, \forall\,y\, \forall\,z ~\sim +(+(xy)z),+(x+(yz))$& and ~
$\forall\,x \, \forall\,y ~\sim +(xy),+(yx)$\,,&\\
& $\forall\,x \, \forall\,y\, \forall\,z ~\sim *(*(xy)z),*(x*(yz))$& and ~
$\forall\,x \, \forall\,y ~\sim *(xy),*(yx)$\,,&\\
\end{tabular}\\

and also the most part of usual number theory. \\

{\bf (3.16) Lemma}

Let $[M;{\cal L}]$ be a mathematical system with the set ${\cal L}$
of restricted argument lists, $F$, $G$ formulas in $[M;{\cal L}]$
and $x,y \in X$.

(a) If $y \notin\mbox{var}(F)$, then

\quad (i)   $\mbox{CF}(F;\,y;\,x)$ ~and ~
\quad (ii)  $\mbox{CF}(F\,\frac{y}{x};\,x;\,y)$ ~and~
\quad (iii) $F\,\frac{y}{x}\frac{x}{y}\,=\,F$\,.

Moreover, the following formulas are provable in $[M;{\cal L}]$:
 
\begin{tabular}{llll}
(b)   & $\to ~ \forall \, x \to F G ~~\to~\forall \, x F ~ \forall \, x G$ & &\\
(c)   & $\leftrightarrow ~ \forall \, x ~ \to F G ~~
        \to ~ F ~ \forall \, x \,G\,,$ & if $x \not\in \mbox{free}(F)$& \\
(d)   & $\leftrightarrow ~ \forall \, x F ~~\, 
      \forall\,y \,F\,\frac{y}{x}\,,$ & if $y \not\in \mbox{var}(F)\,.$& 
\end{tabular}\\

\underline{Proof:} (a) is shown
by induction with respect to the formula $F$
and is needed for part (d) of the Lemma.

(b) From the quantifier axiom (3.11)(a) we know that the formulas\\
$\rightarrow ~ \forall \, x  \to F G ~\to  F G$ and 
$\to \, \forall \, x F  ~ F$ are both provable in $[M;{\cal L}]$.\\
From these formulas and the propositional calculus we can conclude
that $\to ~ \forall \, x \to F G ~\to~\forall \, x F ~ G$
is also provable in $[M;{\cal L}]$. We conclude that
$\forall \, x \, \to ~ \forall \, x \to F G ~\to~\forall \, x F ~ G$
is provable in $[M;{\cal L}]$ due to Rule (3.13)(d),
and on the last formula we can apply
two times the quantifier axiom (3.11)(b) and the propositional calculus  
in order to infer the desired formula.

(c) We must only show the backward implication ``$\leftarrow$'' and suppose
that $x \not\in \mbox{free}(F)$. From the quantifier axiom (3.11)(a) 
and the propositional calculus we can infer the formulas
$\to \, \forall \, x G  ~ G$ and $\to \, \to \, F ~ \forall \, x G ~ \to F G$,
and hence $\forall \, x \, \to \, \to  F ~ \forall \, x G ~ \to F G$
due to Rule (3.13)(d). From the
quantifier axiom (3.11)(b) and the propositional calculus 
we can infer the desired formula
$~\to ~ \to \, F \, \forall \, x G ~ \forall \, x\,\to F G$.

(d) Suppose that $y \notin \mbox{var}(F)$. 
For the forward implication ``$\rightarrow$'' we use 
$\mbox{CF}(F;\,y;\,x)$ from part (a) of the Lemma
and conclude that from $\to \, \forall \, x F  ~ F$ and Rules
(3.13)(c),(d) we can infer $\to \, \forall \, x F  ~ F\frac{y}{x}$
and $\forall y\, \to \, \forall \, x F  ~ F\frac{y}{x}$.
The quantifier axiom (3.11)(b) and the propositional calculus 
admit to infer the formula $\rightarrow ~ \forall \, x F ~~\, 
      \forall\,y \,F\frac{y}{x}$\,. The opposite direction
``$\leftarrow$'' can be shown in the same way, using the remaining
part (a)(ii) and (iii) of the Lemma.     
\dokend \\

{\bf (3.17) Theorem}

Let $[M;{\cal L}]$ be a mathematical system with the set ${\cal L}$
of argument lists.
\begin{itemize}
\item[(a)] Equivalence Theorem\\
Let $H$, $H'$ be formulas in $[M;{\cal L}]$ such
that \,$\leftrightarrow\, H \, H' \in \Pi(M;{\cal L})$.
Let $F,F'$ be any two formulas  in $[M;{\cal L}]$ such that
$F'$ results from $F$ if $H$ is replaced by $H'$ 
at certain places in $F$ where $H$ occurs as a subformula.
Then $\leftrightarrow\, F \, F' \in \Pi(M;{\cal L})$.
\item[(b)] Replacement of bound variables\\
Let $G$ be a formula in $[M;{\cal L}]$. Suppose that $G$ 
contains a subformula of the form $Q\,x\,F$ with 
$Q \in \{\forall\,,\,\exists\}$, $x \in X$.
Let $y$ be a second variable, 
which does not occur in the formula $F$.
Let $G'$ result from $G$ by replacing the subformula $Q\,x\,F$ everywhere or only
at certain places in $G$ by the formula $Q\,y\,\mbox{SbF}(F;y;x)$. \\
Then ~ $\leftrightarrow\,G\,G' \in \Pi(M;{\cal L})$. 
\end{itemize}

\underline{Proof:} We employ induction with respect to the rules
of forming formulas.

(a) Suppose that $\leftrightarrow\,F\,F' \in \Pi(M;{\cal L})$ and that
$\leftrightarrow\,G\,G' \in \Pi(M;{\cal L})$ for formulas 
$F, F', G, G'$ in $[M;{\cal L}]$. This is automatically satisfied for
$F=F'$, $G=G'$. Let be $J \in [\,\to\,;\,\leftrightarrow\,;\,\&\,;\,\vee\,]$.
Then we can first state due to the propositional calculus that the formulas

\begin{tabular}{llll} 
 &$\to ~ \leftrightarrow F F' ~ \leftrightarrow \neg\, F \neg\, F'$ & ~and~
 $\leftrightarrow \neg\, F \neg\,F'$\,,&\\ 
\end{tabular}\\
\begin{tabular}{llll}
 &$\to ~ \leftrightarrow F F' ~ \to ~ \leftrightarrow G G' 
 ~\leftrightarrow J\, F \, G ~ J\, F' \, G'$ & ~and~
 $\leftrightarrow J\, F \, G ~ J\, F' \, G'$&\\ 
\end{tabular}

also belong to $\Pi(M;{\cal L})$. 
There remains the more interesting induction step for the quantifiers. 

We obtain $\to F F' \in \Pi(M;{\cal L})$ as well as 
$\forall x \,\to F F' \in \Pi(M;{\cal L})$ 
due to the assumption $\leftrightarrow F F' \in \Pi(M;{\cal L})$, 
the axioms of the propositional calculus and 
due to the Rules (3.13)(a),(b),(d). 
Therefore we can infer from Lemma (3.16)(b) and Rule (3.13)(b) that 
$\to\,\forall x \,F \, \forall x \,F' \in \Pi(M;{\cal L})$.
The propositional calculus implies that not only
$\to F F' \in \Pi(M;{\cal L})$ but also $\to F' F \in \Pi(M;{\cal L})$, 
and therefore we can repeat the arguments above with interchanged roles of 
$F$ and $F'$ to obtain 
$\to\,\forall x \,F' \, \forall x \,F \in \Pi(M;{\cal L})$.
Applying again the propositional calculus on $\to\,\forall x \,F \, \forall x \,F'$
and $\to\,\forall x \,F' \, \forall x \,F$ we conclude
that $\leftrightarrow\,\forall x \,F \, \forall x \,F' \in \Pi(M;{\cal L})$.

Finally we have to show that 
$\leftrightarrow\,\exists x \,F \, \exists x \,F' \in \Pi(M;{\cal L})$.\\ 
From $\leftrightarrow \neg\,F \neg\,F' \in \Pi(M;{\cal L})$ we obtain that 
$\leftrightarrow\,\forall x \,\neg\,F \, \forall x \,\neg\,F' \in \Pi(M;{\cal L})$
by the induction step for the $\forall$-quantifier proven above.
The propositional calculus and the quantifier axiom (3.11)(c) imply that 
$\leftrightarrow\,\neg\,\forall x \,\neg\,F \, \neg\,
\forall x \,\neg\,F' \in \Pi(M;{\cal L})$ and
$\leftrightarrow\,\exists x \,F \,\exists x \,F' \in \Pi(M;{\cal L})$.
Thus we have shown the first part. 

(b) The proof is clear for $Q=\forall$ due to Lemma (3.16)(d) 
and part (a). For $Q=\exists$ we replace in Lemma (3.16)(d) 
the formula $F$ by $\neg\,F$ and conclude that 
$\leftrightarrow ~ \forall \, x \,\neg\,F ~~\, 
      \forall\,y \,\neg\,\mbox{SbF}(F;y;x)$
and 
$\leftrightarrow ~ \neg\,\forall \, x \,\neg\,F ~~\, 
      \neg\,\forall\,y \,\neg\,\mbox{SbF}(F;y;x)$
are both members of $\Pi(M;{\cal L})$. 
On the last formula we can apply the quantifier-axiom
(3.11)(c) and the propositional calculus in order to obtain that\\
$\leftrightarrow ~ \exists \, x \,F ~~\, 
      \exists\,y \,\mbox{SbF}(F;y;x)$
is a member of $\Pi(M;{\cal L})$. 
In this case we can also apply part (a).\dokend 

{\bf (3.18) Proposition}

Let $[M;{\cal L}]$ be a mathematical system with restricted argument lists.
The following formulas are provable in $[M;{\cal L}]$ for
any formulas $F$ and $G$ in $[M;{\cal L}]$ and $x,y \in X$ 

\begin{tabular}{llll}
(1)   & $\to ~ \forall \, x F  ~~ F$ & & \\
(2)   & $\to ~ F  ~~ \exists \, x F$ & & \\
\end{tabular}\\
\begin{tabular}{llll}
(3)   & $\leftrightarrow ~ \neg\,\forall \, x \, \neg\,F ~~\exists \, x \,F$& &\\
(4)   & $\leftrightarrow ~ \neg\,\exists \, x \, \neg\,F ~~\forall \, x \,F$& &\\
\end{tabular}\\
\begin{tabular}{llll}
(5)   & $\leftrightarrow ~ \forall \, x F ~\,F\,,$ & if $x \not\in \mbox{free}(F)$&\\
(6)   & $\leftrightarrow ~ G ~~\exists \, x G \,,$  & if $x \not\in \mbox{free}(G)$&\\
\end{tabular}\\
\begin{tabular}{llll}
(7)   & $\leftrightarrow ~ \forall \, x F ~~\, 
      \forall\,y \,\mbox{SbF}(F;y;x)\,,$ & if $y \not\in \mbox{var}(F)$&\\
(8)   & $\leftrightarrow ~ \exists \, x F ~~\, 
      \exists\,y \,\mbox{SbF}(F;y;x)\,,$ & if $y \not\in \mbox{var}(F)$&\\
\end{tabular}\\
\begin{tabular}{llll}
(9)   & $\leftrightarrow ~ \forall \, x \, \forall \, y \, F ~~ 
      \forall \, y \, \forall \, x \, F$ & &\\
(10)   & $\leftrightarrow ~ \exists \, x \, \exists \, y \, F ~~ 
      \exists \, y \, \exists \, x \, F$ & &\\
\end{tabular}\\
\begin{tabular}{llll}
(11)   & $\to ~ \forall \, x \to F G ~~\to ~\forall \, x F ~ \forall \, x \,G$ & &\\
(12)   & $\to ~ \forall \, x \to F G ~~\to ~\exists \, x F ~ \exists \, x \,G$ & &\\
\end{tabular}\\
\begin{tabular}{llll}
(13)   & $\to ~ \& \,\exists \, x F ~ \forall \, x \, G~~ \exists \, x \, \& F G$ & &\\
(14)   & $\to ~ \forall \, x \vee F G ~~ \vee \, \forall \, x F ~ \exists \, x \,G$ & &\\
\end{tabular}\\
\begin{tabular}{llll}
(15)   & $\leftrightarrow ~ \forall \, x ~ \& F G ~~ 
      \& ~ \forall\,x \,F ~ \forall\,x \,G$ & &\\
(16)   & $\leftrightarrow ~ \exists \, x \, \vee F G ~ 
      \vee \exists\,x \,F ~ \exists\,x \,G$ & &\\
\end{tabular}\\
\begin{tabular}{llll}
(17)   & $\leftrightarrow ~ \exists \, x \, \to F G ~ 
      \to \forall\,x \,F ~ \exists\,x \,G$ & &\\
\end{tabular}\\
\begin{tabular}{llll}
(18)   & $\leftrightarrow ~ \forall \, x ~ \to F G ~~
        \to ~ \exists \, x \,F ~ G\,,$ & 
        if $x \not\in \mbox{free}(G)$& \\
(19)   & $\leftrightarrow ~ \exists \, x ~ \to F G ~~
        \to ~ \forall \, x \,F ~ G\,,$ & 
        if $x \not\in \mbox{free}(G)$& \\
\end{tabular}\\
\begin{tabular}{llll}
(20)   & $\leftrightarrow ~ \forall \, x ~ J F G ~~
        J ~ F ~ \forall \, x \,G\,,$ & 
        if $x \not\in \mbox{free}(F)$ and $J \in \{\,\to\,;\,\vee\,;\,\& \,\}$& \\
(21)   & $\leftrightarrow ~ \exists \, x ~ J F G ~~
        J ~ F ~ \exists \, x \,G\,,$ & 
        if $x \not\in \mbox{free}(F)$ and $J \in \{\,\to\,;\,\vee\,;\,\& \,\}$\,.& \\
\end{tabular}\\

\underline{Proof:}
In order to check that these formulas are provable in $[M;{\cal L}]$ 
we use former results like
Lemma (3.16) and Theorem (3.17).\dokend\\ 

{\bf (3.19) Proposition}

Let $[M;{\cal L}]$ be a mathematical system with restricted argument lists
and let $F$ be a formula in $[M;{\cal L}]$.
Then for $\mbox{CF}(F;\lambda;x)$ the formulas
\[
\to~\forall\,x\,F~\mbox{SbF}(F\,;\lambda\,;x) \qquad \mbox{and} \qquad 
\to~\mbox{SbF}(F\,;\lambda\,;x)~\exists\,x\,F
\]
are provable in $[M;{\cal L}]$, provided that $\lambda \in {\cal L}$.

\underline{Proof:}~ The formulas $\to~\forall\,x\,F~F$ and
$\to ~ F  ~~ \exists \, x F$ are provable in $M$ due to 
Proposition (3.18). 
Due to Definition (3.7) there hold the conditions
$\mbox{CF}(\forall\,x\,F;\lambda;x)$, $\mbox{CF}(\exists\,x\,F;\lambda;x)$,
$\mbox{CF}(\to~\forall\,x\,F~F;\lambda;x)$ and $\mbox{CF}(\to~F~\exists\,x\,F;\lambda;x)$\,.
The application of the Substitution Rule (c) on $\to~\forall\,x\,F~F$ and 
$\to~F~\exists\,x\,F$ thus gives the proof of Proposition (3.19).\dokend \\

{\bf (3.20) Proposition} (Skolem's normal form)

Let $[M;{\cal L}]$ be a mathematical system with restricted argument lists
and let $F$ be a formula in $[M;{\cal L}]$. Then there are
quantifiers $Q_1$,...,$Q_n$ and variables $x_1$,...,$x_n$ $(n \geq 0)$
as well a formula $G$ in $[M;{\cal L}]$ without quantifiers and 
without the symbols $\leftrightarrow$, $\&$, $\vee$ such that
$$
\leftrightarrow ~ F ~ Q_1 x_1...Q_n x_n G ~\, \in \Pi(M;{\cal L})\,.
$$
\underline{Remark:}~The formula $Q_1 x_1...Q_n x_n G$ 
has Skolem's normal form.

\underline{Proof:}~In the first step we replace $F$ by an equivalent formula
$F'$ such that $F'$ does not contain the symbols
$\leftrightarrow$, $\&$, $\vee$ and such that 
$\leftrightarrow ~ F \, F' ~\, \in \Pi(M;{\cal L})$\,.
This can be done by using the propositional calculus and 
Theorem (3.17)(a) in order to remove subsequently the symbols 
$\leftrightarrow$, $\&$, $\vee$. 

Next we use Theorem (3.17)(b) in order to
construct from $F'$ another formula $F''$ by replacing all bound variables 
in $F'$ by new ones which are not present in $\mbox{free}(F')$ such that
$\leftrightarrow ~ F \, F'' ~\, \in \Pi(M;{\cal L})$\,. 

In the third and last step we use the Proposition (3.18), namely the 
parts (20), (21) for $J = \, \rightarrow$ and (18),(19),(3),(4), and
Theorem (3.17)(a) in order to pull all quantifiers of $F''$
in front of the formula. There finally results the desired formula
$Q_1 x_1...Q_n x_n G$ which has Skolem's normal form. \dokend

\section{The Deduction Theorem and\\ $Z$-homomorphisms}
In this section we first prove the Deduction Theorem, define the $Z$-homo\-mor\-phisms
in a mathematical system and develop the Theorem for $Z$-homo\-mor\-phisms. 
These theorems will be used in order to derive several other results
like the formal proofs by contradiction, the generalization of new constants
and the proofs in mathematical systems with reduced structure.
In the following we fix a formal mathematical system $M = [S; A_M; P_M; B_M]$.

{\bf (4.1)  Definition of statements in M}

A \underline{statement} in $M$ is a formula in $M$ without free variables.

{\bf (4.2)  Extensions of the mathematical system M}
\begin{itemize}
\item[(a)]
Let $\varphi$ be a statement in $M$ and $B_M(\varphi) := B_M \cup \{\varphi\}$.
Then the mathematical system $M(\varphi)$ defined by
$M(\varphi) := [S; A_M; P_M; B_M(\varphi)]$ is called a \underline{simple extension}
of $M$.
\item[(b)]
Let $\Phi$ be a set of statements in $M$ and $B_M(\Phi) := B_M \cup \Phi$.
Then the mathematical system $M(\Phi)$ defined by
$M(\Phi) := [S; A_M; P_M; B_M(\Phi)]$ is called an \underline{extension} of $M$.
\item[(c)]
Let $c$ be a new symbol, which does not occur in $A_M \cup P_M \cup X \cup E$.
Then the \underline{simple symbol-extension} 
$M_c := [S; A_M \cup \{c\}; P_M; B_M]$ of $M$ is also a mathematical system.
\item[(d)]
Let $A \supseteq A_M$ be a set of symbols 
with $z \notin P_M \cup X \cup E$ for all symbols $z \in A$.
Then the \underline{symbol-extension} 
$M_A := [S; A; P_M; B_M]$ of $M$ is also a mathematical system.
\end{itemize}

\underline{Remarks:}
\begin{itemize}
\item[(i)] Note that the extensions of axioms and symbols defined in (4.2)
leave the recursive system $S$ untouched. 
\item[(ii)] $[M(\Phi) ;{\cal L}]$ is a mathematical system with restricted argument
lists if and only if this is the case for $[M;{\cal L}]$
and if in addition the argument lists of all formulas in $\Phi$
are members of ${\cal L}$.
\end{itemize}

\newpage
{\bf (4.3)  The Deduction Theorem, first version}

Let $[M(\varphi) ;{\cal L}]$ be a mathematical system with restricted argument
lists and with a statement $\varphi$. 
Then for every proof $[\Lambda]$ in $[M(\varphi) ;{\cal L}]$ one can construct
a proof $[\Lambda']$ in $[M ;{\cal L}]$ such that $\to \, \varphi \, F ~ \in [\Lambda']$ 
for every $F \in [\Lambda]$. 

\underline{Proof:} We employ induction with respect to the
rules of inference. First we note that for the ``initial proof"
$[\Lambda]=[\,]$ we can also choose $[\Lambda']=[\,]$.

In the following $[\Lambda]$ denotes a proof in $[M(\varphi) ;{\cal L}]$ 
and $[\Lambda']$ a proof in $[M ;{\cal L}]$ such that 
$\to \, \varphi \, F ~ \in [\Lambda']$ 
for every $F \in [\Lambda]$, 
i.e. we assume that the required proof $[\Lambda']$ has already 
been constructed from the proof $[\Lambda]$.
\begin{itemize}
\item[(a)]
Let $F$ be an axiom in $[M(\varphi) ;{\cal L}]$. Then 
the extension $[\Lambda_{*}] = [\Lambda\,;\,F\,]$ is also a proof in 
$[M(\varphi) ;{\cal L}]$
due to Rule (a). If $F=\varphi$, then we put \\
$[\Lambda'_{*}] = [\Lambda'\,;\to \, \varphi \, \varphi]$
for the proof in $[M;{\cal L}]$, 
otherwise $F$ is also an axiom in $[M ;{\cal L}]$, and we put
$
[\Lambda'_{*}] = 
[\Lambda'\,;\,F\,;\,\to ~F ~ \to \, \varphi \, F\,;\,\to \, \varphi \, F\,]
$
for the proof in $[M;{\cal L}]$\,.
\item[(b)]
Let $F$, $G$ be formulas and $F$, $\to ~ F ~ G$ both steps 
of the proof $[\Lambda]$. Then 
$[\Lambda_{*}] = [\Lambda\,;\,G\,]$ is also a proof in $[M(\varphi) ;{\cal L}]$ 
due to Rule (b).
Since $\to \, \varphi \, F ~ \in [\Lambda']$ and 
$\to \, \varphi \, \to F\,G ~ \in [\Lambda']$, we put due to (3.9) and Rule (a),(b)\\
\begin{tabular}{llll}
$[\Lambda'_{*}] = 
[\Lambda'\,;$ & $\to ~ \to \varphi F ~ \to ~ \to \, \varphi \, \to F\,G ~
\to \, \varphi \, G; $&&\\
& $\to ~ \to \, \varphi \, \to F\,G ~\to \, \varphi \, G;$ &&\\
& $\to \, \varphi \, G\,]\,.$ &&\\
\end{tabular}
\item[(c)]
Let $F \in [\Lambda]$, $x \in X$ and $\lambda \in {\cal L}$. 
Suppose that there holds the condition $\mbox{CF}(F;\lambda;x)$.
Then $[\Lambda_{*}] = [\Lambda\,;\,\mbox{SbF}(F;\lambda;x)\,]$ is also a proof 
in $[M(\varphi) ;{\cal L}]$ due to Rule (c). Due to $x \notin \mbox{free}(\varphi)$
there hold the conditions $\mbox{CF}(\to \varphi \, F;\lambda;x)$ and
$\mbox{SbF}(\to \varphi \, F;\lambda;x)~=
~\to \varphi \, \mbox{SbF}(F;\lambda;x)$. 
Since there holds $\to \, \varphi \, F ~ \in [\Lambda']$, we put 
$[\Lambda'_{*}] = [\Lambda'\,;\to \varphi \, \mbox{SbF}(F;\lambda;x)]$.
\item[(d)]
Let $F \in [\Lambda]$ and $x \in X$. Then 
$[\Lambda_{*}] = [\Lambda\,;\,\forall\,x\,F\,]$ is also a proof in $[M(\varphi) ;{\cal L}]$
due to Rule (d), and we put\\
\begin{tabular}{llll}
$[\Lambda'_{*}] = 
[\Lambda'\,;$ & $\forall \, x ~ \to \varphi F\,; $&&\\
& $\rightarrow ~ \forall \, x ~ \to \varphi F ~~
\to ~ \varphi ~ \forall \, x \,F\,;$ &&\\
& $\to ~ \varphi ~ \forall \, x \,F\,]\,.$ &&\\
 \end{tabular}
 
The first new step of the extended proof $[\Lambda'_{*}]$ results from the
assumption $\to \, \varphi \, F ~ \in [\Lambda']$ and Rule (d), the second
step is due to (3.11)(b) and Rule (a) since $\varphi$ has no free variables,
and the third step due to Rule (b).
\item[(e)] 
In the following we fix a predicate symbol $p \in P_S$,
a list $x_1,...,x_i$ of $i\geq 0$ distinct variables
and a formula $G$ in $[M(\varphi);{\cal L}]$. 
Here $x_1,...,x_i$ and the variables of $G$
are not involved in $B_S$.

Then to every R-formula $F$ of $B_S$ there corresponds exactly one formula $F'$ 
of the mathematical system, which is obtained if we replace in $F$ each 
$i-$ary subformula $p \, \lambda_1,...,\lambda_i$, 
where $\lambda_1,...,\lambda_i$ are lists, by the formula 
$G ~ \frac{\lambda_1}{x_1}...\frac{\lambda_i}{x_i}$.
Note that in this case $\lambda_1,...,\lambda_i \in {\cal L}$ is guaranteed.
 
If $F'$ is a step of $[\Lambda]$ for all R-formulas $F \in B_S$ 
for which $p$ occurs $i-$ary in the R-conclusion of $F$, 
then $[\Lambda_*]=[\Lambda; ~ \to ~p\,x_1,...,x_i ~G]$ 
is also a proof in $[M(\varphi) ;{\cal L}]$ due to Rule (e).

(i) First we replace due to Theorem (3.17)(b) the bound variables 
of the statement $\varphi$ subsequently by new ones which are not involved
in $B_S$. There results a proof $[\Lambda'_1]$ in $[M;{\cal L}]$, which is an
extension of $[\Lambda']$ and ends with an equivalence 
$\leftrightarrow \, \varphi \, \psi$, where $\psi$ is a statement in $[M;{\cal L}]$
such that $\mbox{var}(\psi)$ and $\mbox{var}(B_S)$ are disjoint.

(ii) In the next step we consider all R-formulas $F^{(1)},...,F^{(d)} \in B_S$, 
$d \geq 0$, 
for which $p$ occurs $i-$ary in the R-conclusion and note that in this case
${F^{(1)}}',...,{F^{(d)}}'$ are steps of $[\Lambda]$. Let for $1 \leq k \leq d$
the formula ${F^{(k)}_{\psi}}'$ result from $F^{(k)}$ by 
replacing in $F^{(k)}$ each $i-$ary subformula $p \, \lambda_1,...,\lambda_i$, 
where $\lambda_1,...,\lambda_i$ are lists, by the formula 
$\to \, \psi \, G ~ \frac{\lambda_1}{x_1}...\frac{\lambda_i}{x_i}$.
Recall that $\psi$ has no free variables which are available for substitutions.

Then due to the axioms of the propositional calculus the following formulas 
are generally valid
\begin{align}
\to ~ \to \, \psi \, {F^{(1)}}' {F^{(1)}_{\psi}}'\,,~
...~\,,\,
\to ~ \to \, \psi \, {F^{(d)}}' {F^{(d)}_{\psi}}'\,.\nonumber
\end{align}

The formulas 
$\to \, \varphi \, {F^{(1)}}'$,\,...,\,$\to \, \varphi \, {F^{(k)}}' $ 
and the equivalence $\leftrightarrow\,\varphi\,\psi$ are 
steps of the proof $[\Lambda'_1] \supseteq [\Lambda']$, and therefore we can 
use the propositional calculus in order to derive the formulas
${F^{(1)}_{\psi}}'$,...,${F^{(d)}_{\psi}}'$ in a proof 
$[\Lambda'_2]$ in $[M;{\cal L}]$ which is an extension of
$[\Lambda'_1]$.

(iii) Since the variables of the formula 
$\to\,\psi\,G$ are not involved in $B_S$,
we can apply Rule (e) and replace afterwards $\psi$ by $\varphi$
in order to obtain that

\begin{tabular}{llll}
$[\Lambda'_{*}] = 
[\Lambda'_2\,;$ & $\to ~p\,x_1,...,x_i ~\to\,\psi\,G\,;$&&\\
& $\to ~ \to ~p\,x_1,...,x_i ~\to\,\psi\,G $&&\\
& $\to ~\leftrightarrow \, \varphi \, \psi$ &&\\
& $\qquad   \to \varphi ~\to~ p\,x_1,...,x_i ~G;$ &&\\
& $\to ~\leftrightarrow \, \varphi \, \psi$ &&\\
& $\qquad   \to \varphi ~\to~ p\,x_1,...,x_i ~G;$ &&\\
& $\to \varphi ~\to~ p\,x_1,...,x_i ~G\,]\,,$ &&\\
\end{tabular}

is the desired proof in $[M;{\cal L}]$ which satisfies all the required properties.
\end{itemize}
Thus we have proved the first version of the Deduction Theorem.\dokend \\

{\bf (4.4) Corollary, proof by contradiction, first version} 

Let $[M(\neg\,\varphi) ;{\cal L}]$ be a mathematical system with restricted argument
lists and with a statement $\varphi$. 
If $[M(\neg \, \varphi);{\cal L}]$ is contradictory,
i.e. if there is a proof $[\Lambda]$ in $[M(\neg\,\varphi) ;{\cal L}]$ 
which contains a formula $F$ as well as its negation $\neg \, F$, 
then $\varphi \in \Pi(M;{\cal L})$.

\underline{Proof:} ~
Let $[\Lambda]$ be a proof in $[M(\neg\,\varphi) ;{\cal L}]$ which contains a formula $F$
as well as its negation $\neg \, F$. Then

\begin{tabular}{llll}
$[\Lambda_{*}] = 
[\Lambda\,;$ & $\to ~F ~ \to ~\neg\,F ~\varphi~;\,
\to ~\neg\,F ~\varphi~;\,
\varphi\,]$ &&\\
\end{tabular}

is a proof of $\varphi$ in the contradictory system $[M(\neg\,\varphi) ;{\cal L}]$.
From this proof we construct a proof $[\Lambda'_*]$ in $[M ;{\cal L}]$ 
according to the first version of the Deduction Theorem. 
Then $\to ~\neg\,\varphi ~\varphi \in [\Lambda'_{*}]$, 
and we obtain from $[\Lambda'_*]$ the 
extended proof ~ 
$[\Lambda'_*\,;\,\to ~ \to ~\neg\,\varphi ~\varphi~~\varphi\,;\,\varphi\,]$
of $\varphi$ in $[M ;{\cal L}]$\,.\dokend \\

{\bf (4.5) The Deduction Theorem, second version} 

Let $\Phi$ be a set of statements in the mathematical system
$[M;{\cal L}]$ with restricted argument lists. 
For any formula $F$ in $[M;{\cal L}]$ there holds\\ 
~$\to \,\varphi_1\,...\to\,\varphi_m\,F \in \Pi(M;{\cal L})$~ for finitely many
statements $\varphi_1\,,...\,,\,\varphi_m \in \Phi$\\ if and only if
$F \in \Pi(M(\Phi);{\cal L})$.

\underline{Proof:} ~
The ``$\Leftarrow$" direction of the proof is clear, since we can subsequently
apply the Modus Ponens Rule (b) on 
$\to \,\varphi_1\,...\to\,\varphi_m\,F \in \Pi(M;{\cal L}) \subseteq \Pi(M(\Phi);{\cal L})$
for finitely many statements $\varphi_1\,,...\,,\,\varphi_m \in \Phi$
in order to infer $F$ in $\Pi(M(\Phi);{\cal L})$.

For a formula $F$ there holds $F \in \Pi(M(\Phi);{\cal L})$ if and only if 
it is a step of a proof $[\Lambda]$ in $[M(\Phi);{\cal L}]$. 
We define the set $\Gamma = \{\varphi_1\,,...,\,\varphi_m\}$
of all steps in $[\Lambda]$ which are statements in $\Phi$. 
We consider the mathematical systems $M_0 := M$ and
$M_i := M(\{\varphi_1\,,...,\,\varphi_i\})$ for $1 \leq i \leq m$.
Since $F \in \Pi(M_m;{\cal L})$, we conclude from the first version 
of the Deduction Theorem that $\to\,\varphi_m\,F \in \Pi(M_{m-1};{\cal L})$. 
If there is still $m-1>0$, then we infer from $\to\,\varphi_m\,F \in \Pi(M_{m-1};{\cal L})$ that
$\to\,\varphi_{m-1}\,\to\,\varphi_m\,F \in \Pi(M_{m-2};{\cal L})$, using again (4.3), and so on.
After we have applied this procedure $m$-times we conclude that
$\to \,\varphi_1\,...\to\,\varphi_m\,F \in \Pi(M;{\cal L})$.\dokend \\

The theorem for Z-homomorphisms which will be proved in the sequel
is very important in order to obtain proofs
in mathematical systems with certain restricted structure. \\

{\bf (4.6) Definition of a Z-homomorphism}

Let $M = [S; A_M; P_M; B_M]$ be a mathematical system such that
$[M;{\cal L}]$ and $[M;{\cal L}']$ are mathematical systems
with restricted argument lists, and let $Z \subseteq X$ be a (usually finite) subset of variables, 
which may be empty. We consider a mapping $\overline{\,\cdot\,}$
which assigns to each $A_M$-list $\lambda \in {\cal L}$ 
an $A_M$-list $\overline{\lambda} \in {\cal L}'$ such that for all 
$A_M$-lists $\lambda, \mu \in {\cal L}$ the following conditions are satisfied

\begin{tabular}{llll}
(ZH1)   & $\overline{y}=y$ &\mbox{if~} $y \in X$\,,&
\vspace{0.1cm}\\
(ZH2)   & $\overline{\lambda\,\frac{\mu}{x}} = 
          \overline{\lambda}\,\frac{\overline{\mu}}{x}$
&\mbox{if~} $x \in X \setminus Z$\,,& 
\vspace{0.1cm}\\
(ZH3)   & $\mbox{var}(\overline{\mu}) \subseteq \mbox{var}(\mu) \cup Z\,.$&&\\
\end{tabular}

Next we define a natural extension of the mapping $\overline{\,\cdot\,}$ 
to the formulas of $[M;{\cal L}]$. 
Let $F$ be any formula in $[M;{\cal L}]$ such that the variables of $Z$ 
are not occurring bound in $F$, i.e. $F$ does not contain a subformula 
of the form $Q\,z\,F'$, where $Q \in \{\forall, \exists \}$ and $z \in Z$.
For abbreviation we denote the set of all these formulas $F$ by $\Phi_{M;{\cal L}}^Z$.
We replace in $F \in \Phi_{M;{\cal L}}^Z$ all the argument lists $\lambda$ by 
$\overline{\lambda}$\,. 
There results a formula $\overline{F}$ of $[M;{\cal L}']$.

In the following we suppose in addition that there hold the two conditions

\begin{tabular}{llll}
(ZH4)   & $\overline{F} = F$ &for all $F \in B_M \cap \Phi_{M;{\cal L}}^Z$\,.&\\
(ZH5)   & $Z \cap \mbox{var}(B_S)=\{\,\}\,,$ 
&i.e. the variables of $Z$ are not involved in $B_S$\,.&\\
\end{tabular}

Then the extended mapping $\overline{\,\cdot\,}$
is called a $Z$-homomorphism from $[M;{\cal L}]$ in $[M;{\cal L}']$.
If ${\cal L}={\cal L}'$, then $\overline{\,\cdot\,}$
is just called a $Z$-homomorphism in $[M;{\cal L}]$.
Note that especially $F \in \Phi_{M;{\cal L}}^Z$ for all $F \in B_S$.\\

{\bf (4.7) Lemma}

We consider the mapping $\overline{\,\cdot\,}$ from Definition (4.6), 
which satisfies the conditions (ZH1)-(ZH3), 
and its extension to the formulas $F \in \Phi_{M;{\cal L}}^Z$.
Let $F \in \Phi_{M;{\cal L}}^Z$.
Then for every list $\mu \in {\cal L}$ and for all variables $x \in X \setminus Z$ with 
$CF(F; \mu; x)$ there holds the condition $\mbox{CF}(\overline{F};\,\overline{\mu};\,x)$ 
and the equation
\[
\overline{F\,\frac{\mu}{x}}=\overline{F}\,\frac{\overline{\mu}}{x}\,.
\]

\underline{Proof:}

We use induction with respect to the rules for generating 
formulas in $[M;{\cal L}]$.
The variable $x \in X \setminus Z$ and the list $\mu \in {\cal L}$ are arbitrary,
but will be fixed in the following.
For any formula $F$ in $[M;{\cal L}]$ we define the

Condition $(*)$:

If $F \in \Phi_{M;{\cal L}}^Z$ and if $\mbox{CF}(F;\,\mu;\,x)$,
then there holds the condition $\mbox{CF}(\overline{F};\,\overline{\mu};\,x)$ 
and the equation
$\overline{F\,\frac{\mu}{x}}=\overline{F}\,\frac{\overline{\mu}}{x}$\,.

We prove that Condition $(*)$ is satisfied for all formulas $F$ in $[M;{\cal L}]$. 
We use the definitions (3.6) and (3.7) and the notations occurring there 
by treating the corresponding cases (a)-(d). 
\begin{itemize}
\item[(a)] If $F$ is a prime formula in $[M;{\cal L}]$, 
then $\overline{F}$ is a prime formula in $[M;{\cal L}']$.
In this case we obtain $\mbox{CF}(F;\,\lambda;\,x)$
as well as $\mbox{CF}(\overline{F};\,\overline{\lambda};\,x)$. 
We define for $p \in P_M$ and $\lambda_1\,,\,\lambda_2,... \in {\cal L}$

$F_1 \,= ~\sim \lambda_1\,,\,\lambda_2~,~
F_2 \,= p~,~
F_3 \,= p\, \lambda_1~,~
F_4 \,= p\, \lambda_1,...,\lambda_i\,,$

and can apply (ZH2) due to $x \in X \setminus Z$ to obtain 

$\overline{F_1\,\frac{\mu}{x}} \,= 
~\sim \overline{\lambda_1\,\frac{\mu}{x}}\,,\,
      \overline{\lambda_2\,\frac{\mu}{x}}\,=
~\sim \overline{\lambda_1}\,\frac{\overline{\mu}}{x}\,,\,
      \overline{\lambda_2}\,\frac{\overline{\mu}}{x}\,=
\overline{F_1}\,\frac{\overline{\mu}}{x}$

and

$\overline{F_2\,\frac{\mu}{x}} \,= p\,=\overline{F_2}\,\frac{\overline{\mu}}{x}~,~
\overline{F_3\,\frac{\mu}{x}} \,= p\, \overline{\lambda_1\,\frac{\mu}{x}}\,=
p\, \overline{\lambda_1}\,\frac{\overline{\mu}}{x}\,=
\overline{F_3}\,\frac{\overline{\mu}}{x}$~,\\ 
$\overline{F_4\,\frac{\mu}{x}} \,= 
p\, \overline{\lambda_1\,\frac{\mu}{x}},...,\overline{\lambda_i\,\frac{\mu}{x}}\,=
p\, \overline{\lambda_1}\,\frac{\overline{\mu}}{x},...,
    \overline{\lambda_i}\,\frac{\overline{\mu}}{x}\,=
\overline{F_4}\,\frac{\overline{\mu}}{x}$\,.

We have thus confirmed Condition $(*)$ for the prime formulas.
\item[(b)] We assume that Condition $(*)$ is satisfied for a $M$-formula $F$,
that $\neg\,F \in \Phi_{M;{\cal L}}^Z$ and that there holds the condition 
$\mbox{CF}(\neg \, F; \mu; x)$. Then $F \in \Phi_{M;{\cal L}}^Z$, and there holds
the condition $\mbox{CF}(F; \mu; x)$. Since $F$ satisfies Condition $(*)$,
we conclude that $\mbox{CF}(\overline{F}; \overline{\mu}; x)$ and
$\mbox{CF}(\overline{\neg \, F}; \overline{\mu}; x)$ are valid and that
the equations
\[
\overline{\mbox{SbF}(\neg \, F; \mu; x)} =
\overline{\neg \, F\,\frac{\mu}{x}} = 
\neg\,\overline{F\,\frac{\mu}{x}} =
\neg\,\overline{F}\,\frac{\overline{\mu}}{x} =
\mbox{SbF}(\overline{\neg \, F}; \overline{\mu}; x)
\]
are satisfied. Thus we have confirmed Condition $(*)$ for $\neg \, F$.
\item[(c)] We assume that Condition $(*)$ is satisfied for the
$M$-formulas $F,G$,
that $J \, F G \in \Phi_{M;{\cal L}}^Z$ and that $\mbox{CF}(J \, F G; \mu; x)$ holds. 
We obtain $F \in \Phi_{M;{\cal L}}^Z$ and $G \in \Phi_{M;{\cal L}}^Z$,
and there hold $\mbox{CF}(F; \mu; x)$ and $\mbox{CF}(G; \mu; x)$. 
Since $F$ and $G$ satisfy Condition $(*)$,
we conclude that $\mbox{CF}(\overline{F}; \overline{\mu}; x)$ and 
$\mbox{CF}(\overline{G}; \overline{\mu}; x)$ are both valid.
Therefore $\mbox{CF}(J\,\overline{F}\,\overline{G}; \overline{\mu}; x)$, 
which is equivalent to
$\mbox{CF}(\overline{J\,F G}; \overline{\mu}; x)$, is also satisfied.
Since $F$ and $G$ satisfy Condition $(*)$, we obtain
\begin{align}
\overline{\mbox{SbF}(J \, F G; \mu; x)} =
\overline{J\,F\,\frac{\mu}{x} ~ G\,\frac{\mu}{x}} =
J\,\overline{F\,\frac{\mu}{x}} ~\, \overline{G\,\frac{\mu}{x}}\nonumber \\ = 
J\,\overline{F}\,\frac{\overline{\mu}}{x} \,\overline{G}\,\frac{\overline{\mu}}{x} = 
\mbox{SbF}(\overline{J \, F G}; \overline{\mu}; x)\,,\nonumber
\end{align}
i.e. Condition $(*)$ is satisfied for $J \, F G$.
\item[(d)]
We assume that $(*)$ is satisfied for an $M$-formula $F$,
that moreover $Q\,y\,F \in \Phi_{M;{\cal L}}^Z$
and that there holds $\mbox{CF}(Q\,y\,F; \mu; x)$.
It follows that $y \notin Z$, since $y$ is bound in $Q\,y\,F$.
Note that $\mbox{free}(\overline{F}) \subseteq \mbox{free}(F) \cup Z$.

If $x \notin \mbox{free}(F)\setminus \{y\}$, then we obtain 
$\mbox{CF}(\overline{Q\,y\,F};\,\overline{\mu};\,x)$ with $\overline{Q\,y\,F}=Q\,y\, \overline{F}$\\
and $\overline{\mbox{SbF}(Q\,y\,F;\,\mu;\,x)}=
\overline{Q\,y\,F}=
\mbox{SbF}(\overline{Q\,y\,F};\,\overline{\mu};\,x)$\,.

Otherwise we use that
$\mbox{CF}(Q\,y\,F; \mu; x)$ is satisfied with $x \neq y$ in order to conclude that
$y \notin \mbox{var}(\mu)$ and 
$y \notin \mbox{var}(\overline{\mu})\subseteq \mbox{var}(\mu) \cup Z$
due to $y \notin Z$ and that $\mbox{CF}(F; \mu; x)$.
But $F$ satisfies the Condition $(*)$ and $Q\,y\,F \in \Phi_{M;{\cal L}}^{Z}$,
and therefore $\mbox{CF}(\overline{F}; \overline{\mu}; x)$. 
From $y \notin \mbox{var}(\overline{\mu})$
and $\mbox{CF}(\overline{F}; \overline{\mu}; x)$ we conclude that 
$\mbox{CF}(Q\,y\,\overline{F}; \overline{\mu}; x)$, i.e.
$\mbox{CF}(\overline{Q\,y\,F}; \overline{\mu}; x)$ is again satisfied.
Since $F$ satisfies the Condition $(*)$, we finally conclude due to $x \neq y$ that
\[
\overline{\mbox{SbF}(Q \, y \, F ; \mu ; x)}=
Q \, y \, \mbox{SbF}(\overline{F} ; \overline{\mu} ; x)=
\mbox{SbF}(\overline{Q \, y \, F} ; \overline{\mu} ; x)\,,
\]
i.e. Condition $(*)$ is satisfied for $Q \, y \, F$.
\end{itemize}
Thus we have proved Lemma (4.7).\dokend \\

{\bf (4.8) Theorem for $Z$-homomorphisms, Z-Theorem}

We consider a $Z$-homomorphism $\overline{\,\cdot\,}$ from 
$[M;{\cal L}]$ in $[M;{\cal L}']$ with the assumptions given in (4.6).
Suppose that $[\Lambda]=[F_1;...;F_l]$ is a proof in $[M;{\cal L}]$ 
and that the variables of $Z$ are not involved in $F_1,...,F_l$.
Then we conclude that $F_1,...,F_l \in \Phi_{M;{\cal L}}^Z$, and 
$[\overline{\Lambda}]=[\overline{F}_1;...;\overline{F}_l]$ is 
a proof in $[M;{\cal L}']$.


\underline{Proof:} We employ induction with respect to the
rules of inference. First we note that for the ``initial proof"
$[\Lambda]=[\,]$ we can also choose $[\overline{\Lambda}]=[\,]$.

In the following we assume that $[\Lambda]$ is a proof in $[M;{\cal L}]$,
that the variables of $Z$ are not involved in $[\Lambda]$
and that 
$[\overline{\Lambda}]=[\overline{F}_1;...;\overline{F}_l]$ 
is a proof in $[M;{\cal L}']$.

\begin{itemize}
\item[(a)]
Let $H$ be an axiom in $[M;{\cal L}]$ which does not contain
any $z \in Z$. Then $[\Lambda_{*}] = [\Lambda\,;\,H\,]$ is also a proof in $[M;{\cal L}]$ 
due to Rule (a). We note that $H \in \mlz$. 
Therefore it is sufficient to show that $\overline{H}$ is an axiom 
in $\mlp$. For this purpose we distinguish four cases.

1.) ~ Let $\alpha=\alpha(\xi_1,...,\xi_j)$ be an identically true
propositional function of the distinct propositional variables $\xi_1,...,\xi_j$, $j \geq 1$.
We suppose without loss of generality that all $j$ propositional variables occur in $\alpha$.
If $H_1$,...,$H_j$ are any formulas in $\ml$ with $H = \alpha(H_1,...,H_j)$,
then the variables of $Z$ are not involved in $H_1$,...,$H_j$, and
$\overline{H} = \alpha(\overline{H_1},...,\overline{H_j})$ 
is an axiom of the propositional calculus in $\mlp$. 
Therefore $[\overline{\Lambda_{*}}] = [\overline{\Lambda}\,;\,\overline{H}\,]$ 
is a proof in $\mlp$ due to Rule (a).

2.) ~If $H$ is an axiom of equality in $\ml$ according to (3.10)(a,c), then 
$\overline{H}=H$ due to (ZH1), i.e. $\overline{H}$
is also an axiom of equality in $\mlp$.
If $H = ~ \to ~ \mbox{SbF}(\,\sim \lambda,\mu\,;\,x\,;\,y\,) ~ 
         \to ~ \sim x,y ~ 
         \sim \lambda,\mu$ 
is an axiom of equality in $\ml$ according to (3.10)(b), 
then $\overline{H}$ is an axiom of equality in $\mlp$ of the form (3.10)(b)
due to (ZH2), since the variables of $Z$ are not involved in $H$.

3.) ~The quantifier axioms (3.11) can be handled very easily since we suppose
that $Z$ and $\mbox{var}(F) \cup \mbox{var}(G) \cup \{x\}$ are disjoint.

4.) ~For $H \in B_M$ we obtain $H \in B_M \cap \mlz$ from $\mbox{var}(H)\cap Z=\{\}$,
and therefore $\overline{H} \in B_M$ due to (ZH4). Then
$[\overline{\Lambda_{*}}] = [\overline{\Lambda}\,;\,\overline{H}\,]$ 
is a proof in $\mlp$ due to Rule (a).
\item[(b)]
Let $F$, $G$ be two formulas in $\ml$ and $F$, $\to F\,G$ both steps of the proof 
$[\Lambda]$. Then 
$[\Lambda_{*}] = [\Lambda\,;\,G\,]$ is also a proof in $\ml$ due to Rule (b),
which does not contain a variable $z \in Z$.
It follows that $\overline{F}$ and 
$\overline{\to F\,G}\,=\,\to \overline{F}\,\overline{G}$
are both steps of the proof $[\overline{\Lambda}]$ due to our assumptions,
and due to Rule (b) we put 
$[\overline{\Lambda}_{*}] = 
[\overline{\Lambda}\,;\overline{G}\,]$
for the required proof in $\mlp$. 
\item[(c)]
Let $F \in [\Lambda]$, $x \in X$ and $\lambda \in {\cal L}$.
Suppose that there holds the condition $\mbox{CF}(F;\lambda;x)$.
Then $[\Lambda_{*}] = [\Lambda\,;\,F\,\frac{\lambda}{x}\,]$ is also a proof 
in $\ml$ due to Rule (c). We suppose that $x \in \mbox{free}(F)$ without loss of generality. 
Then the condition that  $[\Lambda_{*}]$ does not contain any variable in $Z$
is equivalent to $z \notin \mbox{var}(\lambda)$ for all $z \in Z$, which will be assumed
here. Note that $F \in \mlz$ due to $F \in [\Lambda]$ and 
$z \notin \mbox{var}([\Lambda])$ for all $z \in Z$. 
Moreover, we know that $x \in X \setminus Z$,
since $x \in \mbox{free}(F)$ occurs in $[\Lambda]$ due to $F \in [\Lambda]$.
Therefore we obtain due to Lemma (4.7) that there holds the condition
$\mbox{CF}(\overline{F};\,\overline{\lambda};\,x)$ and the equation 
$\overline{F\,\frac{\lambda}{x}}=\overline{F}\,\frac{\overline{\lambda}}{x}$.
Since $\overline{F} \in [\overline{\Lambda}]$ we conclude that
$[\overline{\Lambda}_{*}] = 
[\overline{\Lambda}\,;\overline{F\,\frac{\lambda}{x}}\,]$
is a proof in $\mlp$ due to Rule (c).
\item[(d)]
Let $F \in [\Lambda]$ and $x \in X$. Then 
$[\Lambda_{*}] = [\Lambda\,;\,\forall\,x\,F\,]$ is also a proof in $\ml$
due to Rule (d). The condition that the variables of $Z$ are not involved in 
$[\Lambda_{*}]$ is equivalent to $x \notin Z$, which will be assumed here. 
Since $F \in [\Lambda]$ implies $\overline{F} \in [\overline{\Lambda}]$ and since 
$\overline{\forall\,x\,F}=\forall\,x\,\overline{F}$, we can apply Rule (d) on 
$[\overline{\Lambda}]$, $\overline{F}$ in order to conclude that
$[\overline{\Lambda}_{*}] = 
[\overline{\Lambda}\,;\overline{\forall\,x\,F}\,]$
is a proof in $\mlp$.  
\item[(e)] 
In the following we fix a predicate symbol $p \in P_S$,
a list $x_1,...,x_i$ of $i\geq 0$ distinct variables
and a formula $G$ in $\ml$. 
We suppose that $x_1,...,x_i$ and the variables of $G$
are not involved in $B_S$.

Then to every R-formula $F$ of $B_S$ there corresponds exactly one formula $F'$ 
of the mathematical system, which is obtained if we replace in $F$ each 
$i-$ary subformula $p \, \lambda_1,...,\lambda_i$, 
where $\lambda_1,...,\lambda_i$ are lists, by the formula 
$G ~ \frac{\lambda_1}{x_1}...\frac{\lambda_i}{x_i}$.
Note that in this case $\lambda_1,...,\lambda_i \in {\cal L}$
due to (ZH4).
 
If $F'$ is a step of $[\Lambda]$ for all R-formulas $F \in B_S$ 
for which $p$ occurs $i-$ary in the R-conclusion of $F$, 
then $[\Lambda_*]=[\Lambda; ~\to ~p\,x_1,...,x_i ~G]$ 
is also a proof in $\ml$ due to Rule (e).

The condition that the variables of $Z$ are not involved in $[\Lambda_{*}]$ 
implies that $z \notin \{\,x_1,...,x_i\,\} \cup \, \mbox{var}(G)$ for all $z \in Z$,
which will be assumed here. 

To every R-formula $F$ of $B_S$ there corresponds the formula $F''$, 
which is obtained if we replace in $F$ each 
$i-$ary subformula $p \, \lambda_1,...,\lambda_i$, 
where $\lambda_1,...,\lambda_i$ are lists, by the formula 
$\overline{G\, \frac{\lambda_1}{x_1}...\frac{\lambda_i}{x_i}}$.
Due to our assumption that $\overline{F}=F$ for all $F \in B_S$
it follows that $\overline{\lambda}=\lambda \in {\cal L}\cap {\cal L}'$ 
for all argument lists $\lambda$
which occur in the formulas of $B_S$. Since 
the variables of $Z$ are not occurring among the bound variables in $G$,
since $x_1,...,x_i \in X \setminus Z$ and since the variables in 
$\lambda_1,...,\lambda_i$ are not occurring among the bound variables in $G$, 
we can i-times apply Lemma (4.7) in order to conclude that
\[
\overline{G\, \frac{\lambda_1}{x_1}...\frac{\lambda_i}{x_i}}=
\overline{G}\, 
\frac{\overline{\lambda_1}}{x_1}...
\frac{\overline{\lambda_i}}{x_i}=
\overline{G}\, 
\frac{\lambda_1}{x_1}...
\frac{\lambda_i}{x_i}\,.
\]
But $F''=\overline{F'}$, and $F''$ is a step of $[\overline{\Lambda}]$ 
for all R-formula $F$ of $B_S$ for which $p$ occurs $i-$ary in the R-conclusion of $F$.
Moreover, the variables of $\overline{G}$ are
not involved in $B_S$ due to $\mbox{var}(G) \cap \mbox{var}(B_S) = \{\,\}$ 
and (ZH3), (ZH5). Thus we can apply Rule (e) on $[\overline{\Lambda}]$ 
and conclude that
\[[\overline{\Lambda}_{*}] = 
[\overline{\Lambda}\,;\overline{~\to ~p\,x_1,...,x_i ~G}\,]\]
is a proof in $\mlp$. 
\end{itemize}
Thus we have proved the Theorem for $Z$-homomorphisms.\dokend \\

Often in mathematical arguments we say ``let $n$ be an arbitrary but fixed integer".
Then we proceed with a proof and come to a certain conclusion $A(n)$.
We can then deduce that $A(n)$ is valid for all integers $n$, since we have
not used special properties of $n$. The next Corollaries show that these 
argumentations can also be done formally in a mathematical system.\\ 

{\bf (4.9) Corollary, generalization of new constants in symbol-extensions}

Let $[M;{\cal L}]$ with $M = [S; A_M; P_M; B_M]$ be a mathematical system
with restricted argument lists. We consider a symbol-extension
$M_A = [S; A; P_M; B_M]$ of $M$ with $A \supseteq A_M$. 
\begin{itemize}
\item[(a)] If the set ${\cal L}_A$ of argument lists in $M_A$ is defined by
\[ {\cal L}_A := \{\, \lambda \frac{c_1}{x_1}...\frac{c_m}{x_m} \,\, | \,\, 
\lambda \in {\cal L},\, x_1,...,x_m \in X,\, c_1,...,c_m \in 
A\setminus A_M \,,\, m \geq 0\}\,,
\]
then $[M_A;{\cal L}_A]$ is a mathematical system with restricted argument lists.
\item[(b)] Suppose that $x_1,...,x_m \in X$ are $m \geq 0$ distinct variables
and that $c_1,...,c_m \in A\setminus A_M$ are $m$ distinct new constants. 
If $F$ is a formula in $[M;{\cal L}]$ such that 
$F\,\frac{c_1}{x_1}...\frac{c_m}{x_m}\in \Pi(M_A;{\cal L}_A)$, 
then $F \in \Pi(M;{\cal L})$ as well as 
$\forall\,x_1...\forall\,x_m\,F \in \Pi(M;{\cal L})$.
\end{itemize}
\underline{Proof:} ~ (a) Choosing $m=0$ we first note that 
${\cal L}_A \supseteq {\cal L}$ is an extension of ${\cal L}$, 
and hence $[M_A;{\cal L}_A]$ to be constructed satisfies (3.15)(i) and (iii).
Note that $x\,\frac{c}{x}=c \in {\cal L}_A$ for any $x \in X$ and all
$c \in A \setminus A_M$\,.
It remains to prove the substitution invariance for ${\cal L}_A$.
Let $\lambda, \mu \in {\cal L}_A$ and $x \in \mbox{var}(\lambda)$.
Let $d_1,...,d_n \in A \setminus A_M$ for $n \geq 0$ be a complete
list of all new symbols occurring in $\lambda$ and $\mu$ and
let $d_1,...,d_n$ be distinct. Choose distinct variables
$y_1,...,y_n \in X$ which are neither occurring in $\lambda$ nor in $\mu$\,.
Since $d_1,...,d_n$ occur only as constant symbols in $\lambda$ and $\mu$,
we can replace them by $y_1,...,y_n$ in order to obtain new
lists $\lambda', \mu' \in {\cal L}$ due to the properties of 
${\cal L}$ and ${\cal L}_A$. We obtain $\lambda'\,\frac{\mu'}{x} \in {\cal L}$ and 
\[
\lambda\,\frac{\mu}{x} = 
\lambda'\,\frac{\mu'}{x}\,\frac{d_1}{y_1}...\frac{d_n}{y_n} \in {\cal L}_A\,.
\]
(b) ~ Suppose without loss of generality that $x_1,...,x_m \in \mbox{free}(F)$.
Let $[\Lambda]$ be a proof of $F\,\frac{c_1}{x_1}...\frac{c_m}{x_m}$
in $[M_A;{\cal L}_A]$ and let $d_1,...,d_n \in A \setminus A_M$ with $n \geq m$
be all distinct new constants occurring in $[\Lambda]$. Choose a set
$Z:=\{z_1,...,z_n \}$ of $n$ distinct variables, which are neither occurring in $[\Lambda]$
nor in $B_S$ and which are distinct from $\mbox{var}(F)$. 
Due to (a) we can define a Z-homomorphism in $[M_A;{\cal L}_A]$
by replacing for $k \leq n$ each occurrence of a new constant $d_k$ as a sublist
in an argument list $\lambda \in {\cal L}_A$ by the variable $z_k$. 
It follows from Theorem (4.8) that $[\overline{\Lambda}]$ is a proof in $[M_A;{\cal L}_A]$ 
which has only formulas with argument lists in ${\cal L}$ and which contains 
the step $F\,\frac{z_{k_1}}{x_1}...\frac{z_{k_m}}{x_m}$\,,
where $z_{k_1},...,z_{k_m} \in Z$ correspond to the new constants
$c_1,...,c_m$, respectively.
Hence we obtain that $[\overline{\Lambda}]$ is already
a proof in $[M;{\cal L}]$ and that 
$F\,\frac{z_{k_1}}{x_1}...\frac{z_{k_m}}{x_m} \in \Pi(M;{\cal L})$\,.
Since $z_{k_1},...,z_{k_m}$ are distinct, 
we can subsequently apply Lemma (3.16)(a) 
and the Substitution Rule on the last formula in order to conclude that $F$ 
and hence $\forall\,x_1...\forall\,x_m\,F$
are provable in $[M;{\cal L}]$. \dokend

{\bf (4.10) Corollary,~proof by contradiction, second version}

Let $[M;{\cal L}]$ with $M = [S; A_M; P_M; B_M]$ be a mathematical system
with restricted argument lists. We consider a symbol-extension
$M_A = [S; A; P_M; B_M]$ of $M$ with $A \supseteq A_M$. 
Define $[M_A;{\cal L}_A]$ as in Corollary (4.9) and suppose that

 \begin{tabular}{llll}  
  (i)  & $c_1,...,c_m \in A \setminus A_M$ &are $m \geq 0$ distinct constants,&\\
  (ii) & $x_1,...,x_m \in X$               &are $m$ distinct variables,&\\
  (iii) & $F$ &is a formula in $[M;{\cal L}]\,,$ &\\
  (iv)  & $F\,\frac{c_1}{x_1}...\frac{c_m}{x_m}$  &is a statement in $[M_A;{\cal L}_A]$\,,&\\
  (v)   & $[M_A(\neg\,F\,\frac{c_1}{x_1}...\frac{c_m}{x_m});{\cal L}_A]$ & is contradictory\,.&\\      
  \end{tabular}

Then $F$ and the statement $\forall\,x_1...\forall\,x_m\,F$ 
are both provable in $\ml$.

\underline{Proof:} Due to Corollary (4.4) we know that the statement
$F\,\frac{c_1}{x_1}...\frac{c_m}{x_m}$ is provable in $[M_A;{\cal L}_A]$, 
and due to Corollary (4.9) we conclude that the formula $F$ 
as well as the statement $\forall\,x_1...\forall\,x_m\,F$ 
are provable in $[M;{\cal L}]$.\dokend 

In the following we consider Z-homomorphisms from
a mathematical system $M=[S;A_M;P_M;B_M]$ without restrictions 
of the argument lists, i.e. formally we can put for ${\cal L}$ 
the set of all $A_M$-lists, to a mathe\-matical systems $[M;{\cal L}']$
with restricted argument lists in ${\cal L}'$.

{\bf (4.11) Corollary,~restriction to special argument lists}

We consider a mathematical system $M = [S; A_M; P_M; B_M]$.
\begin{itemize}
\item[(a)] Let ${\cal L}'=(A_M \cup X)^+$ be the set of all nonempty strings 
with respect to the set $A_M \cup X$. Suppose that $B_M$ has only argument lists
in ${\cal L}'$. Then we have a Z-homomorphism $\overline{\,\cdot\,}$
from $M$ in $[M;{\cal L}']$ erasing operation terms
with $\overline{F}=F$ if $F \in \Phi_{M;{\cal L'}}^Z$ has argument lists 
in ${\cal L}'$.
\item[(b)] Let ${\cal L}'=A_M \cup X$ be the set of all variables
and $A_M$-constants. Suppose that $B_M$ has only argument lists
in ${\cal L}'$. Then one can construct a Z-homomorphism $\overline{\,\cdot\,}$
from $M$ in $[M;{\cal L}']$ erasing all argument lists
which are neither a constant nor a variable symbol
such that $\overline{F}=F$ for any formula $F \in \Phi_{M;{\cal L'}}^Z$ 
with argument lists in ${\cal L}'$.
\item[(c)] Let ${\cal L}'= X$ be the set of all variables. 
Suppose that $B_M$ has only argument lists in ${\cal L}'$. 
Then we can construct a Z-homomorphism $\overline{\,\cdot\,}$
from $M$ in $[M;{\cal L}']$ erasing all non-variable argument lists
such that $\overline{F}=F$ for any formula $F \in \Phi_{M;{\cal L'}}^Z$ with argument lists 
in ${\cal L}'$.
\end{itemize}

\underline{Remark:}~
It follows from this Corollary that in $[M;{\cal L}']$
we can prove all formulas which have only argument lists in ${\cal L}'$
and which are provable in the original mathematical system $M$
without restrictions of the argument lists.

\underline{Proof:} ~ For all three cases we define a mapping $\overline{\,\cdot\,}$ 
which assigns to each $A_M$-list $\lambda$ 
an $A_M$-list $\overline{\lambda} \in {\cal L}'$ such that (ZH1)-(ZH3)
are satisfied. The extension of these mappings to the formulas
$F \in \Phi_{M}^{Z}$ due to Definition (4.6) defines 
the desired Z-homomorphisms from $M$ in $[M;{\cal L}']$ in all three cases. 
This is possible since we take into consideration that $Z \cap \mbox{var}(B_S) = \{ \}$ and 
since $B_M$ has only argument lists in ${\cal L}'$. 
\begin{itemize}
\item[(a)] For any list $\lambda$ in $M$ we replace all the \underline{maximal} 
$a$-subterms in $\lambda$ of the form $a(\mu)$, 
$\mu$ is a list in $M$ and $a \in A_M$, by a 
variable $\delta(a)$ with  $\delta(a) \in X \setminus \mbox{var}(B_S)$.
Note that $\delta$ need not be injective and put $Z=\delta(A_M)$.
There results a list $\overline{\lambda} \in {\cal L}'$ without operation terms, 
and the corresponding mapping $\overline{\,\cdot\,}$ 
can be extended to an Z-homomorphism from $M$ in $[M;{\cal L}']$.
\item[(b)] We put $Z = \{z\}$ for a fixed variable
 $z \in X \setminus \mbox{var}(B_S)$ and define
for any list $\lambda$ in $M$ 
\[
   \overline{\lambda} = \left\{
\begin{array}{r@{\quad,\quad}l}
a &  \lambda = a \in A_M\\
x &  \lambda = x \in X\\ 
z & \mbox{otherwise}\,.
\end{array}
\right.
\]
\item[(c)]We put $Z = \{z\}$ for a fixed variable 
$z \in X \setminus \mbox{var}(B_S)$ and define
for any list $\lambda$ in $M$ 
\[
   \overline{\lambda} = \left\{
\begin{array}{r@{\quad,\quad}l}
x &  \lambda = x \in X\\ 
z & \mbox{otherwise}\,.
\end{array}
\right.
\]
\end{itemize}
Thus we have shown Corollary (4.11). \dokend

\section{Consistency and incompleteness}

Using the Deduction Theorem derived in the last section we have reduced
the question concerning the provability of formulas in an arbitrary 
mathematical system $M$ to the case that $B_M=B_S$. The first simple result 
shows that in these special mathematical systems there cannot appear 
a contradiction. \\

{\bf (5.1) Proposition}

Let $M=[S;\,A_M;\,P_M\,;B_M]$ be a mathematical system with $B_M=B_S$.
Then $M$ is not contradictory, i.e. there is no proof $[\Lambda]$ in $M$
which contains a formula $F$ as well as its negation $\neg\,F$.

\underline{Proof:}

1.) Let $\Gamma$ be a finite set of R-formulas, $p \in P_S$ and $i \geq 0$ an integer
number. We say that the pair $(p,i)$ \underline{fails in $\Gamma$}, if there is no i-ary
R-conclusion $p\,\lambda_1,...,\lambda_i$ in the formulas of $\Gamma$.
Recall that $p\,\lambda_1,...,\lambda_i=p$ for $i=0$.

2.) An R-formula $F \in \Gamma$ is called \underline{spare in $\Gamma$}, if there
is a $p \in P_S$ and an integer number $i \geq 0$ such that an i-ary prime R-formula
$p\,\lambda_1,...,\lambda_i$ occurs as an R-subformula in $F$ and such that $(p,i)$
fails in $\Gamma$. Let $\Gamma' \subseteq \Gamma$ result from $\Gamma$ by
cancelling all the formulas $F \in \Gamma$ which are spare in $\Gamma$.

3.) Let $B_S^{(0)}$ result from $B_S$ by cancelling all the formulas $F \in B_S$
for which there are two i-ary prime R-formulas $p\,\lambda_1,...,\lambda_i$ and
$p\,\lambda_1',...,\lambda_i'$ with the same predicate symbol $p \in P_S$
such that $p\,\lambda_1,...,\lambda_i$ is the R-conclusion of $F$ and 
$p\,\lambda_1',...,\lambda_i'$ an R-premise of $F$. Then we define 
$B_S^{(k+1)} = {B_S^{(k)}}'$ for all integer numbers $k \geq 0$. Since

\[
B_S^{(0)} \supseteq B_S^{(1)} \supseteq B_S^{(2)} \supseteq B_S^{(3)} \supseteq~...
\]

and since $B_S^{(0)}$ is finite, there is a minimal index $k_0 \geq 0$
such that

\[
B_S^{(k_0)} = B_S^{(k_0+1)} = B_S^{(k_0+2)} = B_S^{(k_0+3)} = ~...\,.
\]

4.) Let $\mbox{Prime}\,(p,i)$ for $(p,i) \in P_S \times \N_0$ be the set
of all i-ary prime R-formulas $p\,\lambda_1,...,\lambda_i$ and define
$\chi \, : \bigcup\limits_{(p,i) \in P_S \times \N_0}\mbox{Prime}\,(p,i) \to \{-1,+1\}$ 
by

$
\chi(p\,\lambda_1,...,\lambda_i)=\left\{
\begin{array}{r@{\quad,\quad}l}
+1 &  \mbox{if}~p~\mbox{occurs i-ary in}~B_S^{(k_0)}\\ -1 &  \mbox{otherwise}\,.
\end{array}
\right.
$

Moreover we put $\chi(\sim\lambda_1,\lambda_2)=1$ for all lists $\lambda_1$, $\lambda_2$
and $\chi(F)=-1$ for all prime formulas with a predicate symbol $p \in P_M \setminus P_S$.
Thus $\chi$ defines a sign for all prime formulas in $M$. 

5.) Let $F$, $G$ be formulas in $M$ 
for which $\chi(F)$ and $\chi(G)$ are already declared. Then we put
for $x \in X$ and $Q \in \{ \forall\,,\,\exists\}$

\begin{tabular}{llll}
(i)  &$\chi(\neg\,F)=-\chi(F)$\,, &&\\ \vspace{0.2cm}
(ii) &$\chi(\to\,F G)=\left\{
\begin{array}{r@{\quad,\quad}l}
+1 &  \mbox{if}~ \chi(F)=-1 ~\mbox{or}~ \chi(G)=1 \\ -1 &  \mbox{otherwise}\,,
\end{array}
\right.$ &&\\
\end{tabular}\\
\vspace{0.2cm}
\begin{tabular}{llll}
(iii) &$\chi(\leftrightarrow\,F G)=\left\{
\begin{array}{r@{\quad,\quad}l}
+1 &  \mbox{if}~ \chi(F) = \chi(G)\\ 
-1 &  \mbox{otherwise}\,,
\end{array}
\right.$ &&\\
\end{tabular}\\
\vspace{0.2cm}
\begin{tabular}{llll}
(iv) &$\chi(\&\,F G)=\left\{
\begin{array}{r@{\quad,\quad}l}
+1 &  \mbox{if}~ \chi(F) = \chi(G) = 1\\ 
-1 &  \mbox{otherwise}\,,
\end{array}
\right.$ &&\\
\end{tabular}\\
\vspace{0.2cm}
\begin{tabular}{llll}
(v) &$\chi(\vee\,F G)=\left\{
\begin{array}{r@{\quad,\quad}l}
+1 &  \mbox{if}~ \chi(F) = 1 ~\mbox{or}~\chi(G) = 1\\ 
-1 &  \mbox{otherwise}\,,
\end{array}
\right.$ &&\\\vspace{0.2cm}
(vi) &$\chi(Q\,x\,F)=\chi(F)\,.$ &&\\
\end{tabular}\\
In this way a sign is defined for all formulas of the mathematical system.

6.) Let $F$ be an R-axiom in $B_S$ with the i-ary R-conclusion 
$p\,\lambda_1,...,\lambda_i$. If there is an R-premise $F'$ of $F$ such that $\chi(F')=-1$,
then we obtain immediately that $\chi(F)=1$. Now we suppose that $\chi(F')=1$ for all 
R-premises $F'$ of $F$. If $F$ contains an i-ary R-premise $p\,\lambda_1',...,\lambda_i'$, 
then we obtain again that $\chi(F)=1$. 
Otherwise it can be shown by induction with respect to $k\geq 0$ that
the R-axiom $F$ is contained in all sets $B_S^{(k)}$, 
especially in $B_S^{(k_0)}$, and thus
$\chi(p\,\lambda_1,...,\lambda_i)=1$ since $p$ occurs i-ary in $F$. 
Therefore we obtain also in this case that $\chi(F)=1$. Note that 
$\chi(F)=1$ for all $F$ in $B_S$ with an equation as an R-conclusion.
Therefore $\chi(F)=1$ for all $F$ in $B_S$.

7.) Next we suppose that $[\Lambda]$ is a proof in $M$ and show that $\chi(F)=1$ 
for all $F \in [\Lambda]$. Then it is clear due to $\chi(\neg\,F)=-\chi(F)$ that
$[\Lambda]$ cannot contain a formula $F$ as well as its negation $\neg\,F$.
Now we employ induction with respect to the rules of inference.

The desired statement is true for the empty proof $[\Lambda]=[\,]$\,.
Assume that $\chi(F)=1$ for all steps $F$ of a proof $[\Lambda]$ in $M$.
For any axiom $F$ we obtain $\chi(F)=1$,
which can be seen very easily by using 4.),
5.), 6.) and (3.9)-(3.11). The induction steps with respect to Rules (b)-(d) 
are also straightforward. Thus we will assume that all the conditions 
for the application of Rule (e) given there are satisfied in $[\Lambda]$. 
Moreover we assume that
$\chi(p\,x_1,...,x_i)=1$, because otherwise it is clear that
$\chi(\to\,p\,x_1,...,x_i\,G)=1$. It remains to show
$\chi(G)=1$. 

But $\chi(p\,x_1,...,x_i)=1$ means that
$p$ occurs i-ary in $B_S^{(k_0)}$, and we conclude due to ${B_S^{(k_0)}}'=B_S^{(k_0)}$
that there is an R-formula $H \in B_S^{(k_0)}$ with an i-ary
R-conclusion $p\,\lambda_1,...,\lambda_i$\,. From the definition of $B_S^{(0)}$
and from $B_S^{(0)} \supseteq B_S^{(k_0)}$ we obtain that $p$ does not occur i-ary
in the R-premises of $H$, and from $H \in B_S^{(k_0)}$ we obtain that
all the R-premises of $H$ have a positive sign. Therefore $H'$, which
is a step in $[\Lambda]$ with $\chi(H')=1$ due to the induction assumption, 
has only positive premises and the j-ary conclusion 
$G\,\frac{\lambda_1}{x_1}...\frac{\lambda_i}{x_i}$.
This is only possible if 
\[
\chi(G\,\frac{\lambda_1}{x_1}...\frac{\lambda_i}{x_i})=
\chi(G)=1\,.
\]

Thus we have proved Proposition (5.1). \dokend \\

As a further result we have shown that all provable formulas $F$ of a mathematical
system $M$ with $B_M = B_S$ satisfie $\chi(F)=1$.

In the following we consider the Peano arithmetic $PA = [M';{\cal L}]$
introduced in example 4 in Section 3. Recall the mathematical system $M'$,
the set ${\cal L}$ of numeral terms and the Induction scheme (IS) defined there.
Since the sixth axiom $\forall x  ~\neg \sim s(x),0$ of PA
has a negative sign, Proposition (5.1) is not sufficient in order
to establish the consistency of PA. In the following we will look
for a more general criterion which guarantees the consistency of PA
and of some other kind of mathematical systems.

Before we proceed with a special Lemma, we first start with
a general defi\-nition for a mathematical system $M=[S;A_M;P_M;B_M]$
and for a fixed predicate symbol $p \in P_M$.

Let $F$ be any formula in $M$ and $x_1$,...,$x_n$ with $n \geq 0$
the uniquely determined sequence of the distinct free variables 
in the formula $F$, ordered according to their first occurrence in $F$.
We define $\Gamma_p(F)=~\to~p\,x_1~...~\to~p\,x_n$
for the block of $p$-premises with respect to all free variables
occurring in $F$. For $n=0$ the string $\Gamma_p(F)$ is defined to be empty.

{\bf (5.2) Lemma concerning relative quantification} 

We consider the mathematical system PA and define a
second mathematical system PA$_{N_0}$ which results from PA by the following changes:
We adjoin the single predicate symbol $N_0$ to the empty set $P_{PA}$ of
predicate symbols of PA. The basis axioms of PA$_{N_0}$ consists exactly
on the two formulas $N_0\,0$ and $\to~N_0\,x\,N_0\,s(x)$ with $x \in X$
and on all formulas $\Gamma_{N_0}(F)\,\Psi_{N_0}(F)$, where F is any basis axiom of PA. 
Here $\Psi_{N_0}$ is the following recursively defined map from the set of all 
PA-formulas to the set of formulas in PA$_{N_0}$:

\begin{tabular}{llll}
(a)  &$\Psi_{N_0}(F)= F$\,, 
&$F$ prime formula in PA\,,\\ 
(b)  &$\Psi_{N_0}(\neg\,F) = \neg~\Psi_{N_0}(F)$\,, 
&$F$\, PA-formula\,,\\ 
(c) &$\Psi_{N_0}(J F G) = J\, \Psi_{N_0}(F) \Psi_{N_0}(G)$\,, 
&$F$, $G$\, PA-formulas\,,\\ 
(d)  &$\Psi_{N_0}(\forall x\,F) = \forall x\,\to ~ N_0\,x\,\Psi_{N_0}(F)$\,, 
&$F$\, PA-formula\,,\\
(e)  &$\Psi_{N_0}(\exists x\,F)= \exists x~\,\,\&~\, N_0\,x\,\Psi_{N_0}(F)$\,, 
&$F$\, PA-formula\,.\\ 
\end{tabular}\\

In (c) the symbol $J$ is a member of the set $\{\to;\leftrightarrow;\&;\vee \}$\,
and in (d), (e) let $x \in X$\,. For the system PA$_{N_0}$ 
we will again require the restriction to the set ${\cal L}$ of numeral argument lists. 
Our statements are as follows
\begin{itemize}
\item[(i)] Let $\lambda$ be any numeral term. 
Then $\Gamma_{N_0}(N_0\,\lambda)\,N_0\,\lambda$ is provable in PA$_{N_0}$\,.
\item[(ii)] Let $F$ be any PA-formula, $x \in X$ and $\lambda$ a numeral term.  
Then $\mbox{CF}(F;\lambda;x)$ is true if and only if
$\mbox{CF}(\Psi_{N_0}(F);\lambda;x)$ is true, and in this case
there holds $\Psi_{N_0}(F\,\frac{\lambda}{x})=\Psi_{N_0}(F)\,\frac{\lambda}{x}$.
\item[(iii)] $\Gamma_{N_0}(F)\,\Psi_{N_0}(F) \in \Pi(\mbox{PA}_{N_0})$ 
for all provable PA-formulas $F$. 
\end{itemize}

\underline{Proof:}~
The restriction concerning the numeral terms
for the formulas of PA and for the use of the rules of inference in PA
is essential here.\\
For (i) one has to show first that
\[
\to ~ N_0\,x ~ \to ~ N_0\,y\,~ N_0 +(xy) \,, \qquad
\to ~ N_0\,x ~ \to ~ N_0\,y\,~ N_0 *(xy)
\]
are both provable in $PA_{N_0}$, using the formal induction principle
for $PA_{N_0}$.
From these formulas and the PA$_{N_0}$-axioms $N_0\,0$ and
$\to\,N_0\,x\,N_0\,s(x)$ we can derive that 
$\Gamma_{N_0}(N_0\,\lambda)\,N_0\,\lambda$ 
is provable in PA$_{N_0}$\,.\\
For the proof of (ii) we employ induction with respect to the formula $F$.\\
For the proof of (iii) we employ induction with respect to the rules of inference in $PA$,
using (i) and (ii). \dokend \\

{\bf (5.3) Reduction of the consistency problem for PA}

Let us define the mathematical system $M=[S;A_M;P_M;B_M]$ as follows:

We choose $A_M=A_S=[\,0;\,s\,;\,+\,;\,*\,]$, $P_M=P_S=[\,N_0\,]$ and $B_M=B_S$, 
where the basis R-axioms $B_S$ of the underlying recursive system $S$ are given by

\begin{tabular}{llll}
(1)  & $N_0\,0$&&\\
(2)  & $\to ~ N_0\,x \,~ N_0\,s(x)$&&\\
\end{tabular}\\
\begin{tabular}{llll}
(3)  & $\to ~ N_0\,x ~ \sim +(0x),x$&&\\
(4)  & $\to ~ N_0\,x ~ \to ~ N_0\,y\,~ \sim +(s(x)y),s(+(xy))$&&\\
\end{tabular}\\
\begin{tabular}{llll}
(5)  & $\to ~ N_0\,x ~ \sim *(0x),0$&&\\
(6)  & $\to ~ N_0\,x ~ \to ~ N_0\,y\,~ \sim *(s(x)y),+(*(xy)y)$&&\\
\end{tabular}\\
\begin{tabular}{llll}
(7) & $\to ~ N_0\,x ~ \to ~ N_0\,y ~ \to ~ \sim s(x),s(y)\,~ \sim x,y$\,.&&\\
\end{tabular}

To the mathematical system $M$ we adjoin the single statement

\begin{tabular}{llll}
$(*)$ \qquad && $\forall x \, \to ~ N_0\,x ~\neg \sim s(x),0$&\\
\end{tabular}

in order to define the mathematical system $M_{PA}=[M((*));{\cal L}]$
with argument lists restricted to the numerals ${\cal L}$\,,
where the basis axiom $(*)$ again has a negative sign.
Here $x,y$ denote different variables.

For all $M_{PA}$ formulas $F$ the following 
expression is provable in $M_{PA}$
$$\to ~~ \forall\,x\,\to~ N_0\,x ~~ \& ~ F\frac{0}{x} ~ 
   \to\,F ~\,  F\frac{s(x)}{x} \quad
   \forall\,x\,\to~ N_0\,x ~ F\,,$$  
which states the \underline{Induction Principle} for $M_{PA}$.
It can be shown by using the Induction Rule (e)
in $M_{PA}$. Therefore $M_{PA}$ is at least as strong as the ``$N_0$-relative" 
Peano arithmetic PA$_{N_0}$.

Next we define an extended recursive system $S^*=[A_S;P^*_S;B^*_S]$
with the predicate symbols $P^*_S=[\,N_0\,;\,Contra]$ by adding the
new basis R-axiom 

\begin{tabular}{llll}
(8) & $\to ~ N_0\,x ~ \to ~ \sim s(x),0\,~ Contra$&&\\
\end{tabular}

to the basis R-axioms (1)-(7) of the recursive system $S$.
The list of basis R-axioms (1)-(8) constitutes the list $B^*_S$.
There results a second mathematical system $M^*=[S^*;A_S;P^*_S;B^*_S]$
with $\Pi(M;{\cal L}) \subseteq \Pi(M^*;{\cal L})$.

Now we assume that PA is contradictory. Then $\exists x \, \sim s(x),0$
is provable in PA, and due to Lemma (5.2) we conclude that  
$\exists x \, \& \, N_0\,x \,\sim s(x),0$
is provable in PA$_{N_0}$.
But then the latter statement which contradicts the statement $(*)$ 
is also provable in $M_{PA}$. 
We conclude that in this case $M_{PA}$ is contradictory like PA.

We show as a further consequence of this assumption that
the 0-ary predicate $Contra$ is provable in 
$[M^*;{\cal L}]$. In order to see that this is true we first check that the formula

\begin{tabular}{llll}
(9) & $\to ~~\exists x ~ \& \, N_0\,x \,\sim s(x),0 ~~ Contra$&&
\end{tabular}

is a consequence of axiom (8) and the predicate calculus in 
$[M^*;{\cal L}]$. $M_{PA}$ is equivalent to 
$[M(\neg \,\exists x ~ \& \, N_0\,x \,\sim s(x),0);{\cal L}]$ and contradictory
due to our assumption. Therefore we can apply Corollary (4.4)
in order to conclude that $\exists x ~ \& \, N_0\,x \,\sim s(x),0$ is provable
in $[M;{\cal L}]$. But every proof in $[M;{\cal L}]$ is also a proof 
in $[M^*;{\cal L}]$, and thus we finally obtain
that $Contra$ is provable in $[M^*;{\cal L}]$, despite the fact that $Contra$ is not
R-derivable in $S^*$. 

\underline{Remark:}

Within $[M^*;{\cal L}]$ we can also apply the Induction Rule (e) on (8) 
for the formula $G = \exists z ~ \& \, N_0\,z \,\sim s(z),0$ with a new variable $z \in X$
in order to conclude
that the following formula is provable in $[M^*;{\cal L}]$:

\begin{tabular}{llll}
(10) & $\to ~~ Contra ~~ \exists x ~ \& \, N_0\,x \,\sim s(x),0\,.$&&
\end{tabular}

Combining the formulas (9) and (10) we conclude that 
\begin{align*}
\leftrightarrow ~~ Contra ~~ \exists x ~ \& \, N_0\,x \,\sim s(x),0
\end{align*}
is provable in $[M^*;{\cal L}]$, but this is not needed in the following. \\

Let $[M;{\cal L}]$ with $M=[S;A_M;P_M;B_M]$ be a general mathematical system with
restricted argument lists in ${\cal L}$ and with an underlying recursive
system $S=[A_S;P_S;B_S]$. Now we suppose that 
$$A_M = A_S = [{\bf a_1}\,;\,{\bf a_2}\,;\,...\,;\,{\bf a_k}]\,,$$
define the alphabet
$
\Lambda=[\,a\,;\, v\,; \,'\, ; \,\underline{(}\,;\, \underline{)}\,]
$
and assume without loss of generality that $\Lambda$ 
and the other sets of symbols in $[M;{\cal L}]$ are disjoint.
Using the strings 
$$
a^{(1)}=a'\,,a^{(2)}=a''\,,a^{(3)}=a'''\,,...\,; \quad 
v^{(1)}=v'\,,v^{(2)}=v''\,,v^{(3)}=v'''\,,...
$$
we encode the lists $\lambda \in {\cal L}$ into strings over the alphabet $\Lambda$ as follows:
Let $\tilde{\lambda}$ result from $\lambda$ if we replace each symbol ${\bf a_i}$ in $\lambda$
by $a^{(i)}$, $i=1,...,k$, each variable ${\bf x_j}$ by $v^{(j)}$, $j \in \N$, 
the brackets ``$($" by ``$\underline{(}$" and  ``$)$" by ``$\underline{)}$". We put
$\tilde{{\cal L}}=\{\tilde{\lambda}\,:\,\lambda \in {\cal L}\,\}\,.$
If $\tilde{{\cal L}}$ is recursively enumerable
then we will simply say that ${\cal L}$ is enumerable. In this case
an R-derivation $[\Lambda]$ in $[S;{\cal L}]$ is defined as an R-derivation in $S$
with  the following restrictions: The R-formulas in $[\Lambda]$ and 
the R-formulas $F$, $G$ in (1.11) have only argument lists in ${\cal L}$,
and the use of the Substitution Rule (1.11)(c) is restricted to $\lambda \in {\cal L}$.
Then the R-formulas in $[\Lambda]$ are called R-derivable in $[S;{\cal L}]$.
By $\Pi_R(S;{\cal L})$ we denote the set of all R-derivable R-formulas 
in $[S;{\cal L}]$\,.

We conclude that the consistency of $PA$ and some other formal mathematical systems
of interest is a consequence of the more general\\

{\bf (5.4) Conjecture}

Let $M=[S;A_M;P_M;B_M]$ be a mathematical system with an underlying recursive
system $S=[A_S;P_S;B_S]$ such that $A_M = A_S$, $P_M = P_S$, $B_M = B_S$. 
Suppose that $[M;{\cal L}]$ is a mathematical system with restricted argument lists
in ${\cal L}$ and that ${\cal L}$ is enumerable\,.
Let $p \in P_S$
and $\lambda_1,...,\lambda_i \in {\cal L}$ for $i \geq 0$ be elementary $A_S$-lists. Then

$p\,\lambda_1,...,\lambda_i \in \Pi(M;{\cal L})$ ~ if and only if ~
$p\,\lambda_1,...,\lambda_i \in \Pi_R(S;{\cal L})$\,.

\underline{Remark:}

The acceptance of (5.4) is merely a verification that the axioms and the 
rules of inference (a)-(e) correspond
to correct methods of deduction. 
Though Conjecture (5.4) implies the consistency of the Peano arithmetic PA,
its meaning seems to go beyond this special application.\\

The mathematical system in Conjecture (5.4) is a special case of the so called 
axiomatized mathematical systems which we will define now.\\

{\bf (5.5) Axiomatized mathematical systems}

Now we consider mathematical systems $M = [S; A_M; P_M; B_M]$ 
with the infinite countable alphabets
\begin{enumerate}
\item[(a)]
$A_M = [{\bf a_1}\,;\,{\bf a_2}\,;\,{\bf a_3}\,;\,...\,]$ 
of constants or operation symbols and
\item[(b)]
$P_M = [{\bf p_1}\,;\,{\bf p_2}\,;\,{\bf p_3}\,;\,...\,]$ of predicate symbols.
\end{enumerate}
The underlying recursive system $S = [A_S; P_S; B_S]$ may have the alphabets
$A_S = [{\bf a_1}\,;\,{\bf a_2}\,;\,...\,{\bf a_k}\,]$ and
$P_S = [{\bf p_1}\,;\,{\bf p_2}\,;\,...\,{\bf p_l}\,]$,
which are finite parts of $A_M$ and $P_M$, respectively.
Next we define the alphabet 
\[
A_{17} := [\,a\,; \, v \,;\, p \,;\, \Box \,; \,'\,; \, * \,;
\, \us \, ; \, \uo \, ;\, \uc \, ; \, \uk \, ; \, \ui \, ;
\un \, ; \ue \, ; \ua \, ;\, \uv \, ;
\, \ual \, ;\, \uex \,]
\]
in order to encode the formulas $F$ of $M$ as follows
\begin{enumerate}
\item[(c)]
The symbols of $A_M$ in $F$ are replaced by $a'\,;\,a''\,;\,a'''\,;\,...$, respectively.
\item[(d)]
The symbols of $P_M$ in $F$ are replaced by $p'\,;\,p''\,;\,p'''\,;\,...$, respectively.
\item[(e)]
The variables of $X$ in $F$ are replaced by $v'\,;\,v''\,;\,v'''\,;\,...$, respectively.
\item[(f)]
The symbols of $E = [\, \sim \, ; \, ( \, ;\, ) \, ; \, , \, ; \, \to \, ;
\neg \, ; \leftrightarrow \, ; \& \, ;\, \vee \, ;
\, \forall \, ;\, \exists \,]$ in $F$ are\\ replaced by 
$ \us \, ; \, \uo \, ;\, \uc \, ; \, \uk \, ; \, \ui \, ;
\un \, ; \ue \, ; \ua \, ;\, \uv \, ;\, \ual \, ;\, \uex$\,, respectively.
\end{enumerate}

Let $A^+$ be the set of all finite and nonempty strings with respect to an alphabet $A$. 
Then to every list $\lambda$ and to every formula $F$ in $M$ 
there corresponds exactly one string $\tilde{\lambda} \in A_{17}^+$ and
$\tilde{F} \in A_{17}^+$ respectively, and therefore
we only need the finite alphabet $A_{17}$ of symbols in order to encode
all formulas of the mathematical system $M$,
where we will suppose that the first 17 symbols of $A_M$ in (a)
form the alphabet $A_{17}$, i.e. ${\bf a_1}=a$, ${\bf a_2}=v$, ... ,
${\bf a_{17}}=\uex$.

Recall that the notation for
$\tilde{F}$ is consistent with the corresponding notation introduced in (2.1)
for the encoding of the R-formulas $F$ in a recursive system.

$M$ is called \underline{axiomatized},
if the set ${\tilde B_M}=\{\,\tilde{F}\,|\,F \in B_M\,\}\subseteq A_{17}^+$ 
is recursively enumerable in the sense of definition (1.12)(a). 
The usual requirement that ${\tilde B_M}$ is decidable leads
to a decision procedure for the formal proofs of $M$, but will not be needed
in the following.

If in addition $[M;{\cal L}]$ is a mathematical system with argument lists
restricted to a set ${\cal L}$ such that
${\tilde {\cal L}}=\{\,\tilde{\lambda}\,|\,\lambda \in {\cal L}\,\}\subseteq A_{17}^+$
is recursively enumerable in the sense of definition (1.12)(a), then
$[M;{\cal L}]$ is called an axiomatized mathematical system with restricted argument lists.

Using these definitions, we obtain the following version of 
G\"odel's First Incompleteness Theorem, which is closely related to Theorem (2.6).\\

{\bf (5.6) Theorem}

Let $[M;{\cal L}]$ be an axiomatized mathematical system with restricted argument lists,
where $M=[S;A_M;P_M;B_M]$ is defined as above. Recall that the set\\
${\tilde {\cal L}}=\{\,\tilde{\lambda}\,|\,\lambda \in {\cal L}\,\}\subseteq A_{17}^+$
is recursively enumerable.

\begin{itemize}
\item[(i)]
$
{\tilde \Pi}(M;{\cal L}) := \{ \,\tilde{F}\,|\,F \in \Pi(M;{\cal L}) \, \} 
\subseteq A_{17}^+
$ is recursively enumerable.
\item[(ii)]
We suppose that the first 11 symbols of the alphabet $A_M$ coincide
with the alphabet $A_{11}$ and that ${\cal L} \supset A_{11}^+$. 
Suppose that there is a formula $G$ of $[M;{\cal L}]$ with $\mbox{free}(G)=\{\,x\,\}$ 
such that $G\,\frac{\lambda}{x}$ is provable in $[M;{\cal L}]$ 
for each 1-ary $S_{11}$-theorem $\lambda \in A_{11}^+$ and such that 
$G\,\frac{\lambda}{x}$ is not provable in $[M;{\cal L}]$ for each 1-ary $S_{11}$-statement 
$\lambda \in A_{11}^+$ which is not an $S_{11}$-theorem. 

Then there is a 1-ary $S_{11}$-statement $\lambda \in A_{11}^+$ such that neither 
the statement $G\,\frac{\lambda}{x}$ nor its negation $\neg\,G\,\frac{\lambda}{x}$ 
are provable in $[M;{\cal L}]$.
\end{itemize}

\underline{Proof:}

(i) is merely a consequence of the facts that the $A_{17}$-encoding 
of the axioms of $[M;{\cal L}]$
leads to a recursively enumerable subset of $A_{17}^+$ and that the rules of inference
are constructive. Therefore we can represent all the relations needed for the 
definition of a formal proof and a provable formula given in Section 3 
in a recursive system which uses the alphabet $A_{17}$, extending the strategy 
in Section 2 for the construction of $S_{11}$. 

(ii) We construct a recursive system 
$S'=[A_{17};P_{S'};B_{S'}]$ which depend on $[M;{\cal L}]$ and $G$ 
and has the following properties:
\begin{itemize}
\item[(1)] $S'$ is a conservative extension of the universal recursive 
system $S_{11}$, i.e. all axioms in $B_{S'} \setminus B_S$ have only conclusions
of the form $p\,\lambda_1,...,\lambda_n$ with $p \in P_{S'} \setminus P_S$,
$A_{17}$-lists $\lambda_1,...,\lambda_n$, $n \geq 0$,
and without equations in $B_{S'}$.
\item[(2)] There is a predicate symbol $B_s^{(1)} \in P_{S'}$ such that

\begin{tabular}{llll}
 & $\to ~ RBasis~x  ~ \to ~ P_s\,y,w  ~ \to ~ EL\,z,u ~\, B_s^{(1)}\,xyz$&&\\
\end{tabular}

is the only basis R-axiom of $S'$ which contains this predicate symbol
in its R-conclusion. Here $x,y,w,z,u \in X$ denote distinct variables.
\item[(3)] Due to (i) there is a predicate symbol 
$\Pi_{M;{\cal L}} \in P_{S'}$ such that
$\Pi_{M;{\cal L}} \,\lambda$ is R-derivable in $S'$ if and only if
$\lambda$ represents a provable formula in $[M;{\cal L}]$.  
\item[(4)] There is a predicate symbol $SbF \in P_{S'}$ such that
$SbF \,\alpha,\beta,\gamma, \delta$ is R-derivable in $S'$ if and only if
$\alpha$ represents a formula $F$ in $[M;{\cal L}]$, $\beta$ a list $\lambda \in {\cal L}$,
$\gamma$ a variable $x \in X$ and $\delta$ the formula $F\,\frac{\lambda}{x}$. 
\item[(5)] There is a predicate symbol $G_{11} \in P_{S'}$ such that 
the only basis R-axioms of $S'$ which contain this predicate symbol
in its R-conclusions are given by the axioms (1)-(12) in the proof of Theorem (2.6). 
\item[(6)] There is a predicate symbol $P^- \in P_{S'}$ such that

\begin{tabular}{llll}
 & $\to ~B_s^{(1)}\,y ~\to ~ G_{11}\,y,s ~ \to ~ \Pi_{M;{\cal L}}\,z ~ 
    \to ~ SbF\,\un{\tilde G},s,{\tilde x},z  
     ~\, P^-~y$&&\\
\end{tabular}

is the only basis R-axiom of $S'$ which contains this predicate symbol
in its R-conclusion, where $y,s,z \in X$ denote distinct variables. Here ${\tilde G} \in A_{17}^+$
represents the formula $G$ and ${\tilde x}\in A_{17}^+$ the only free variable
$x$ of $G$. $P^-~\lambda$ is R-derivable in $S'$
if and only if $\lambda$ is a 1-ary $S_{11}$-statement
for which $\neg\,G\,\frac{\lambda}{x}$ is provable in $[M;{\cal L}]$.
\end{itemize}
The set of all 1-ary $S_{11}$-statements $\lambda \in A_{11}^+$ 
for which $P^-~\lambda$ is R-derivable in 
$S'$ may also be denoted by $P^-$. This will not lead to confusions.
Due to our assumptions we first obtain that $[M;{\cal L}]$ is consistent.
Therefore $P^-~\lambda$ is not R-derivable in $S'$ whenever
$\lambda$ is a 1-ary $S_{11}$-theorem, and  
$P^- \subseteq \overline{\Omega}_s^{(1)}$. But due to Theorem (2.6)
the set $\overline{\Omega}_s^{(1)}$ is not recursively enumerable, in contrast to $P^-$.
We conclude that there is a 1-ary $S_{11}$-statement 
$\lambda \in \overline{\Omega}_s^{(1)} \setminus P^-$ for which
neither $G\,\frac{\lambda}{x} \in \Pi(M;{\cal L})$ nor 
$\neg\,G\,\frac{\lambda}{x} \in \Pi(M;{\cal L})$.
\dokend \\

Next we show that it is possible to construct a recursive system $\Sigma_*$
with a 2-ary universal provability predicate $\Pi \, \lambda, \mu$,
where $\lambda$ represents an axiomatized mathematical system $[M;{\cal L}]$
in the sense of definition (5.5) and $\mu=\tilde{F}$ the $A_{17}$-encoding
of any formula $F$ provable in $[M;{\cal L}]$. This construction of $\Pi$
satisfies L\"ob's representation properties and enables the construction
of G\"odel's self referential statement. Therefore the validity
of G\"odel's Second Incompleteness Theorem is guaranteed for all 
axiomatized mathematical systems which are able to simulate R-derivations
in $\Sigma_*$. Next we prepare the construction of $\Sigma_*$, where we make free
use of Church's thesis, which may be eliminated here by giving an explicit
but very long list of basis R-axioms.


\begin{itemize}
\item[(1)] There is a 2-ary r.e. predicate $RB_2 \subseteq [a]^+ \times A_{17}^+$
which assigns to each $\lambda_1 =a^n$, $n \geq 1$, 
exactly one R-basis string $\mu$ such that $RB_2 \,\lambda_1,\mu$. 
Moreover, for every R-basis string $\mu$ one can find an appropriate parameter
$\lambda_1=a^n$ such that $RB_2 \,\lambda_1,\mu$.
Let $RB_2(\lambda_1)=[A_S;P_S;B_S]$ be the recursive system determined
by the R-basis string $\mu$ with $RB_2 \,\lambda_1,\mu$.
We require that $A_S$ is an initial part of $A_M$ in (5.5)(a) and that 
$P_S$ is an initial part of $P_M$ in (5.5)(b).
$RB_2$ can be constructed if we count the R-basis strings in lexicographic order.
\item[(2)] There is a 2-ary r.e. predicate $L_2 \subseteq [a]^+ \times A_{17}^+$
such that for each fixed $\lambda_2 \in [a]^+$ 
there is a set ${\cal L}$ of $A_M$-lists satisfying (3.15) with 
${\tilde {\cal L}} = \{ \mu \in A_{17}^+\,|\,L_2\,\lambda_2,\mu\}$.
Finally, every r.e. set ${\tilde {\cal L}}$ with ${\cal L}$ satisfying (3.15)
is generated in this way by $L_2$ and at least one parameter $\lambda_2 \in [a]^+$.
Let $L_2(\lambda_2)$ be this set of restricted $A_M$-argument lists determined 
by $L_2$ and the parameter $\lambda_2 \in [a]^+$ .
\item[(3)] There is a 3-ary r.e. predicate $L_3 \subseteq ([a]^+)^2 \times A_{17}^+$
such that for each fixed $\lambda_1, \lambda_2 \in [a]^+$ 
there is a set ${\cal L}$ of $A_M$-lists with 
$${\tilde {\cal L}} = \{ \mu \in A_{17}^+\,|\,L_3\,\lambda_1,\lambda_2,\mu\}\,,$$
where ${\cal L}$ is the smallest possible set 
which satisfies (3.15) and contains the set $L_2(\lambda_2)$ 
and the $A_S$-lists with the alphabet $A_S$ of the recursive system 
$RB_2(\lambda_1)$. 
Let $L_3(\lambda_1,\lambda_2)$ be this set of restricted 
$A_M$-argument lists determined by $L_3$ and the parameters 
$\lambda_1, \lambda_2 \in [a]^+$.
\item[(4)] There is a 4-ary r.e. predicate 
$ML_4 \subseteq ([a]^+)^3 \times A_{17}^+$
such that for each fixed $\lambda_1, \lambda_2, \lambda_3 \in [a]^+$ 
there is an axiomatized mathematical system $M=[S;A_M;P_M;B_M]$ defined
in (5.5) with argument lists restricted to ${\cal L}=L_3(\lambda_1,\lambda_2)$ 
such that $S=RB_2(\lambda_1)$ and
$$
{\tilde {B}_M} = \{ \mu \in A_{17}^+\,|\,ML_4\,\lambda_1,\lambda_2,\lambda_3,\mu\}\,.
$$
Moreover, every axiomatized mathematical system $[M;{\cal L}]$,
where ${\cal L}$ contains all $A_S$-lists of the recursive system underlying $M$,
is genera\-ted in this way by $ML_4$ and appropriate parameters 
$\lambda_1,\lambda_2, \lambda_3 \in [a]^+$.
\item[(5)] There is a  4-ary r.e. predicate $N'_4 \subseteq ( [a]^+)^4$
which coincides with a bijective function $N_4 : ( [a]^+)^3 \to [a]^+$
such that there holds for all $\lambda_1, \lambda_2, \lambda_3, \lambda \in [a]^+$
$$N_4(\lambda_1, \lambda_2, \lambda_3)=\lambda \Leftrightarrow
N'_4\,\lambda_1, \lambda_2, \lambda_3, \lambda\,.$$
\end{itemize}

Since $N_4$ is a bijective, recursive function,
there are uniquely determined
recursive functions $N^{-1}_{4,i} : [a]^+ \to [a]^+$ for $i=1,2,3$ such that
$\lambda_i = N^{-1}_{4,i}(\lambda)$ and
$N_4(\lambda_1, \lambda_2, \lambda_3)= \lambda$ for all $\lambda \in [a]^+$.

We conclude that any parameter $\lambda \in [a]^+$ determines 
a mathematical system $[M;{\cal L}]$ due to the
r.e. relations $RB_2, L_2, L_3, ML_4, N'_4$ described in (1)-(5),
where $S = RB_2(N^{-1}_{4,1}(\lambda))$ is the recursive system
underlying $M$. 
In the following we will simply express this fact by saying that the mathematical
system $[M;{\cal L}]$ is determined by a so called basis number $\lambda \in [a]^+$.
Note that in turn $\lambda$ must not be unique.
\begin{itemize}
\item[(6)] There is a 2-ary r.e. predicate $G_{17} \subseteq (A_{17}^+)^2$
such that $G_{17} \,\lambda,\mu$ holds if and only if 
$\mu = \tilde{\lambda}$ due to (5.5) for $\lambda,\mu \in A_{17}^+$.
We require that $G_{17} \,\lambda,\mu$ can be satisfied for all $\lambda \in A_{17}^+$.
\item[(7)] There is a 2-ary r.e. predicate $Form \subseteq [a]^+ \times A_{17}^+$
such that $Form \,\lambda,\mu$ holds if and only if i) the basis number $\lambda$ 
determines the mathematical system $[M;{\cal L}]$ and ii) $\mu = \tilde{F}$
represents a formula $F$ in $[M;{\cal L}]$.
\item[(8)] There is a 2-ary r.e. predicate $\Pi \subseteq [a]^+ \times A_{17}^+$
such that $\Pi \,\lambda,\mu$ holds if and only if i) the basis number $\lambda$ 
determines the mathematical system $[M;{\cal L}]$ and ii) $\mu = \tilde{F}$
represents a formula $F \in \Pi(M;{\cal L})$. 

This property implies that $\Pi$ satisfies 
the so called first L\"ob condition which states that
whenever a formula $F$ is provable in an axiomatized mathematical system $[M;{\cal L}]$
determined by a basis number $\lambda$, then there must hold $\Pi \,\lambda,\tilde{F}$.
\item[(9)] There is a 2-ary r.e. predicate $\Pi RBasis_2 \subseteq [a]^+ \times A_{17}^+$
such that $\Pi RBasis_2 \,\lambda,\mu$ if and only if

i) $\lambda$ is the basis number of a mathematical system $[M;{\cal L}]$
with an underlying recursive system $S = RB_2(N^{-1}_{4,1}(\lambda)) = [A_S;P_S;B_S]$,\\
ii) $\mu$ is the R-basis string of a recursive system 
$\Sigma'=[A_{\Sigma'};P_{\Sigma'};B_{\Sigma'}]$,\\
iii) there holds $A_{\Sigma'}\subseteq A_S$ and $P_{\Sigma'} \subset P_M$ 
with $P_M$ in (5.5)(b).\\
iv) all basis R-axioms in $B_{\Sigma'}$ are provable in the mathematical system 
$[M;{\cal L}]$ described by the basis number $\lambda$. 

These conditions enable the simulation of the recursive system 
$\Sigma'$ within the mathematical system $[M;{\cal L}]$, even if predicates of $\Sigma'$ 
are neither represented in $S$ nor in $[M;{\cal L}]$.
\item[(10)] There is a 3-ary r.e. predicate $Diag \subseteq [a]^+ \times (A_{17}^+)^2$
such that $Diag \,\lambda,\mu,\nu$ if and only if 

i) $\mu = \tilde{F}$ represents a formula $F$ with exactly one free variable $u \in X$ 
in the mathematical system $[M;{\cal L}]$ given by the basis number $\lambda$,\\
ii) $\mu \in {\cal L}$ and
iii) $\nu$ represents the formula $F\,\frac{\mu}{u} = F\,\frac{\tilde{F}}{u}$\,.
\item[(11)] There is a 2-ary r.e. predicate $R \subseteq [a]^+ \times A_{17}^+$
such that $R \,\lambda,\mu$ if and only if there is a string $\nu \in A_{17}^+$
with i) $Diag \,\lambda,\mu,\nu$ and ii) $\Pi\,\lambda,\un \nu$\,.
\end{itemize}

Consider now a recursive system $\Sigma = [A_{\Sigma};P_{\Sigma};B_{\Sigma}]$
which represents the r.e. predicates in (1)-(11). We require that $A_{\Sigma} \supseteq A_{17}$
is an initial part of $A_M$ given in (5.5)(a).

We will suppose that the names
of the r.e. predicates in (1)-(11) represented in $\Sigma$ are given 
by the corresponding predicates symbols in $P_{\Sigma}$ and that 
$P_{\Sigma}$ is an initial part of $P_M$ in (5.5)(b).
For simplicity we will suppose that
there is no equation involved in $B_{\Sigma}$. We can also choose the basis R-axioms
of $\Sigma$ in such a way that the predicates represented in $\Sigma$ will not change
if the alphabet $A_{\Sigma}$ will be enlarged by using appropriate relatives
representing $A_{\Sigma}$-lists.

Let $\tilde{\Sigma}$ be the R-basis string corresponding to $\Sigma$
and $\tilde{\Pi}$ the encoding of the predicate symbol $\Pi$ according to (5.5).
Now we extend $\Sigma$ to a new recursive system 
$\Sigma_*=[A_{\Sigma};P_{\Sigma};B_{\Sigma_*}]$ by appending the following
four basis R-axioms to the list $B_{\Sigma}$, which are written down in column form

\begin{tabular}{llll}
(2nd and 3rd L\"ob condition)  
& $\to ~  Form\,x,y$&&\\
& $\to ~  Form\,x,z$&&\\
& $\to ~ \Pi\,x,\ui yz$&&\\  
& $\to ~ \Pi\,x,y  $&&\\
& $ \quad~~ \Pi\,x,z$\,,&&\\
\end{tabular}
\begin{tabular}{llll}
& $\to ~  \Pi RBasis_2\,x,\tilde{\Sigma}$&&\\
& $\to ~ G_{17}\,x,s  $&&\\
& $\to ~ G_{17}\,y,t  $&&\\
& $\to ~ \Pi\,x,y$&&\\ 
& $ \quad~~ \Pi\,x,\tilde{\Pi}\,s \uk t$\,.&&\\
\end{tabular}

\begin{tabular}{llll}
(two R-axioms for a self-
& $\to ~ Diag\,x,y,z  $&&\\
referential statement) & $\to ~ \Pi\,x,\un z$&&\\
& $ \quad~~ R\,x,y$\,,&&\\
\end{tabular}
\begin{tabular}{llll}
& $\to ~ Diag\,x,y,z  $&&\\
& $\to ~ R\,x,y$&&\\
& $ \quad~~ \Pi\,x,\un z$\,.&&\\
\end{tabular}

Here $x,y,z,s,t \in X$ denote distinct variables.
The added R-axioms above are in accordance with the meaning of the r.e. predicates
described in (1)-(11). Hence $\Sigma$ and $\Sigma_*$
represent exactly the same predicates. \\

{\bf (5.7) Theorem, due to G\"odel's Second Incompleteness Theorem }
\begin{itemize}
\item[(a)]
Suppose that $\lambda \in [a]^+$, $\mu \in A_{17}^+$ 
and consider the recursive system $\Sigma_*$. 
Then there holds $\Pi \,\lambda,\mu \in \Pi_R(\Sigma_*)$ if and only if 
there is a formula $F$ in the mathe\-matical system $[M;{\cal L}]$
determined by the basis number $\lambda$ such that $\mu = \tilde{F}$
and $F \in \Pi(M;{\cal L})$. 
\item[(b)] Suppose that $\lambda \in [a]^+$ 
satisfies the condition $\Pi RBasis_2 \,\lambda,{\tilde{\Sigma}_*}$
and determines its axiomatized mathe\-matical system $[M;{\cal L}]$.
Then the condition $\Pi RBasis_2 \,\lambda,{\tilde{\Sigma}}$ is also satisfied,
and $[M;{\cal L}]$ is able to simulate the R-derivations in $\Sigma$ and $\Sigma_*$.
Let $F_0$ be any refutable statement in $[M;{\cal L}]$, 
for example the statement $F_0 = \, \neg \forall \, x \sim x,x$, 
where $x = {\bf x_1} \in X$. Define the statement
$$
C = \neg \, \Pi\,\lambda,{\tilde F_0} = \neg \, \Pi\,\lambda,\un \ual v' \us v' \uk v'\,.
$$
Then the statement
$$
\to ~ C ~ \neg \, \Pi\,\lambda,\tilde{C}
$$
is provable in $[M;{\cal L}]$. Moreover, if $C \in \Pi(M;{\cal L})$,
then $[M;{\cal L}]$ is contradictory.
\item[(c)] Let $\Lambda^-$ be the set of all basis numbers $\lambda$
such that the corresponding mathematical system $[M;{\cal L}]$ is contradictory,
and $\Lambda^+ = [a]^+ \setminus \Lambda^-$ the set of all basis numbers
which describes the consistent mathematical systems. Then $\Lambda^-$ is recursively
enumerable, but not $\Lambda^+$. 
\end{itemize}
\underline{Remarks:} 
\begin{itemize}
\item
Part (a) states that the recursive systems
$\Sigma$ and $\Sigma_*$ both represent the same 2-ary predicate $\Pi$
described in (8). 
\item
If the mathematical system $[M;{\cal L}]$
in part (b) also represents the predi\-cate $\Pi$ in (8),
then $C$ states that $[M;{\cal L}]$ 
is free from contradictions, but in this case we cannot prove in $[M;{\cal L}]$ 
the formula $C$ expressing the consistency of this mathematical system.
\item
The presentation and proof of this Theorem are completely independent
on Theorem (5.6) and Theorem (2.6).
\end{itemize}
\underline{Proof:}
\begin{itemize}
\item[(a)]
This is clear since we have already noted that the recursive systems 
$\Sigma$ and $\Sigma_*$ represent the same predicates.
\item[(b)] Since $\Sigma_*$ is an extension of $\Sigma$, we first note that 
$\lambda \in [a]^+$ satisfies
$\Pi RBasis_2\,\lambda, \tilde{\Sigma}$, and therefore the mathematical system
$[M;{\cal L}]$ determined by $\lambda$ is able to simulate $\Sigma$ as well as 
$\Sigma_*$ in the sense that any R-derivation in these systems is also a proof 
in $[M;{\cal L}]$. This will be used in the sequel, where $\lambda$ and 
$[M;{\cal L}]$ are fixed.

First we define the function $g_{17}$, which assigns to each formula $F$ 
of an axiomatized mathematical system described in (5.5) the 
$A_{17}$-string $\tilde{F} = g_{17}(F)$. Recall that the alphabet
$A_M$ in (5.5)(a) starts with $A_{17}$.

We will also make use of the following fact:

Let $F_1$,...,$F_n$ for $n \geq 2$ formulas in $[M;{\cal L}]$ and
assume that\\ $\to\,F_1\,...\to\,F_{n-1}\,\,\,F_{n}$ is provable in $[M;{\cal L}]$.
Then

$
1)\qquad
\to\,\Pi\,\lambda,\tilde{F}_1\,...\to\,\Pi\,\lambda,\tilde{F}_{n-1}\,
\,\,\Pi\,\lambda,\tilde{F}_n \in \Pi(M;{\cal L})\,.
$

It is sufficient to prove this for $n=2$. 
From $\to\,F_1\,F_2 \in \Pi(M;{\cal L})$ we obtain that
$\Pi\,\lambda,g_{17}(\to\,F_1\,F_2)$ is R-derivable in $\Sigma_*$
and hence prova\-ble in $[M;{\cal L}]$. The same holds for the R-formulas
$Form\,\lambda,g_{17}(F_1)$, $Form\,\lambda,g_{17}(F_2)$, and therefore
we can infer our statement from the second L\"ob condition. 

Next we introduce a new variable $y' \in X$ and the abbreviation
 
$
2)\qquad
\Omega := R\,\lambda, g_{17}(R\,\lambda,y')
$

and put $x=\lambda$, $y = g_{17}(R\,\lambda,y')$
and $z = g_{17}(R\,\lambda,g_{17}(R\,\lambda,y')) = g_{17}(\Omega)$ in the
last two R-axioms of $\Sigma_*$ to conclude

$
3)\qquad
\leftrightarrow \Omega ~ \Pi\,\lambda,g_{17}(\neg \Omega) ~ \in ~ \Pi(M;{\cal L})\,.
$

Therefore we obtain from 1)

$
4)\qquad
\to \Pi\,\lambda,g_{17}(\Pi\,\lambda,g_{17}(\neg \Omega))
~ \Pi\,\lambda,g_{17}(\Omega)~ \in ~ \Pi(M;{\cal L})\,.
$

We can also apply the third L\"ob condition to infer

$
5)\qquad
\to \Pi\,\lambda,g_{17}(\neg \Omega)
~ \Pi\,\lambda,g_{17}(\Pi\,\lambda,g_{17}(\neg \Omega))~ \in ~ \Pi(M;{\cal L})\,.
$

Using 3), 4) and 5) we conclude

$
6)\qquad
\to \Omega
~ \Pi\,\lambda,g_{17}(\Omega)~ \in ~ \Pi(M;{\cal L})\,.
$

Since $\to \Omega \to \neg \Omega ~ F_0$ with the refutable formula $F_0$
is an axiom of the propositional calculus, we obtain from 1) with 
$F_1 = \Omega$, $F_2 = \neg \Omega$, $F_3 = F_0$ that

$
7)\qquad
\to \Pi\,\lambda,g_{17}(\Omega)
\to \Pi\,\lambda,g_{17}(\neg \Omega) ~
\Pi\,\lambda,g_{17}(F_0) 
~\in ~ \Pi(M;{\cal L})\,.
$

The propositional calculus yields, 
if applied on 3), 6) and 7)

$
8)\qquad
\to \Omega ~
\Pi\,\lambda,g_{17}(F_0) 
~\in ~ \Pi(M;{\cal L})\,.
$

Since $\to F_0 ~\neg \Omega$ is provable in $[M;{\cal L}]$, 
we obtain from 1) with $F_1 = F_0$, $F_2 = \neg \Omega$ and 3) that

$
9)\qquad
\to \Pi\,\lambda,g_{17}(F_0) ~ \Omega ~
~\in ~ \Pi(M;{\cal L})\,.
$

Using $C = \neg \, \Pi\,\lambda,{\tilde F_0}$, we may rewrite 8) as

$
10)\qquad
\to C ~ \neg \Omega ~
~\in ~ \Pi(M;{\cal L})\,,
$

and applying 1) on 10) with $F_1 = C$, $F_2 = \neg \Omega$
regarding 3) leads to

$
11)\qquad
\to \Pi\,\lambda,g_{17}(C) ~ \Omega ~
~\in ~ \Pi(M;{\cal L})\,.
$

From 10) and 11) we finally obtain the desired result

$
12)\qquad
\to C ~ \neg \Pi\,\lambda,g_{17}(C) 
~\in ~ \Pi(M;{\cal L})\,.
$

Assume that $C \in \Pi(M;{\cal L})$. Then $\Pi\,\lambda,g_{17}(C) \in \Pi(M;{\cal L})$
and 12) would cause a contradiction in  $[M;{\cal L}]$.
\item[(c)] That $\Lambda^-$ is r.e. can be seen by adding with $x \in X$ the R-axiom
$$
\to \Pi\,x,\tilde{F}_0 ~\, \Lambda^-\,x
$$ 
with a refutable formula $F_0$ and the new predicate symbol $\Lambda^-$ to $\Sigma$.

Assume now that $\Lambda^+$ is also r.e., and let 
$S = [A_{S};P_{S};B_{S}]$ be any recursive system which represents
$\Lambda^+$ and all predicates of $\Sigma_*$ and which extends $\Sigma_*$ without using equations such that $A_{S} \supseteq A_{\Sigma}$ and $P_{S} \supseteq P_{\Sigma}$, 
$B_{S} \supseteq B_{\Sigma_*}$. Recall that we have chosen the basis R-axioms
of $\Sigma$ and $\Sigma_*$ in such a way that the predicates represented in 
these systems will be unchanged by extending the set of symbols $A_{\Sigma}$ to $A_{S}$.
We require that $A_S$ is an initial part of $A_M$ in (5.5)(a) and that
$P_S$ is an initial part of $P_M$ in (5.5)(b).

Consider the mathematical system $M_0 = [S_0;A_{M};P_{M};B_{S}]$ with\\
$S_0 = [A_{S};[\,];[\,]\,]$, and adjoin the single axiom
$$
{\cal A} = \forall x ~\to \Lambda^+\,x \to \Pi RBasis_2\,x,\tilde{\Sigma}_* ~
\neg \Pi\,x,\tilde{F}_0
$$
to obtain the new system $M_0({\cal A})$. Note that we have supressed 
the use of the Induction Rule (e) in $M_0$ and $M_0({\cal A})$ 
due to our choice of $S_0$.
Let ${\cal L}_0$ be the set of all $A_S$-lists 
and assume that $[M_0({\cal A});{\cal L}_0]$ is free from contradictions. 
Due to the construction we can find a basis number $\lambda_0$ 
generating $[M_0({\cal A});{\cal L}_0]$ such that $\Lambda^+\,\lambda_0$
and $\Pi RBasis_2\,\lambda_0,\tilde{\Sigma}_*$ are both satisfied.
Therefore $\neg \Pi\,\lambda_0,\tilde{F}_0 \in \Pi(M_0({\cal A});{\cal L}_0)$
can be infered from ${\cal A}$, which contradicts the part (b) of this Theorem.

We conclude that $[M_0({\cal A});{\cal L}_0]$ is contradictory, 
and due to the Deduction Theorem the formula 
$$
\exists x \, \& ~ \Lambda^+\,x ~ \& ~
\Pi RBasis_2\,x,\tilde{\Sigma}_* ~ \Pi\,x,\tilde{F}_0
$$
is provable in $[M_0;{\cal L}_0]$. Since the basis axioms of $[M_0;{\cal L}_0]$ consist only
on the quantifier free positive horn formulas in $B_S$, we obtain with 
a slight modi\-fi\-ca\-tion of Herbrand's Theorem adapted
for use of argument lists that 
$\Lambda^+\,\mu_0$ and $\Pi\,\mu_0,\tilde{F}_0$ 
are both R-derivable in $[S;{\cal L}_0]$ and hence in $S$
for some appropriate $\mu_0 \in [a]^+$,
which is again a contradiction. We conclude that  
$\Lambda^+$ is not r.e. \dokend
\end{itemize}

\section{Outlook}
We have obtained a unified treatment for the generation of languages 
in recursive systems closely related to formal grammars and for 
the predicate calculus in combination with a constructive induction principle.
Thus we hope that this paper may lead to a discussion and further development
of the methods for applications in mathematical logic and computer science.

Complexity results like Theorem (2.10) for certain recursive systems 
and the characterization of special recursive predicates, for example 
by using formal grammars, require an own study which may be of interest 
in computer science.

Special topics of linguistics include the study of a language 
by using formal grammars and languages, see Chomsky \cite{Ch2}, 
Haegeman \& Gueron \cite{HG},
Meyer \cite{My} and Montague \cite{Mn1}, \cite{Mn2}. 
The use of recursive systems may lead to an alternative approach.
 
A further study is necessary to investigate 
additional interesting examples of formal mathematical systems 
which are consistent as a consequence of Conjecture (5.4) 
and to look for a constructive proof of this conjecture. 
Such a study will be related to results given by Gentzen in 
\cite{Gn3}, \cite{Gn4} for the consistency of PA.
But it may also lead to some kind of generalized Herbrand Theorem
in the mathematical systems which are using the Induction Rule.
This generalized Herbrand Theorem should characterize the formulas 
derivable in a mathematical system $[M;{\cal L}]$ satisfying the assumptions
of Conjecture (5.4), at least under additional restrictions,
for example for a restricted use of the Induction Rule (e).
A study of the classical characterization problem due to Herbrand can be found 
in the textbooks of Shoenfield \cite{Sh} and in Heijenoort's
collection of original papers \cite{Hn}.

Kirby \& Paris \cite{KP}, Paris \cite{Pr2} and  Paris \& Harrington \cite{PH}
have presented examples for simple number-theoretical and 
combinatorial statements which are true
but not provable in PA. These statements do not rely on the encodings 
of the logical syntax used by G\"odel in \cite{Gd1} and \cite{Gd2}
for the construction of his famous undecidable formulas,
see also Simpson \cite{Sm1}, \cite{Sm2} and Simpson \& Sch\"utte \cite{ScSm}.
The construction of interesting undecidable combinatorial statements 
for certain mathematical systems besides PA which are consistent 
as a consequence of Conjecture (5.4) may also be a future task.


\end{document}